\documentclass[a4paper]{article}

\usepackage{amsmath}
\usepackage{amsthm}
\usepackage{amssymb}

\usepackage[latin1]{inputenc}

\usepackage[a4paper=true,%
breaklinks=true,%
colorlinks=true,%
linkcolor={blue},%
citecolor={red},%
pdftitle={Scattering, determinants, hyperfunctions in
  relation to Gamma(1-s)/Gamma(s)},%
pdfauthor={Jean-Francois Burnol},%
pdfkeywords={Hankel transform; Scattering; Fredholm determinants},%
pdfsubject={The idea of co-Poisson},%
pdfstartview=FitH,
]{hyperref}

\def\ladate{February 19, 2006.}
\def\ladatedeux{January 2, 2008.}
\let\etendu\diamond

\newtheorem{theorem}{Theorem}
\newtheorem{lemma}[theorem]{Lemma}

\newtheorem{proposition}[theorem]{Proposition}

\theoremstyle{definition}

 \newcommand{\NN}{{\mathbb N}}
 \newcommand{\ZZ}{{\mathbb Z}}
 \newcommand{\QQ}{{\mathbb Q}}
 \newcommand{\RR}{{\mathbb R}}
 \newcommand{\CC}{{\mathbb C}}
 \newcommand{\HH}{{\mathbb H}}

\newcommand{\cA}{{\mathcal A}}
\newcommand{\cB}{{\mathcal B}}
\newcommand{\cC}{{\mathcal C}}
\newcommand{\cD}{{\mathcal D}}
\newcommand{\cE}{{\mathcal E}}
\newcommand{\cF}{{\mathcal F}}
\newcommand{\cH}{{\mathcal H}}
\newcommand{\cX}{{\mathcal X}}
\newcommand{\cY}{{\mathcal Y}}

\newcommand{\Un}{{\mathbf 1}}

\DeclareMathOperator{\Res}{Res}

\let\wh=\widehat
\let\wt=\widetilde

\begin{document}

\title{Scattering,
determinants, hyperfunctions in 
relation to $\frac{\Gamma(1-s)}{\Gamma(s)}$}

\author{Jean-Fran\c cois Burnol}

\date{\ }

\maketitle
\medskip

\begin{abstract}  The method of realizing certain
self-reciprocal transforms as (absolute) scattering,
previously presented in summarized form in the case of the
Fourier cosine and sine transforms, is here applied to the
self-reciprocal transform $f(y)\mapsto {\mathcal H}(f)(x) =
\int_0^\infty J_0(2\sqrt{xy})f(y)\,dy$, which is
isometrically equivalent to the Hankel transform of order
zero and is related to the functional equations of
the Dedekind zeta functions of imaginary quadratic
fields. This also allows to re-prove and to extend 
theorems of de~Branges and V.~Rovnyak regarding square
integrable functions which are self-or-skew reciprocal under
the Hankel transform of order zero. Related integral
formulae involving various Bessel functions are all
established internally to the method. Fredholm determinants
of the kernel $J_0(2\sqrt{xy})$ restricted to finite
intervals $(0,a)$ give the coefficients of first and second
order differential equations whose associated scattering is
(isometrically) the self-reciprocal transform $\mathcal H$,
closely related to the function
$\frac{\Gamma(1-s)}{\Gamma(s)}$. Remarkable distributions
involved in this analysis are seen to have most natural
expressions as (difference of) boundary values
(\textit{i.e.} hyperfunctions.) The present work is
completely independent from the previous study by the author
on the same transform ${\mathcal H}$, which centered around
the Klein-Gordon equation and relativistic causality. In an
appendix, we make a simple-minded observation regarding the
resolvent of the Dirichlet kernel as a Hilbert space
reproducing kernel.
\end{abstract}

\medskip
\begin{footnotesize}
\begin{quote}
keywords:
  Hankel transform; Scattering; Fredholm determinants.

\medskip

Universit\'e Lille 1\\ 
UFR de Math\'ematiques\\ 
Cit\'e scientifique M2\\ 
F-59655 Villeneuve d'Ascq\\ 
France\\
burnol@math.univ-lille1.fr\\
\bigskip

\ladate \\  \ladatedeux: there was an $i$ in \eqref{eq:57beta} and other
equations leading to Theorem \ref{thm:expansion} which should not have been
there. This only affected equation \eqref{eq:Hoperator}.  This is all folks. I
could possibly have proposed other improvements if only the referee who
kept my paper hostage for most of 2006 and 2007 had actually read it. Not
reading it did not prevent from commenting upon it, unfortunately the
substance of those inspired five lines is hard to transfer beneficially to my
wide readership.
\end{quote}
\end{footnotesize}

\pagebreak

\begin{small}

\tableofcontents

\end{small}

\setlength{\normalbaselineskip}{16pt}
\baselineskip\normalbaselineskip
\setlength{\parskip}{6pt}

\section{Introduction (the idea of co-Poisson)}

We explain the underlying framework and the general contours
of this work. Throughout the paper, we have tried to
formulate the theorems in such a form that one can, for most
of them, read their statements without having studied the
preceeding material in its entirety, so a sufficiently
clear idea of the results and methods is easily
accessible. Setting up here all notations and necessary
preliminaries for stating the results would have taken up
too much space.

The Riemann zeta function $\zeta(s)= \frac1{1^s} +
\frac1{2^s} + \frac1{3^s}
+ \dots$ is a meromorphic function in the entire complex
plane with a simple pole at $s=1$, residue $1$. Its
functional equation is usually written in one of the
following two forms:
  \begin{subequations}
\begin{align}
    \pi^{-\frac s2}\Gamma(\frac s2)\zeta(s) &=
    \pi^{-\frac{1-s}2}\Gamma(\frac{1-s}2)\zeta(1-s)\\
\zeta(s) &= \chi_0(s)\zeta(1-s)\qquad \chi_0(s) =
    \pi^{s-\frac12}\frac{\Gamma(\frac{1-s}2)}{\Gamma(\frac
    s2)} \label{eq:0}
\end{align}
  \end{subequations}
The former is related to the expression of $\pi^{-\frac
s2}\Gamma(\frac s2)\zeta(s)$ as a left Mellin
transform\footnote{in the left Mellin transform we use
$s-1$, in the right Mellin transform we use $-s$.}
 and to the Jacobi
identity:
  \begin{subequations}
\begin{align}
    \pi^{-\frac
s2}\Gamma(\frac s2)\zeta(s) &= \frac12\int_0^\infty
(\theta(t) - 1) t^{\frac s2 -1}\,dt &\qquad(\Re(s)>1)\\
&= \frac12\int_0^\infty (\theta(t) - 1 - \frac1{\sqrt t})
t^{\frac s2 -1}\,dt &\qquad(0<\Re(s)<1)\\
  \theta(t)
  &= 1 + 2 \sum_{n\geq1} q^{n^2} \qquad q = e^{-\pi t}\qquad
\theta(t)= \frac1{\sqrt t}\theta(\frac1t)&
\end{align}
  \end{subequations}
The latter form of the functional equation is related to the
expression of $\zeta(s)$ as the right Mellin transform of a
tempered distribution with support in $[0,+\infty)$, which
is self-reciprocal under the Fourier cosine
transform:\footnote{of course, $\delta_m(x) = \delta(x-m)$.}
\begin{subequations}
  \begin{gather}
    \zeta(s) = \int_0^\infty (\sum_{m\geq1} \delta_m(x) - 1)
    x^{-s} dx\\
\int_0^\infty 2\cos(2\pi xy) (\sum_{n\geq1} \delta_n(y) - 1)\,dy
    =  \sum_{m\geq1} \delta_m(x) - 1\qquad (x>0)
  \end{gather}
\end{subequations}
This last identity may be written in the more familiar form:
\begin{equation}\label{eq:05}
  \int_\RR e^{2\pi i xy} \sum_{n\in\ZZ} \delta_n(y) dy =
  \sum_{m\in\ZZ} \delta_m(x) 
\end{equation}
which expresses the invariance of the ``Dirac comb''
distribution $\sum_{m\in\ZZ} \delta_m(x)$ under the Fourier
transform. As a linear functional on Schwartz functions
$\phi$ , the invariance of $\sum_{m\in\ZZ} \delta_m(x)$
under Fourier is expressed as the
Poisson summation formula:
\begin{equation} \label{eq:1}
  \sum_{n\in\ZZ} \wt \phi(n) = \sum_{m\in\ZZ} \phi(m)
\qquad \qquad \wt \phi(y) =  \int_\RR e^{2\pi i xy}\phi(x)\,dx
\end{equation}
The Jacobi identity is the special instance with $\phi(x) = \exp(-\pi t
x^2)$, and conversely the validity of \eqref{eq:1} for
Schwartz functions (and more) may be seen as a corollary to
the Jacobi identity. 

The \emph{idea of co-Poisson} \cite{hab} leads to
\emph{another} formulation of the functional equation as an
identity
  involving functions. The co-Poisson identity
  (\eqref{eq:3} below) appeared in the work of Duffin and
Weinberger \cite{duffwein}. In one of the
approaches to this identity, we start with a function $g$ on
the positive half-line such that both $\int_0^\infty
g(t)\,dt$ and $\int_0^\infty g(t)t^{-1}\,dt$ are absolutely
convergent. Then we consider the averaged distribution
$g*D(x) = \int_0^\infty g(t)D(\frac xt)\,\frac{dt}t$ where
$D(x) = \sum_{n\geq1} \delta_n(x)
- \Un_{x>0}(x)$. This gives (for $x>0$):
\begin{equation}
  g*D(x) = \sum_{n=1}^\infty \frac{g(x/n)}n - \int_0^\infty \frac{g(1/t)}t\,dt
\end{equation}
If $g$ is smooth with
 support in $[a,A]$, $0<a<A<\infty$, then the co-Poisson sum $g*D$
 has Schwartz decrease at $+\infty$ (easy from applying
 the Poisson formula to $\frac{g(1/t)}t$; \textit{cf.}
 \cite[4.29]{copoisson} for a general statement). The right
 Mellin transform $\wh{g*D}(s)$ is related to the right Mellin
 transform $\wh g(s)$ of $g$ via the
 identity:
 \begin{equation}
   \wh{g*D}(s) = \int_0^\infty (g*D)(x)x^{-s}\,dx = \zeta(s) \int_0^\infty
   g(x) x^{-s}\,dx = \zeta(s) \wh g(s)
 \end{equation}
This is because the right Mellin transform of a
multiplicative convolution is the product of the right
Mellin transforms. The necessary calculus of tempered
distributions needed for this and other statements in this
paragraph is detailed in \cite{copoisson}. The functional
equation in the form of \eqref{eq:0} gives:\footnote{one
observes that $\wh{I(g)}(s) = \wh g(1-s)$.}
\begin{equation}\label{eq:2-1}
  \wh{g*D}(s) = \chi_0(s) \zeta(1-s) \wh g(s) = \chi_0(s)
  \wh{I(g)*D}(1-s)\qquad I(g)(t) = \frac{g(1/t)}t 
\end{equation}
One may reinterpret this in a manner involving the cosine
transform $\cC$ acting on $L^2(0,+\infty;dx)$. The Mellin
transform of a function $f(x)$ in
$L^2(0,\infty;dx)$ is a function $\wh f(s)$ on $\Re(s) = \frac12$  which
is nothing else than the Plancherel Fourier
transform of $e^{\frac12 u}f(e^u)$: $\wh f(\frac12 +
i\gamma) = \int_0^\infty f(x)x^{-\frac12 - i\gamma}\,dx =
\int_{-\infty}^\infty f(e^u)e^{\frac u2} e^{-i\gamma
u}\,du$, $\int_0^\infty |f(x)|^2\,dx = \int_{-\infty}^\infty
|f(e^u)e^{\frac u2}|^2\,du =
\frac1{2\pi}\int_{\Re(s)=\frac12} |\wh f(s)|^2 |ds|$. The
unitary operator $\cC I$ is scale invariant hence it is
diagonalized by the Mellin transform: $\wh{\cC I(f)}(s) =
\chi_0(s)\wh{f}(s)$, $\wh{\cC(f)}(s) =
\chi_0(s)\wh{f}(1-s)$, where $\chi_0(s)$ is obtained for
example using $f(x) =
e^{-\pi x^2}$ and coincides with the chi-function defined in
\eqref{eq:0}.  It has modulus $1$ on the critical line as
$\cC$ is unitary. So \eqref{eq:2-1} says that the \emph{co-Poisson
intertwining identity} holds:
\begin{equation}\label{eq:2}
  \cC(g*D) = I(g)*D
\end{equation}
The
co-Poisson intertwining \eqref{eq:2} or explicitely:
\begin{equation}\label{eq:3}
  \int_0^\infty 2\cos(2\pi xy)\left(\sum_{m=1}^\infty\frac{g(x/m)}m
  - \int_0^\infty \frac{g(1/t)}t\,dt\right)dx = \sum_{n=1}^\infty\frac{g(n/y)}y
  - \int_0^\infty g(t)\,dt\qquad(y>0)
\end{equation}
is, when $g$ is smooth with support in $[a,A]$,
$0<a<A<\infty$, an identity of (even) Schwartz functions. If
$g$ is only supposed to be such that $\int_0^\infty
|g(t)|(1+\frac1t)\,dt<\infty$ then the co-Poisson
intertwining $ \cC(g*D) = I(g)*D$ holds as an identity of
distributions (either considered even or with support in
$[0,\infty)$). Sufficient conditions for pointwise validity
have been established \cite{copoisson}. The general
statement of the intertwining is $\cC(g*E) = I(g)*\cC(E)$
where $E$ is an arbitrary tempered distribution with support
on $[0,\infty)$ (see footnote\footnote{both sides in fact
depend only on $E(x)+E(-x)$ as a distribution on the line,
which may be identically $0$, and this happens exactly when
$E$ is a linear combination of odd derivatives of the delta
function.}) and it is proven directly. The co-Poisson
identity \eqref{eq:3} is another manner, not identical with
the Poisson summation formula, to express the invariance of
$D$ under the cosine transform, or the invariance of the
Dirac comb under the Fourier transform.

If the integrable function $g$
has its support in $[a,A]$, $0<a<A<\infty$, then $g*D$ is
constant in $(0,a)$ and its cosine transform (thanks to the
co-Poisson intertwining) is constant in
$(0,A^{-1})$. Up to a rescaling we may take $A = a^{-1}$,
and then $a<1$ (if a non zero example is wanted.)
Appropriate modifications allow to construct non zero even
Schwartz functions constant in $(-a,a)$ and with Fourier
transform again constant in $(-a,a)$ where $a>0$ is
arbitrary \cite{copoisson}.

Schwartz functions are square-integrable so here we have
made contact with the investigation of de~Branges
\cite{bra64}, V~Rovnyak \cite{rov1} and J. and V.~Rovnyak
\cite{rov2,rov3} of square integrable functions on
$(0,\infty)$ vanishing on $(0,a)$ and with Hankel transform
of order $\nu$ vanishing on $(0,a)$. For $\nu=-\frac12$ the
Hankel transform of order $\nu$ is $f(y)\mapsto
\sqrt\frac2\pi\int_0^\infty \cos(xy)f(y)\,dy$ and up to a
scale change this is the cosine transform considered
above. The co-Poisson idea allows to attach  the
zeta function to, among the spaces defined by de~Branges
\cite{bra64}, the spaces associated with the cosine
transform: it has allowed the definition of some novel
Hilbert spaces \cite{cras2001} of entire functions  in
relation with the Riemann zeta function and
Dirichlet $L$-functions (the co-Poisson idea is in
\cite{hab} on the adeles of an arbitrary algebraic number
field $K$; then, the study of the related Hilbert spaces was
begun for $K=\QQ$. Further results were obtained in
\cite{twosys}.)

The study of the function $\chi_0(s) =
\pi^{s-\frac12}\frac{\Gamma(\frac{1-s}2)}{\Gamma(\frac
s2)}$, of unit modulus on the critical line, is interesting.
We proposed to realize the $\chi_0$ function as a
``scattering matrix''. This is indeed possible and has been
achieved in \cite{cras2003}. The distributions, functions,
and differential equations involved are all related to, or
expressed by, the Fredholm determinants of the finite cosine
transform, which in turn are related to the Fredholm
determinants of the finite Dirichlet kernels
$\frac{\sin(t(x-y))}{\pi(x-y)}$ on $[-1,1]$. The study of
the Dirichlet kernels is a topic with a vast literature. A
minor remark will be made in an appendix.

We mentioned the Riemann zeta function and how it relates to
$\chi_0(s) =
    \pi^{s-\frac12}\frac{\Gamma(\frac{1-s}2)}{\Gamma(\frac
    s2)}$ and to the cosine transform. Let us now briefly
    consider the Dedekind zeta function of the Gaussian
    number field
$\QQ(i)$ and how it relates to $\chi(s) =
\frac{\Gamma(1-s)}{\Gamma(s)}$ and to the $\cH$
transform. The $\cH$ transform is
\begin{equation}
  \label{eq:8}
  \cH(g)(y) = \int_0^\infty J_0(2\sqrt{xy})
  g(x)\,dx\qquad\qquad  J_0(2\sqrt{xy}) = \sum_{n=0}^\infty
  (-1)^n \frac{x^n y^n}{n!^2}
\end{equation}
Up to the unitary transformation
$g(x)=(2x)^{-\frac14}f(\sqrt{2x})$, $\cH(g)(y) =
(2y)^{-\frac14}k(\sqrt{2y})$, it becomes the Hankel
transform of order zero $k(y) = \int_0^\infty
\sqrt{xy}J_0(xy)f(x)\,dx$. It is a self-reciprocal, unitary,
scale reversing operator ($\cH(g(\lambda x))(y) =
\frac1\lambda \cH(g)(\frac y\lambda)$). We shall also extend
  its action to tempered distributions on $\RR$ with support
  in $[0,+\infty)$. At the level of right Mellin transforms
  of elements of $L^2(0,\infty;dx)$
it acts as:
\begin{equation}
  \label{eq:4}
\wh{\cH(g)}(s) = \chi(s)\wh g(1-s)\qquad \chi(s) =
\frac{\Gamma(1-s)}{\Gamma(s)}\qquad \Re(s) = \frac12
\end{equation}
It has $e^{-x}\Un_{x\geq0}(x)$ as one among its
self-reciprocal functions, as is verified directly by series
expansion $\int_0^\infty J_0(2\sqrt{xy})e^{-y}\,dy =
  \sum_{n=0}^\infty \frac{(-1)^n}{n!^2} x^n \int_0^\infty
  y^n e^{-y}\,dy = e^{-x}$. The identity
\begin{equation}
  \label{eq:158}
  \int_0^\infty J_0(2\sqrt t) t^{-s}\,dt = \chi(s) =
  \frac{\Gamma(1-s)}{\Gamma(s)} 
\end{equation}
is equivalent to a special case of well-known formulas of Weber, Sonine
and Schafheitlin \cite[13.24.(1)]{watson}. Here we have an
absolutely convergent integral for $\frac34 < \Re(s) < 1$
and in that range the identity may be proven as in: $e^{-x}
= \int_0^\infty
  J_0(2\sqrt{xy}) e^{-y}\,dy = \int_0^\infty J_0(2\sqrt
  y)\frac1x e^{-\frac yx}\,dy$, $\Gamma(1-s) = \int_0^\infty
  J_0(2\sqrt y)(\int_0^\infty x^{-s-1}e^{-\frac yx}\,dx)\,dy
  = \Gamma(s)\int_0^\infty J_0(2\sqrt y)y^{-s}\,dy$. The
  integral is semi-convergent for $\Re(s)>\frac14$, and of
  course \eqref{eq:158} still holds. In particular on the
  critical line and writing $t=e^{u}$, $s=\frac12 +
  i\gamma$, we obtain the identities of tempered
  distributions $\int_\RR e^{\frac12 u} J_0(2e^{\frac12 u})
  e^{-i\gamma u}\,du = \chi(\frac12+i\gamma)$, $e^{\frac12
  u} J_0(2e^{\frac12 u}) = \frac1{2\pi} \int_\RR
  \chi(\frac12+i\gamma) e^{+i\gamma u}\,du$. 

We have $\zeta_{\QQ(i)}(s) = \frac14
\sum_{(n,m)\neq(0,0)} \frac1{(n^2+m^2)^s} = \frac1{1^s} +
\frac1{2^s} + \frac1{4^s} + \frac2{5^s} + \frac1{8^s} +
\dots = \sum_{n\geq1} \frac{c_n}{n^s}$ and it is a
meromorphic function in the entire complex plane with a
simple pole at $s=1$, residue $\frac\pi4$. Its functional
equation assumes at least two convenient well-known forms:
  \begin{subequations}
\begin{align}
    (\sqrt4)^s (2\pi)^{-s}\Gamma(s)\zeta_{\QQ(i)}(s) &=
    (\sqrt4)^{1-s}
    (2\pi)^{-(1-s)}\Gamma(1-s)\zeta_{\QQ(i)}(1-s) \\
(\frac1\pi)^{s}\zeta_{\QQ(i)}(s) &=
    \chi(s)(\frac1\pi)^{1-s}\zeta_{\QQ(i)}(1-s)\qquad \chi(s) = 
    \frac{\Gamma(1-s)}{\Gamma(s)}  
\end{align}
  \end{subequations}
  \begin{subequations}
The former is related to the expression of
$\pi^{-s}\Gamma(s)\zeta_{\QQ(i)}(s)$ as a left Mellin 
transform:
\begin{align}
    \pi^{-s}\Gamma(s)\zeta_{\QQ(i)}(s) &= \frac14\int_0^\infty
(\theta(t)^2 - 1) t^{s -1}\,dt &\qquad(\Re(s)>1)\\
&= \frac14\int_0^\infty (\theta(t)^2 - 1 - \frac1t)
t^{s -1}\,dt &\qquad(0<\Re(s)<1)\\
 \theta(t)^2
&= \frac1{t}\theta(\frac1t)^2& 
\end{align}
  \end{subequations}
The latter form of the functional equation is related to the
expression of $(\frac1\pi)^{s}\zeta_{\QQ(i)}(s)$ as the right Mellin
transform of a tempered distribution which is supported in
$[0,\infty)$ and which is self-reciprocal under the $\cH$-transform:
\begin{subequations}
  \begin{gather}
   (\frac1\pi)^{s}\zeta_{\QQ(i)}(s)= \int_0^\infty
    (\sum_{m\geq1} c_m\delta_{\pi m}(x) - \frac14) x^{-s}
    dx\\
\int_0^\infty J_0(2\sqrt{xy}) (\sum_{n\geq1} c_n\delta_{\pi
n}(y) - \frac14)\,dy
    = \sum_{m\geq1} c_m\delta_{\pi m}(x) -
    \frac14\Un_{x>0}(x) = E(x) \qquad (x>0)
  \end{gather}
\end{subequations}
The invariance of $E$ under the
$\cH$-transform is equivalent  to
the validity of the functional equation of
$(\frac1\pi)^{s}\zeta_{\QQ(i)}(s)$ and it having a pole with
residue $\frac14$ at $s=1$.
The
co-Poisson intertwining
becomes the assertion:
\begin{equation}
  \label{eq:7}
  y>0\implies \int_0^\infty J_0(2\sqrt{xy})\left(\sum_{m=1}^\infty
  c_m\frac{g(x/\pi m)}{\pi m} - \frac14 \int_0^\infty
  g(\frac1t)\,\frac{dt}t\right)\,dx = \sum_{n=1}^\infty
  c_n\frac{g(\pi n/y)}{y} - \frac14 \int_0^\infty
  g(t)\,dt
\end{equation}
If $g$ is smooth with support in $[b,B]$, $0<b<B<\infty$,
then we have on the right hand side a function of Schwartz
decrease at $+\infty$ (compare to Theorem \ref{theo:2}), and
its $\cH$-transform is also of Schwartz decrease at
$+\infty$.  The former is constant for $0<y<\pi B^{-1}$ and
the latter is constant for $0<x<\pi b$. The supremum of the
values obtainable for the product of the lengths of the
intervals of constancy is $\pi^2$. But, as for the cosine
and sine transforms, there does exist smooth functions
which are constant on a given $(0,a)$ for arbitrary $a>0$
with an $\cH$ transform again constant on $(0,a)$ and have
Schwartz decrease at $+\infty$ (the two constants being
arbitrarily prescribed.)

De~Branges and V.~Rovnyak have obtained \cite{bra64, rov1}
rather complete results in the study of the Hankel transform
of order zero $f(x)\mapsto g(y)=\int_0^\infty \sqrt{xy}
J_0(xy)f(x)\,dx$ from the point of view of understanding the
support property of being zero and with transform again zero
in a given interval $(0,b)$. They obtained an isometric
expansion (Theorem \ref{thm:brarov} of section \ref{sec:2})
and also the detailed description of the related spaces of
entire functions (\cite{bra64}). The more complicated case
of the Hankel transforms of non-zero integer orders was
treated by
J. and V.~Rovnyak \cite{rov2, rov3}. These, rather complete,
results are an indication that the Hankel transform of order
zero or of integer order is easier to understand than the
cosine or sine transforms, and that doing so thoroughly
could be useful to better understand how to try to
understand the cosine and sine transforms.

The kernel $J_0(2\sqrt{uv})$ of the $\cH$-transform
satisfies  the Klein-Gordon equation in
the variables $x= v-u$, $t = v+u$:
\begin{equation}
  \label{eq:5} (\frac{\partial^2}{\partial u\partial v} +
1) J_0(2\sqrt{uv}) = 
  (\square + 1) J_0(2\sqrt{uv})  =
  (\frac{\partial^2}{\partial t^2} -
  \frac{\partial^2}{\partial x^2} + 1) J_0(\sqrt{t^2 - x^2})
  = 0
\end{equation}
It is a noteworthy fact
  that the support condition, initially considered by
  de~Branges and V.~Rovnyak, and  which, nowadays, is also  seen
  to be in relation
with the co-Poisson identities, has turned out to be related
to the relativistic causality governing the propagation of
solutions to the Klein-Gordon equation. This
has been established in \cite{einstein} where we obtained as
an
  application of this idea the isometric expansion of
  \cite{bra64,rov1} in a novel manner. It was furthermore
  proven in
\cite{einstein} that
    the $\cH$ transform is indeed an (absolute) scattering,
    in fact the scattering from the past boundary to the
    future boundary of the Rindler wedge $0<|t|<x$ for
    solutions of a first order, two-component (``Dirac''),
    form of the KG equation.

In the present paper, which is completely independent from
\cite{einstein}, we shall again study the $\cH$-transform
and show in particular how to recover in yet a different way
the earlier results of \cite{bra64,rov1} and also we shall
extend them. This will be based on the methods from
\cite{cras2002,cras2003}, and uses the techniques motivated
by the study of the co-Poisson idea \cite{copoisson}. Our
exposition will thus give a fully detailed account of the
material available in summarized form in
\cite{cras2002,cras2003}. Then we proceed with a development
of these methods to provide the elucidation of the (two
dimensions bigger) spaces of functions constant in
$(0,a)$ and with $\cH$-transforms constant in $(0,a)$.

The use of tempered distributions is an important point of
our approach\footnote{at the bottom of page 456 of
\cite{bra64} the formulas given for $A(a,z)$ and $B(a,z)$ as
completed Mellin transforms are lacking terms which would
correspond to Dirac distributions; possibly related to this,
the isometric expansion as presented in Theorem II of
\cite{bra64} is lacking corresponding terms. The exact
isometric expansion appears in \cite{rov1} and the exact
formulas for $A(a,z)$ and $B(a,z)$ as completed Mellin
transforms appear, in an equivalent form, in
\cite[eq.(37)]{rov3}.}; also one may envision the co-Poisson
idea as asking not to completely identify a distribution
with the linear functional it ``is''. In this regard it is
of note that the distributions which arise following the
method of \cite{cras2002} are seen in the present case of
the study of the $\cH$-transform to have a very natural
formulation as differences of boundary values of analytic
functions, that is, as hyperfunctions \cite{morimoto}. We do
not use the theory of hyperfunctions as such, but could not
see how not to mention that this is what these distributions
seem to be in a natural manner.

The paper contains no number theory. And, the reader will
need no prior knowledge of \cite{Bra}; some familiarity with
the $m$-function of Hermann Weyl \cite{coddlevinson,
levsarg, reedsimon} is necessary at one stage of the
discussion (there is much common ground, in fact, between
the properties of the $m$-function and the axioms of
\cite{Bra}). The reproducing kernel in any space with the
axioms of \cite{Bra} has a specific appearance (equation
\eqref{eq:80} below) which has been used as a guide to what
we should be looking for. The validity of the formula is
re-proven in the specific instance considered
here\footnote{the critical line here plays the rôle of the
real axis in \cite{Bra}, $s$ is $\frac12 - iz$ and the use
of the variable $s$ is most useful in distinguishing the
right Mellin transforms which need to be completed by a
Gamma factor from the left Mellin transforms of
``theta''-like functions.}.  Regarding the differential
equations governing the deformation, with respect to the
parameter $a>0${\hskip3pt}\footnote{the $a$ here corresponds
to $\frac12 a^2$ in \cite{bra64}.}, of the Hilbert spaces,
we depart from the general formalism of \cite{Bra} and
obtain them in a canonical form, as defined in
\cite[\textsection3]{levsarg}.  Interestingly this is
related to the fact that the $A$ and $B$ functions
(connected to the reproducing kernel, equation
\eqref{eq:80}) which are obtained by the method of
\cite{cras2002} turn out not to be normalized according to
the rule in general use in \cite{Bra}. Each rule of
normalization has its own advantages; here the equations are
obtained in the Schrödinger and Dirac forms familiar from
the spectral theory of linear second order differential
equations \cite{coddlevinson, levsarg, reedsimon}. This
allows to make reference to the well-known
Weyl-Stone-Titchmarsh-Kodaira theory \cite{coddlevinson,
levsarg, reedsimon},
and to understand $\cH$ as a scattering. Regarding spaces
with the axioms of \cite{Bra}, the articles of Dym
\cite{dym} and Remling \cite{remling} will be useful to the
reader interested in second order linear differential
equations. And we refer the reader with number theoretical
interests to the recent papers of Lagarias \cite{laghouches,
lagacta}.

The author has been confronted with a dilemma: a substantial
portion of the paper (most of chapters \ref{sec:5},
\ref{sec:6}, \ref{sec:8}) has a general validity for
operators having a kernel of the multiplicative
type $k(xy)$ possessing certain properties in common with
the cosine, sine or $\cH$ transforms. But on the other hand
the (essentially) unique example 
where all quantities
arising may be computed is the $\cH$ transform (and
transforms derived from it, or closely related to it, as the Hankel
transforms of integer orders). We have tried to give proofs
whose generality is obvious, but we also  made full use of
distributions, as this allows to give to the quantities
arising very natural expressions.  Also we never hesitate
using arguments of  analyticity although
for some topics (for example, some aspects involving certain
integral equations and
Fredholm determinants) this is certainly  not really needed.

\section{Hardy spaces and the de~Branges-Rovnyak isometric
  expansion}
\label{sec:2}

Let us state the isometric expansion of \cite{bra64,rov1}
regarding the square integrable Hankel transforms of order
zero. We reformulate the theorem to express it with
the $\cH$ transform \eqref{eq:8} rather than the Hankel
transform of order zero.

  \begin{subequations}

\begin{theorem}[\cite{bra64}, \cite{rov1}]:\label{thm:brarov}
Let $k\in L^2(0,\infty;dx)$. The functions $f_1$
and $g_1$, defined as the following integrals:
\begin{align}
    f_1(y) &= 
  \int_{y}^\infty J_0(2\sqrt{y(x-y)})k(x)\,dx\;,\\
g_1(y) &= k(y) -
  \int_{y}^\infty
  \sqrt{\frac y{x-y}}J_1(2\sqrt{y(x-y)})
  k(x)\,dx\;,
\end{align}
exist in $L^2$ in the sense of mean-square convergence, and
they verify:
\begin{equation}
  \int_0^\infty |f_1(y)|^2 +
  |g_1(y)|^2\,dy = \int_0^\infty |k(x)|^2\,dx\;.
\end{equation}
The function $k$ is given in terms of the pair
$(f_1,g_1)$ as:
\begin{equation}
  k(x)  = g_1(x) + \int_0^{x} 
  J_0(2\sqrt{y(x - y)}) f_1(y) \,dy  - \int_0^{x}
  \sqrt{\frac y{x-y}}{J_1(2\sqrt{y(x 
  - y)})}g_1(y) \,dy
\end{equation}
The assignment $k\mapsto (f_1,g_1)$ is a unitary
equivalence of $L^2(0,\infty;dx)$ with $L^2(0,\infty;dy) \oplus
L^2(0,\infty;dy)$ such that the $\cH$-transform acts as
$(f_1,g_1)\mapsto (g_1,f_1)$. Furthermore $k$ and $\cH(k)$ both
identically 
vanish in $(0,a)$ if and only if $f_1$ and $g_1$ both
identically vanish in $(0,a)$.
\end{theorem}
  \end{subequations}

Let us mention the following (which follows from the proof
we have given of Thm. \ref{thm:brarov} in \cite{einstein}):
\emph{if $f_1$, $f_1'$, $g_1$, $g_1'$ are in $L^2$ then
$k$, $k'$ and $\cH(k)'$ are in $L^2$. Conversely if $k$,
$k'$ and $\cH(k)'$ are in $L^2$ then the integrals defining
$f_1(y)$ and $g_1(y)$ are convergent for each $y>0$ as
improper Riemann integrals, and $f_1'$ and $g_1'$ are in
$L^2$.}

It will prove convenient
to work with  $(f(x),g(x)) =\frac12
(g_1(\frac x2)+f_1(\frac x2), g_1(\frac x2)-f_1(\frac x2))$:
\begin{subequations}
  \begin{align}
    f(y) &=   \frac12 k(\frac{y}2)
 + \frac12 \int_{y/2}^\infty
  \left(J_0(\sqrt{y(2x-y)}) -
  \sqrt{\frac y{2x-y}}J_1(\sqrt{y(2x-y)}) \right)
  k(x)\,dx \label{eq:9a}\\    
g(y) &=  \frac12 k(\frac{y}2)
 - \frac12 \int_{y/2}^\infty
  \left(J_0(\sqrt{y(2x-y)}) +
  \sqrt{\frac y{2x-y}}J_1(\sqrt{y(2x-y)}) \right)
  k(x)\,dx \label{eq:9b}\\
  k(x)  &= f(2x) + \frac12\int_0^{2x} 
  \left(J_0(\sqrt{y(2x - y)}) -
  \sqrt{\frac y{2x-y}}{J_1(\sqrt{y(2x 
  - y)})}\right)f(y) \,dy \notag\\
&+ g(2x) - \frac12\int_0^{2x} 
  \left(J_0(\sqrt{y(2x - y)}) +
  \sqrt{\frac y{2x-y}}J_1(\sqrt{y(2x 
  - y)})\right)g(y) \,dy \label{eq:9c}\\ \label{eq:9d}
\int_0^\infty &|k(x)|^2\,dx = \int_0^\infty |f(y)|^2 +
  |g(y)|^2 \,dy
\end{align}
\end{subequations}
The $\cH$ transform on $k$ acts as $(f,g)\mapsto(f,-g)$. The
pair $(k,\cH(k))$ identically vanishes on $(0,a)$ if and
only if the pair $(f,g)$ identically vanishes on
$(0,2a)$. The structure of the formulas is more apparent
after observing ($x,y>0$):
\begin{equation}
  \label{eq:10}
  \frac\partial{\partial x} \Big(\frac12\, J_0(\sqrt{y(2x -
  y)})\Un_{0<y<2x}(y)\Big) = \delta_{2x}(y) - \frac12\,\sqrt{\frac
  y{2x-y}}J_1(\sqrt{y(2x-y)}) \Un_{0<y<2x}(y)
\end{equation} 
In this section I shall prove the existence of an isometric
expansion $k\leftrightarrow(f,g)$ having the stated support
properties and relation to the $\cH$-transform; that this
construction does give the equations \eqref{eq:9a},
\eqref{eq:9b}, \eqref{eq:9c},  will only be
established in the last section (\ref{sec:9}) of the paper.
The method followed in this section coincides partly with
the one of V.~Rovnyak \cite{rov1}; we try to produce the
most direct arguments, using the commonly known facts on
Hardy spaces. The reader only interested in Theorem
\ref{thm:brarov} is invited after having read the present
section to then jump directly to section
\ref{sec:9} for the conclusion of the proof.

To a function $k\in L^2(0,\infty;dx)$ we associate the
analytic function 
\begin{equation}
  \label{eq:11}
  \wt k(\lambda) = \int_0^\infty e^{i\lambda x}k(x)\,dx\qquad (\Im(\lambda)>0)
\end{equation}
with boundary values for $\lambda\in\RR$ again written
$\wt k(\lambda)$, which defines an element of
$L^2(\RR,\frac{d\lambda}{2\pi})$, the assignment $k\mapsto
\wt k$ being unitary from $L^2(0,\infty;dx)$ onto
$\HH^2(\Im(\lambda>0), \frac{d\lambda}{2\pi})$. Next we have
the conformal equivalence and its associated unitary map
from $\HH^2(\Im(\lambda>0), \frac{d\lambda}{2\pi})$ to
$\HH^2(|w|<1, \frac{d\theta}{2\pi})$:
\begin{equation}
  \label{eq:12}
  w = \frac{\lambda - i}{\lambda + i}\qquad K(w) =
  \frac1{\sqrt2}\frac{\lambda + i}i\; \wt k(\lambda)
\end{equation}
It is well known that this indeed unitarily identifies the
two Hardy spaces. With $k_0(x) = e^{-x}$, $\wt k_0(\lambda)
= \frac{i}{\lambda + i}$, $K_0(w) = \frac1{\sqrt2}$, and $\|
k_0\|^2 = \int_0^\infty e^{-2x} dx = \frac12 =
\|K_0\|^2$. The functions $\wt k_n(\lambda) = (\frac{\lambda
- i}{\lambda + i})^n \frac{i}{\lambda + i}$ correspond to
$K_n(w) = \frac1{\sqrt2} w^n $. To obtain explicitely the
orthogonal basis $(k_n)_{n\geq0}$, we first observe that $w
= 1 - 2\frac{i}{\lambda +
  i}$, so as a unitary operator it acts as:
\begin{equation}
  \label{eq:13}
  w\cdot k(x) = k(x) - 2 \int_0^x e^{-(x-y)}k(y)\,dy = k(x)
  - e^{-x} 2\int_0^x e^y k(y)\,dy
\end{equation}
Writing $k_n(x) = P_n(x) e^{-x}$ we thus obtain $P_{n+1}(x)
= P_n(x) - 2\int_0^x P_n(y)\,dy$:
\begin{equation}
  \label{eq:14}
  P_n(x) = \left(1 - 2\int_0^x\right)^n\cdot 1 = \sum_{j=0}^n
  {n \choose j}\frac{(-2x)^j}{j!} 
\end{equation}
So as is well-known $P_n(x) = L_n^{(0)}(2x)$ (in the
notation of \cite[\textsection5]{szego}) where the Laguerre
polynomials $L_n^{(0)}(x)$ are an  orthonormal
system for the weight $e^{-x}dx$ on $(0,\infty)$.

One of the most common manner to be led to the
  $\cH$-transform is to define it from the two-dimensional
  Fourier transform as:
  \begin{equation}
    \label{eq:15}
\begin{split}    \cH(f)(\frac12 r^2) = \frac1{2\pi}\iint
  e^{i(x_1y_1+x_2y_2)}f(\frac{y_1^2 + y_2^2}2)dy_1dy_2 =
  \int_0^\infty \left(\int_0^{2\pi}
  e^{irs\cos\theta}\frac{d\theta}{2\pi}\right) f(\frac 12
  s^2)sds\\ \cH(f)(\frac12 r^2) = \int_0^\infty J_0(rs)  f(\frac 12
  s^2)\,sds \qquad r^2 = x_1^2 + x_2^2, s^2 = y_1^2 + y_2^2
\end{split}
  \end{equation}
which proves its unitarity, self-adjointness, and
self-reciprocal character and the fact that it has $e^{-x}$
has self-reciprocal function. The direct verification of
  $\cH(k_0) = k_0$ is
immediate: $\cH(k_0)(x) =
  \int_0^\infty J_0(2\sqrt{xy})e^{-y}\,dy =
  \sum_{n=0}^\infty \frac{(-1)^n}{n!^2} x^n \int_0^\infty y^n
  e^{-y}\,dy = e^{-x}$. Then, $\cH(e^{-tx}) = t^{-1}e^{-\frac
  xt}$ for each $t>0$. So $\int_0^\infty
  e^{-tx}\cH(k)(x)\,dx = t^{-1}\int_0^\infty e^{-\frac1t
  x}k(x)\,dx$ hence:
\begin{equation}
  \label{eq:16}
  \forall k \in L^2(0,\infty;dx)\qquad \wt{\cH(k)}(\lambda)
  = \frac i\lambda \wt k(\frac{-1}\lambda)
\end{equation} 
With the notation $\cH(K)$
for the function in $\HH^2(|w|<1)$ corresponding to
$\cH(k)$, we obtain from \eqref{eq:12}, \eqref{eq:16}, an
extremely simple result:\footnote{we also take note of the operator
identity $\cH\cdot w = - w\cdot \cH$.}
\begin{equation}
  \label{eq:17}
  \cH(K)(w) = K(-w)
\end{equation}
This obviously leads us to associate to $K(w) =
\sum_{n=0}^\infty c_n w^n$ the functions:
\begin{subequations}
 \begin{align} \label{eq:18a} 
  F(w) &:= \sum_{n=0}^\infty c_{2n} w^{n}\\
   \label{eq:18b}  
G(w) &:= \sum_{n=0}^\infty c_{2n+1} w^n \\
K(w) &= F(w^2) + w G(w^2)
 \end{align}
\end{subequations}
and to $k$ the functions $f$ and $g$ in $L^2(0,+\infty;dx)$
corresponding to $F$ and $G$. Certainly, $\|k\|^2 = \|f\|^2
+ \|g\|^2$, and the assignment of $(f,g)$ to $k$ is an
isometric identification. Furthermore, certainly the $\cH$
transform acts in this picture as $(f,g)\mapsto (f,-g)$.
Let us now check the support properties. Let $\alpha(m)$ be
the leftmost point of the (essential) support of a given
$m\in L^2(0,\infty;dx)$. As is well-known,
\begin{equation}
  \label{eq:19}
  - \alpha(m) = \limsup_{t\to+\infty} \frac1t\log|\wt m(it)|\;,
\end{equation}
If $w$ corresponds to $\lambda$ via
\eqref{eq:12} then $w^2$ corresponds to $\frac12(\lambda -
\frac1\lambda)$, so if to a function $f$ with corresponding
$F(w)$ we associate the function
$\psi(f)\in  L^2(0,\infty;dx)$ which corresponds to
$F(w^2)$, 
\begin{equation}
  \label{eq:21}
  (t+1)\wt{\psi(f)}(i\,t) = (\frac{t+\frac1t}2 +
  1)\wt{f}(i\,\frac{t+\frac1t}2)\;, 
\end{equation}
then we have the identity:
\begin{equation}
  \label{eq:20}
  \alpha(\psi(f)) = \frac12 \alpha(f)
\end{equation}
Returning to  $F$ (resp. $f$) and $G$
(resp. $g$) associated via \eqref{eq:18a}, \eqref{eq:18b},
to $K$ (resp. $k$) we thus have $k = \psi(f) + w\cdot
\psi(g)$, $\cH(k) = \psi(f) - w\cdot \psi(g)$, hence if the
pair $(f,g)$ vanishes on $(0,2a)$ then the pair $(k,\cH(k))$
vanishes on $(0,a)$ (clearly the unitary operator of
multiplication by $w = \frac{\lambda-i}{\lambda+i}$ does not
affect $\alpha(m)$.) Conversely, as $\alpha(f) = 2 \alpha( k
+ \cH(k))$ and $\alpha(g) = 2 \alpha( k - \cH(k))$, if the
pair $(k,\cH(k))$ vanishes on $(0,a)$ then the pair $(f,g)$
vanishes on $(0,2a)$.

We have thus established the existence of an isometric
expansion, its support properties, and its relation to the
$\cH$-transform. That there is indeed compatibility of
\eqref{eq:9a} and \eqref{eq:9b} with \eqref{eq:18a} and
\eqref{eq:18b}, and with \eqref{eq:9c}, will be established
in the last section (\ref{sec:9}) of the paper with a direct study of
\eqref{eq:21}. In the meantime equations \eqref{eq:9a},
\eqref{eq:9b}, \eqref{eq:9c} and \eqref{eq:9d}
will have been confirmed in another manner. Yet another proof
of the isometric expansion has been given in \cite{einstein}.

\section{Tempered distributions and their $\cH$ and Mellin transforms}

Any distribution $D$ on $\RR$ has a primitive. If the closed
support of $D$ is included in $[0,+\infty)$, then it has a
unique primitive, which we will denote $\int_0^x D(x)\,dx$,
or, more safely, $D^{(-1)}$, which also has its support in
$[0,+\infty)$. The temperedness of such a $D$ is equivalent
to the fact that $D^{(-N)}$ for $N\gg0$ is a continuous
function with polynomial growth. With $D^{(-N)}(x) = (1+
x^2)^M g_{(N,M)}(x)$, $M\gg 0$, we can express $D$ as
$P(x,\frac d{dx})(g)$ where $P$ is a polynomial and $g\in
L^2(0,\infty;dx)$. Conversely any such expression is a
tempered distribution vanishing in $(-\infty,0)$.  The
Fourier transforms of such tempered distributions $\wt
D(\lambda)$ appear thus as the boundary values of $Q(\frac
d{d\lambda},\lambda) f(\lambda)$ where $Q$ are polynomials
and the $f$'s belong to $\HH^2(\Im(\lambda)>0)$. As taking
primitives is allowed we know without further ado that this
class of analytic functions
is the same thing as the space of
functions $g(\lambda) = R(\frac d{d\lambda},\lambda,
\lambda^{-1}) f(\lambda)$, $R$ a polynomial and
$f\in\HH^2$. It is thus clearly left stable by the
operation:
\begin{equation}
  \label{eq:25}
  g\mapsto \cH(g)(\lambda) := \frac{i}\lambda
  g(\frac{-1}\lambda)\qquad(\Im(\lambda)>0) 
\end{equation}
which will serve to define the action of $\cH$ on tempered
distributions with support in $[0,+\infty)$.

Let us also use \eqref{eq:25}, where
now
$\lambda\in\RR$, to define $\cH$ as a unitary operator on
$L^2(-\infty,+\infty;dx)$. It will anti-commute with
$f(x)\to f(-x)$ so:
\begin{equation}
  \label{eq:26}
  \cH(f)(x) = \int_{-\infty}^\infty (
J_0(2\sqrt{xy})\Un_{x>0}(x)\Un_{y>0}(y) -
J_0(2\sqrt{xy})\Un_{x<0}(x)\Un_{y<0}(y))f(y)\,dy
\end{equation}

Useful operator identities  are easily established from \eqref{eq:25}:
\begin{subequations}
\begin{align}
  \label{eq:31a}
  x \frac{d}{dx}\cdot \cH &= - \cH \cdot  \frac{d}{dx}x
&\quad\text{and}&\quad
  \frac{d}{dx}x\cdot \cH = - \cH \cdot  x\frac{d}{dx}\\
  \label{eq:34}
  \frac{d}{dx}\cdot \cH &= \cH\cdot\int_0^x
&\quad\text{and}&\quad
  \int_0^x\cdot\,\cH = \cH\,\cdot \frac{d}{dx}
\end{align}
\begin{align}
\label{eq:31b}
  x\cdot \cH &= - \cH\cdot \frac{d}{dx}\,x\, \frac{d}{dx}
&\quad\text{and}&\quad
  \cH\cdot x = - \frac{d}{dx}\,x\, \frac{d}{dx}\cdot \cH
\end{align}
\end{subequations}
It is important that $\frac{d}{dx}$ is always taken
in the distribution sense. It would actually be possible to
define the action of $\cH$
on distributions supported in $[0,+\infty)$ without mention
of the Fourier transform, because these identities uniquely
determine $\cH(D)$ if $D$ is written $(\frac{d}{dx})^N
(1+x)^M g_{N,M}(x)$ with $g_{N,M}\in L^2(0,\infty;dx)$. But
the proof needs some organizing then as it is necessary to
check independence from the choice of $N$ and $M$, and also
to establish afterwards all identities above. So
\eqref{eq:25} provides the easiest road. Still, in this
context, let us mention the following which relates to the
restriction of $\cH(D)$ to $(0,+\infty)$:

\begin{lemma}\label{lem:1}
  Let $k$ be smooth on $\RR$ with compact support in
  $[0,+\infty)$. Then $\cH(k)$ is the restriction
  to $[0,+\infty)$ of an entire function $\gamma$ which has Schwartz
  decrease as $x\to+\infty$. For any tempered distribution
  $D$ with support in $[0,+\infty)$, there holds
  \begin{equation}\label{eq:40}
  \int_0^\infty \cH(D)(x)k(x)\,dx = \int_0^\infty D(x)\gamma(x)\,dx\;,
  \end{equation}
  where in the right hand side in fact one has
  $\int_{-\epsilon}^\infty D(x)\gamma(x)\theta(x)\,dx$ where
  the smooth function $\theta$ is $1$ for $x\geq -
  \frac\epsilon3$ and $0$ for $x\leq-\frac\epsilon2$ and is
  otherwise arbitrary (as is $\epsilon$).
\end{lemma}

Let us suppose $k=0$ for $x>B$. Defining:
\begin{equation}
  \label{eq:30}
  \gamma(x) = \int_0^B J_0(2\sqrt{xy})\, k(y)\,dy
\end{equation}
we obtain an entire function and, according to our
definitions,
  $\cH(k)(x) = \gamma(x)\Un_{x>0}(x)$ as a distribution or a
square-integrable function.  Using \eqref{eq:31b} ($\cH =
-\frac1x \cH\cdot \frac{d}{dx}\,x\, \frac{d}{dx}$
  for $x>0$) and bounding $J_0$ by $1$ we see (induction) that
$\gamma$ is $O(x^{-N})$ for any $N$ as $x\to+\infty$ , and
using \eqref{eq:31a} ($\frac{d}{dx}\cdot \cH = - \frac1x \cH
  \cdot  \frac{d}{dx}x$ for $x>0$) the same applies to its derivative and
also to its higher derivatives. So it is of the Schwartz
  class for $x\to+\infty$.

Replacing $D$ by $\cH(D)$ in \eqref{eq:40} it will be more convenient to
prove:
\begin{equation}
  \label{eq:33}
  \int_0^\infty D(x)k(x)\,dx = \int_0^\infty \cH(D)(x)\gamma(x)\,dx
\end{equation}
If \eqref{eq:33} holds for $D$ (and all $k$'s) then $<D',k>
= -<D,k'> = - <\cH(D),-\theta(x)\int_x^\infty\gamma(y)dy>$
(observe that $\int_0^x \gamma(y)\,dy = \cH(k')(x)$ vanishes
at $+\infty$) so $<D',k> = +<\int_0^x \cH(D),\theta\gamma> =
<\nobreak\cH(D'),\theta\gamma>$ hence \eqref{eq:33} holds as well
for $D'$ (and all $k$'s). So we may assume $D$ to be a
continuous function of polynomial growth. It is also checked
using \eqref{eq:31b} that if \eqref{eq:33} holds for $D$ it
holds for $xD$. So we may reduce to $D$ being
square-integrable, and the statement then follows from the
self-adjointness of $\cH$ on $L^2$ (or we reduce to Fubini).

The behavior of $\cH$ with respect to the translations
$\tau_a: f(x)\mapsto f(x-a)$ is important. For
$f\in L^2(\RR;dx)$ the value of $a$ is arbitrary and we can define
\begin{subequations}
\begin{gather}
  \label{eq:27a}
  \tau_a^\# : =\cH\,\tau_a\,\cH\\
  \wt{\tau_a(f)}(\lambda) = e^{ia\lambda}\wt{f}(\lambda)\\
  \wt{\tau_a^\#(f)}(\lambda) = e^{ia\frac{-1}\lambda}\wt{f}(\lambda)
\end{gather}
\end{subequations}
We observe the remarkable commutation relations
(which would fail for the cosine or sine transforms):
\begin{equation}
\label{eq:28}
  \forall a,b\qquad \tau_a\tau_b^\# = \tau_b^\#\tau_a
\end{equation}
For a distribution $D$ the action of $\tau_a^\#$ is here
defined only for $a\geq -\alpha(\cH(D))$, where $\alpha(E)$
is the leftmost point of the closed support of the
distribution $E$. On this topic from the validity of
\eqref{eq:19} when $f\in L^2(0,\infty;dx)$, and invariance
of $\alpha$ under derivation\footnote{It is important in
order to avoid a possible confusion to insist on the fact
that $\frac{d}{dx}$ is always taken in the distribution
sense so for example $\frac{d}{dx} \Un_{x>0} = \delta(x)$
indeed has the leftmost point of its support not affected by
$\frac{d}{dx}$.}, integration, and multiplication by $x$,
one has:
\begin{equation}
  \label{eq:19distrib}
  -\alpha(E) = \limsup_{t\to+\infty} \frac1t\,\log|\wt E(it)|
\end{equation}
We thus have the property, not shared by the cosine or
sine transforms:
\begin{equation}
  \label{eq:29}
  a\geq -\alpha(\cH(D))\implies \alpha(\tau_a^\#(D)) = \alpha(D)
\end{equation}

We now consider $D$ with $\alpha(D)>0$ and
$\alpha(\cH(D))>0$ and prove that its Mellin transform is an
entire function with trivial zeros at $0$, $-1$, $-2$,
\dots, following the method of regularization by
multiplicative convolution and co-Poisson intertwining
from \cite{copoisson}. The other, very classical in spirit,
proof  shall be presented later. The latter method is shorter but
the former provides complementary information.

In
\cite[\textsection4.A]{copoisson} the detailed explanations
relative to the notion of multiplicative convolution are
given:
\begin{equation}
  \label{eq:22}
  (g*D)(x)\ \text{``}=\text{''} \int_\RR g(t) D(\frac xt)\,\frac{dt}{|t|}\;,
\end{equation}
where we will in fact always take $g$ to have compact
support in $(0,+\infty)$. It is observed that
  \begin{equation}\label{eq:23}
    g*xD = x(\frac{g}x*D)\qquad (g*D)' = \frac{g}x*D'
  \end{equation}
The notion of right Mellin transform
$\int_0^\infty D(x)x^{-s}dx$ is developed in
\cite[\textsection4.C]{copoisson}, for $D$ with support in
$[a,+\infty)$, $a>0$:
\begin{equation}
  \label{eq:23a}
  \wh D(s) = s(s+1)\cdots(s+N-1)\wh{D^{(-N)}}(s+N)\;,
\end{equation}
where $N\gg0$.  The meaning of $\wh D$ is as the maximal
possible analytic continuation to a half-plane
$\Re(s)>\sigma$, where $\sigma$ is as to the left as is
possible. The notion is extended\footnote{if $D$ is near the
origin a function with an analytic character, then
straightforward elementary arguments allow a complementary
discussion. However if $D$ is just an element of
$L^2(0,\infty;dx)$ then $\wh D$ is a square-integrable
function on the critical line, and nothing more nor less.} in
\cite[\textsection4.F]{copoisson} to the case where the
restriction of $D$ to $(-a,a)$ is ``quasi-homogeneous''. For
example, if $\left. D\right|_{(-a,a)} = \Un_{0<x<a}$
(resp. $\delta$), then $\wh D$ is defined as $\wh D_1$ with
$D_1 = D - \Un_{0<x<\infty}$ (resp. $D- \delta$.) Then, also
in the extended case, the following holds:
\begin{equation}
  \label{eq:24}
  \wh{g*D}(s) = \wh g(s)\wh D(s)
\end{equation}
where $g$ in an integrable function with compact support in
$(0,\infty)$ and $\wh g(s)$ is the entire function
$\int_0^\infty g(t)t^{-s}\,dt$. We then have the following
theorem:

\begin{theorem}
  \label{theo:2}
Let $D$ a tempered distribution with support in
$[a,+\infty)$, $a>0$ and such that $\cH(D)$ also has a
positive leftmost point of support. Let $g$ be a smooth
function with compact support in $(0,\infty)$. Then the
multiplicative convolution $g*D$ belongs to the Schwartz
class.
\end{theorem}

This is the analog of \cite[Thm 4.29]{copoisson}. The function
$k(t) = (Ig)(t) = \frac{g(1/t)}t$ is defined and it is written
as $k = \cH(\gamma\Un_{x>0})$ where $\gamma$ is the entire
function, of Schwartz decrease at $+\infty$ such that
$\cH(k) = \gamma\cdot\Un_{x>0}$. Then it is observed that  
\begin{equation}
  \label{eq:32}
t>0\implies\quad (g*D)(t) = \int_0^\infty D(x)\frac{k(x/t)}t
\,dx = \int_0^\infty \cH(D)(x)\gamma(tx)\,dx
\end{equation}
We have used Lemma \ref{lem:1}. Then the Schwartz decrease of
$\int_0^\infty \cH(D)(x)\gamma(tx)\,dx$ as $t\to+\infty$ is
established as is done at the end of the proof of \cite[Thm
4.29]{copoisson}, integrating by parts enough times to
transform $\cH(D)$ into a continuous function of polynomial
growth, identically zero on $[0,c]$, $c>0$.

\begin{theorem}
  \label{theo:3}
Let $D$ a tempered distribution with a positive lefmost
point of support and such that  $\cH(D)$ also has a
positive leftmost point of support. Then $\wh D(s)$ and
$\Gamma(s)\wh D(s)$ are entire functions and:
\begin{equation}
  \Gamma(s) \wh D(s) = \Gamma(1-s) \wh{\cH(D)}(1-s)
\end{equation}
\end{theorem}

We first establish:

\begin{theorem}[``co-Poisson intertwining'']
Let $D$ be a tempered distribution supported in
$[0,+\infty)$ and let $g$ be an integrable function with
compact support in $(0,\infty)$. Then, with  $(Ig)(t) = \frac{g(1/t)}t$:
\begin{equation}\label{eq:41}
  \cH(g*D) = (Ig)*\cH(D)
\end{equation}
\end{theorem}

Let us first suppose that $D$ is an $L^2$ function. In that
case, we will use the Mellin-Plancherel transform $f\mapsto
\wh f(s) = \int_0^\infty f(t)t^{-s}\,dt$, for $f$ square integrable and
$\Re(s) = \frac12$. Then $\wh{g*f}$ is,
changing variables,  the Fourier transform of an additive
convolution where one of the two has compact support, well
known to be the product $\wh g\cdot \wh f$. We need also to
understand the Mellin transform of $\cH(f)$. Let us suppose
$f_t(x) = \exp(-tx)$. Then $\cH(f_t) = \frac1t f_{\frac1t}$
has Mellin transform $\wh{\cH(f_t)}(s) = t^{-s}\Gamma(1-s)$
and $\wh{f_t}(s) = t^{s-1}\Gamma(1-s)$, so we have the
identity for such $f$'s:
\begin{equation}\label{eq:42}
  \wh{\cH(f)}(s) = \frac{\Gamma(1-s)}{\Gamma(s)}\wh f(1-s)
\end{equation}
The linear combinations of the $f_t$'s are dense in $L^2$,
so \eqref{eq:41} holds for all $f$'s as an identity of
square integrable functions on the critical line. We are now
in a position to check the intertwining: $\wh{\cH(g*f)}(s) =
\frac{\Gamma(1-s)}{\Gamma(s)} \wh g(1-s)\wh f(1-s) =
\wh{Ig}(s)\wh{\cH(f)}(s) = \wh{Ig*\cH(f)}(s)$.

For the case of an arbitrary distribution it will then be
sufficient to check that if \eqref{eq:41} holds for $D$ it
holds for $xD$ and for $D'$. This is easily done using
\eqref{eq:23}. We have $g*(D') = (xg * D)'$, so $\cH(g*D') =
\int_0^x \cH(xg * D) = \int_0^x (\frac{Ig}x * \cH(D)) =
Ig * (\int_0^x \cH(D)) = Ig *\cH(D')$. A similar proof is
done for $xD$. This completes the proof of the intertwining.

The theorem \ref{theo:3} is then established as is \cite[Thm
4.30]{copoisson}. We pick an arbitrary $g$ smooth with compact
support in $(0,\infty)$. We know by theorem \ref{theo:2} that $g*D$
is a Schwartz function as $x\to+\infty$, and certainly it
vanishes identically in a neighborhood of the origin, so
$\wh{g*D}(s) = \wh g(s)\wh D(s)$ is an entire function. So
$\wh D(s)$ is a meromorphic function in the entire complex
plane, in fact an entire function as $g$ is arbitrary. We
then use the intertwining and \eqref{eq:42} for square
integrable functions. This gives $\wh
g(1-s)\wh{\cH(D)}(s) = \wh{Ig*\cH(D)}(s) = \wh{\cH(g*D)}(s)
= \frac{\Gamma(1-s)}{\Gamma(s)} \wh g(1-s)\wh
D(1-s)$. Hence, indeed, after replacing $s$ by $1-s$:
\begin{equation}
  \Gamma(s) \wh D(s) = \Gamma(1-s) \wh{\cH(D)}(1-s)
\end{equation}
The left-hand side may have poles only at $0$, $-1$, \dots,
and the right-hand side only at $1$, $2$, \dots. So both
sides are entire functions and $\wh D(s)$ has trivial zeros
at $0$, $-1$, $-2$, \dots

We now give another proof of Theorem \ref{theo:3}, which is
more classical, as it is the descendant of the second of
Riemann's proof, and is the familiar one from the theory of 
theory of $L$-series and modular functions. The existence of
two complementary proofs is  instructive, as it helps to
better understand the rôle of the right Mellin
transform $\int_0^\infty f(x)x^{-s}\,dx$ vs. the left Mellin
transform $\int_0^\infty \theta(it)t^{s-1}\,dt$. 

To the distribution $D$ we associate its ``theta''
function\footnote{the author hopes to be forgiven this
temporary terminology in a situation
where only the behavior under
$\lambda\mapsto \frac{-1}\lambda$ is at work.}
$\theta_D(\lambda) = \wt D(\lambda) = \int_0^\infty
e^{i\lambda x}D(x)\,dx$, which is an analytic function for
$\Im(\lambda)>0$ \footnote{we adopt the usual notation, and
consider $\theta_D$ as a function of $it$ rather than $t$.}.  Right from
the beginning we have:
\begin{equation}
  \theta_{\cH(D)}(it) = \frac1t \theta_D(\frac it)
\end{equation}
If the leftmost point of the support of $D$ is positive then
$\theta_D(it)$ has exponential decrease as $t\to+\infty$ and
$\int_1^\infty \theta_D(it)t^{s-1}\,dt$ is an entire
function. If also the leftmost point of support of $\cH(D)$
is positive then $\theta_{\cH(D)}(it)$ has exponential
decrease as $t\to+\infty$ and $\int_0^1
\theta_D(it)t^{s-1}\,dt = \int_1^\infty \theta_{\cH(D)}(it)
t^{-s}\,dt$ is an entire function. So, under the support
property considered in Theorem \ref{theo:3} $\cD(s) :=
\int_0^\infty \theta_D(it)t^{s-1}\,dt$ is indeed an entire
function, and the functional equation is
\begin{equation}
  \cD(s) = \cD^*(1-s)
\end{equation}
with $\cD^*(s) = \int_0^\infty
\theta_{\cH(D)}(t)t^{s-1}\,dt$.

To conclude we also need to establish:
\begin{equation}\label{eq:45}
  \cD(s) = \Gamma(s)\wh D(s)
\end{equation}
We shall prove this for $\Re(s)\gg0$ under the hypothesis
that $D$ has support in $[a,+\infty)$, $a>0$ (no hypothesis
on $\cH(D)$). In that case, as $\theta_D(it)$ is $O(t^{-N})$
for a certain $N$ as $t\to 0$ ($t>0$), and is of exponential
decrease as $t\to+\infty$, we can define $\cD(s) =
\int_0^\infty \theta_D(it)t^{s-1}\,dt$ as an analytic
function for $\Re(s)\gg0$. Let us suppose that $D$ is a
continuous function which is $O(x^{-2})$ as
$x\to+\infty$. Then, for, $\Re(s)>0$, the identity
\eqref{eq:45} holds as an application of the Fubini
theorem. We then apply our usual method to check that if
\eqref{eq:45} holds for $D$ it also holds for $xD$ and for
$D'$. For this, obviously we need things such as $\wh{D'}(s)
= s\wh D(s+1)$ \cite[4.15]{copoisson} and $\wh{xD}(s) = \wh
D(s-1)$, the formulas $\theta_{D'} = -i\lambda\theta_D$,
$\theta_{xD} = -i\frac{\partial}{\partial\lambda}\theta_D$,
and $\Gamma(s+1) = s\Gamma(s)$. The verifications are then
straightforward.

In summary we have seen how the support property for $D$ and
$\cH(D)$ is related in two complementary manners to the
functional equation, one using the right Mellin transform
$\wh D(s)$ of $D$ and the idea of co-Poisson, the other
using the left Mellin transform $\cD(s)$ of the ``theta''
function $\theta_D$ associated to $D$ as an analytic
function on the upper half-plane and the behavior of
$\theta_D(it)$ under $t\mapsto \frac1t$. It is possible to
push further the analysis and to characterize the class of
entire functions $\cD(s) = \Gamma(s)\wh D(s)$, as has been
done in \cite{copoisson} in the case of the cosine and sine
transforms. It is also explained in \cite{copoisson} how the
discussion extends to allow finitely many poles. The proofs
and statements given there are easily adapted to the case of
the $\cH$ transform. Only the case of poles at $1$ and $0$
will be needed here and this corresponds, either to the
condition that $D$ and $\cH(D)$ both restrict in $(-a,a)$
for some $a>0$ to multiples of the Dirac delta function, or,
that they are both constant in $[0,a)$ for some $a>0$. We
recall that the Mellin transform $\wh D(s)$ is defined in
such a manner, that it is not affected from either
substracting $\delta$ or $\Un_{x>0}$ from $D$.


\section{A group of distributions and related integral formulas}

We now derive some integral identities which will prove
central. The identities will be re-obtained later as the
outcome of a less direct path. We are interested in the
tempered distribution $g_a(x)$ whose Fourier transform is
$\exp(i a \frac{-1}\lambda)$. Indeed $\tau_a^\#(f)$
(equation \eqref{eq:27a}) is the
additive convolution of $f$ with $g_a$: we note that
$g_a$
    differs from $\delta(x)$ by a square integrable function
    as $1 - \exp(-ia\lambda^{-1})=
    O_{|\lambda|\to\infty}(|\lambda|^{-1})$; so there is
a convolution formula
  $\tau_a^\#(f) = f - f_a*f$ for a certain square integrable
function $f_a$. For $f\in L^2$, the convolution $f_a*f$ as the Fourier
transform of an $L^1$-function is continuous on $\RR$.
Starting from the identity
    $\exp(i a
\frac{-1}\lambda) = - i \lambda \frac{i}\lambda\exp(i a
\frac{-1}\lambda)$ we identify $g_a$ for $a\geq0$ as
$\frac{\partial}{\partial x} \cH \delta_a$. It is important
that $\frac{\partial}{\partial x}$ is taken in the
distribution sense. So we have, simply:
\begin{equation}
  \label{eq:142}
  g_a(x) = \delta(x) -
  \frac{a J_1(2\sqrt{ax})}{\sqrt{ax}}\Un_{x>0}(x)\qquad
  (a\geq 0)
\end{equation}
If $a<0$ then $g_a(x) = g_{-a}(-x)$, $f_a(x) = f_{-a}(-x)$. So:
\begin{equation}
  \label{eq:143}
  g_{-a}(x) = \delta(x) -
  \frac{ a J_1(2\sqrt{-ax})}{\sqrt{-ax}}\Un_{x<0}(x)\qquad
  (-a\leq 0)
\end{equation}
The group property under the additive convolution $g_a*g_b =
g_{a+b}$ leads to remarkable integral identities $f_{a+b} =
f_a + f_b - f_a * f_b$ involving the Bessel
functions. The pointwise validity
is guaranteed by continuity; the Plancherel
identity confirms  the identity, where $f_a(x) =  \frac{a
J_1(2\sqrt{ax})}{\sqrt{ax}}\Un_{x>0}(x)$ for $a\geq 0$ and
$f_{-a}(x) = f_a(-x)$:
\begin{equation}
  \label{eq:155}
  f_{a+b} = f_a + f_b - f_a * f_b
\end{equation}
At $x=0$ the pointwise identity is obtained by continuity
    from either $x>0$ or $x<0$. We have essentially two cases: $g_a*g_b$ for
$a,b\geq0$ and
    $g_a*g_{-b}$ for $a\geq b \geq0$. The following is obtained:
\begin{subequations}
\begin{proposition}
  Let $a\geq b\geq 0$ and $x\geq0$. There holds:
  \begin{equation}
    \label{eq:144}
    \frac{(a+b)J_1(2\sqrt{(a+b)x})}{\sqrt{(a+b)x}} =
    \frac{a J_1(2\sqrt{ax})}{\sqrt{ax}} + 
    \frac{b J_1(2\sqrt{bx})}{\sqrt{bx}} - \int_0^x \frac{a
    J_1(2\sqrt{ay})}{\sqrt{ay}}\;\frac{ b
    J_1(2\sqrt{b(x-y)})}{\sqrt{b(x-y)}}\,dy
  \end{equation}
  \begin{equation}
    \label{eq:145}
     \frac{(a-b)J_1(2\sqrt{(a-b)x})}{\sqrt{(a-b)x}} =
    \frac{a J_1(2\sqrt{ax})}{\sqrt{ax}}  - \int_x^\infty \frac{a
    J_1(2\sqrt{a y})}{\sqrt{ay}}\;\frac{ b
    J_1(2\sqrt{b(y-x)})}{\sqrt{b(y-x)}}\,dy
  \end{equation}
  \begin{equation}
    \label{eq:146}
     0 =
    \frac{b J_1(2\sqrt{bx})}{\sqrt{bx}}  - \int_0^\infty \frac{a
    J_1(2\sqrt{a y})}{\sqrt{ay}}\;\frac{ b
    J_1(2\sqrt{b(y+x)})}{\sqrt{b(y+x)}}\,dy
  \end{equation}
\end{proposition}
\end{subequations}

Exchanging $a$ and $b$ and changing variables we combine
    \eqref{eq:145} and \eqref{eq:146} into one single
    equation for $x\geq0$ and $a,b\geq0$:
\begin{equation}
  \label{eq:147}
    \frac{(a-b)J_1(2\sqrt{(a-b)x})}{\sqrt{(a-b)x}}\Un_{a-b\geq0}(a-b)=
    \frac{a J_1(2\sqrt{ax})}{\sqrt{ax}}  - \int_0^\infty \frac{a
    J_1(2\sqrt{a(y+x)})}{\sqrt{a(y+x)}}\;\frac{ b
    J_1(2\sqrt{b y})}{\sqrt{b y}}\,dy
\end{equation}
The formula for $x=0$ in \eqref{eq:147} is obtained by
continuity. It is equivalent to
\begin{equation}
  \label{eq:156}
  \int_0^\infty J_1(u)J_1(cu)\,\frac{du}u = \frac12\min(c,\frac1c)\qquad(c>0)
\end{equation}
which is a very special case of formulas of Weber, Sonine
and Schafheitlin (\cite[13.42.(1)]{watson}). Another
interesting special case of \eqref{eq:147} is for $a=b$. The
formula becomes
\begin{equation}
  \label{eq:157}
  \frac{J_1(2\sqrt x)}{\sqrt x} = \int_0^\infty
  \frac{J_1(2\sqrt y)}{\sqrt y}
  \frac{J_1(2\sqrt{x+y})}{\sqrt{x+y}} \,dy
\end{equation}
which is equivalent to a special case  of a formula of Sonine
(\cite[13.48.(12)]{watson}).

We already mentioned the equation
    $\frac{\partial^2}{\partial u\partial v} J_0(2\sqrt{uv})
    = - J_0(2\sqrt{uv})$. New identities are obtained from
    \eqref{eq:147} or \eqref{eq:144} after taking either the
    $a$ or the $b$ derivative. We investigate no further
    \eqref{eq:147} as
    the corresponding  semi-convergent
    integrals, in a form or another, are
    certainly among the formulas of
    \cite[\textsection13]{watson}. Let us rather focus more
    closely on the case $a,b\geq0$ (\eqref{eq:144}.) We have a
    function which is entire in $a$, $b$, and $x$ and the
    identity holds for all complex values of $a$, $b$, and
    $x$. Let us take the derivative with respect to $a$:
    \begin{equation}
      \label{eq:148}
    J_0(2\sqrt{(a+b)x})=
    J_0(2\sqrt{ax})   - \int_0^x J_0(2\sqrt{ay}) \;\frac{ b
    J_1(2\sqrt{b(x-y)})}{\sqrt{b(x-y)}}\,dy
    \end{equation}
We replace $b$ by $-b$ and then set $x=b$. This gives:
\begin{equation}
  \label{eq:149}
    I_0(2\sqrt{b(b-a)})=
    J_0(2\sqrt{ba})   + \int_0^b J_0(2\sqrt{ay}) \;\frac{ b
    I_1(2\sqrt{b(b-y)})}{\sqrt{b(b-y)}}\,dy
\end{equation}
We take the derivative of \eqref{eq:148} with respect to
    $b$:
    \begin{equation}
      \label{eq:150}
    -\frac{x J_1(2\sqrt{(a+b)x})}{\sqrt{(a+b)x}} =
      - \int_0^x J_0(2\sqrt{ay}) J_0(2\sqrt{b(x-y)})\,dy
    \end{equation}
Then we replace $b$ by $-b$ and set $x=b$:
\begin{equation}
  \label{eq:151}
    \frac{b I_1(2\sqrt{b(b-a)})}{\sqrt{b(b-a)}} =
       \int_0^b J_0(2\sqrt{ay}) I_0(2\sqrt{b(b-y)})\,dy
\end{equation}
Combining \eqref{eq:149} and \eqref{eq:151} by  addition and 
substraction  we discover that we have solved certain
    integral equations:
\begin{subequations}
  \begin{equation}
    \label{eq:153a}
    \phi_b^+(x) = I_0(2\sqrt{b(b-x)}) - \frac{b
    I_1(2\sqrt{b(b-x)})}{\sqrt{b(b-x)}} = (1 +
    \frac\partial{\partial x}) I_0(2\sqrt{b(b-x)}) 
  \end{equation}
  \begin{equation}
    \label{eq:153b}
    \phi_b^-(x) = I_0(2\sqrt{b(b-x)}) + \frac{b
    I_1(2\sqrt{b(b-x)})}{\sqrt{b(b-x)}}  =  (1 -
    \frac\partial{\partial x}) I_0(2\sqrt{b(b-x)})
  \end{equation}
  \begin{equation}
    \label{eq:153c}
    \phi_b^+(x) + \int_0^b J_0(2\sqrt{xy})\phi_b^+(y)\,dy = J_0(2\sqrt{bx})
  \end{equation}
  \begin{equation}
    \label{eq:153d}
    \phi_b^-(x) - \int_0^b J_0(2\sqrt{xy})\phi_b^-(y)\,dy = J_0(2\sqrt{bx})
  \end{equation}
\end{subequations}
The significance will appear later in the paper and we leave
the matter here. The method was devised after the importance
of solving equations \eqref{eq:153c} and \eqref{eq:153d} had
    emerged and after the solutions \eqref{eq:153a} and
    \eqref{eq:153b} had been obtained as the outcome of a
    more indirect path. Of course, direct verification by
    replacement of the Bessel functions by their series
    expansions is possible and easy.

\section{Orthogonal projections and Hilbert space
evaluators}
\label{sec:5}

Let $a>0$ and let $P_a$ be the orthogonal projection on
$L^2(0,a;dx)$ and $Q_a = \cH P_a \cH$ the orthogonal
projection on $\cH(L^2(0,a;dx))$ and let $K_a\subset
L^2(0,\infty;dx)$ be the Hilbert space of square integrable
functions $f$ such that both $f$ and $\cH(f)$ have their
supports in $[a,\infty)$. Also we shall write $H_a = P_a\cH
P_a$. Also we shall very often use $D_a = H_a^2 = P_a \cH
P_a \cH P_a$.  Using:
\begin{equation}
  J_0(2\sqrt{xy}) = \sum_{n=0}^\infty (-1)^n \frac{x^n y^n}{n!^2}\;,
\end{equation}
we exhibit $H_a = P_a \cH P_a$ as a limit in operator norm
of finite rank operators so $P_a\cH P_a$ is a compact
(self-adjoint) operator. It is not possible for a non zero
$f\in L^2(0,a;dx)$ to be such that $\|H_a(f)\| = \|f\|$, as
this would imply that $H_a(f)$ vanishes identically for
$x>a$, but $H_a(f)$ is an entire function. So the operator
norm of $H_a$ is strictly less than one, and $1\pm H_a$ as
well as $1 -D_a$ are invertible.  We consider the
equation
\begin{equation}
  \phi = u + \cH(v)\qquad u,v \in L^2(0,a;dx)
\end{equation}
Hence:
\begin{subequations}
  \begin{align}
    u + H_a(v) &= P_a(\phi)\\
   H_a(u) + v &= P_a(\cH(\phi))\\
   u &= (1 - D_a)^{-1} (P_a(\phi) - H_aP_a\cH(\phi))\\
   v &= (1 - D_a)^{-1} (-H_aP_a(\phi) + P_a\cH(\phi))
  \end{align}
\end{subequations}
Then if $\phi_n = u_n + \cH(v_n)$ is $L^2$-convergent,
$(u_n)$ and $(v_n)$ will be convergent, and the vector space
sum $L^2(0,a;dx)+ \cH(L^2(0,a;dx))$ is closed. Its elements
are analytic functions for $x>a$ so certainly this is a
proper subspace of $L^2$. Hence we obtain that each $K_a$ is
not reduced to $\{0\}$ and 
\begin{equation}
  K_a^\perp = L^2(0,a;dx)+\cH(L^2(0,a;dx))
\end{equation}
We also mention that $\cup_{a>0} K_a$ is dense in but not
equal to $L^2(0,\infty;dx)$, more generally that $\cup_{a>b}
K_a$ is dense in but not equal to $K_b$, and also obviously
$\cap_{a<\infty} K_a = \{0\}$, $\cap_{a<b} K_a = K_b$.

In this section $a>0$ will be fixed (all defined quantities
and functions will depend on $a$, but this will not always
be explicitely indicated.) We shall be occupied
with understanding the vectors $X_s^a\in K_a$ such that
\begin{equation}
  \forall f\in K_a\qquad \int_a^\infty f(x)X_s^a(x)\,dx = \wh
  f(s) = \int_a^\infty f(x) x^{-s}\,dx
\end{equation}
and in particular we are interested in computing
\begin{equation}
  X_a(s,z) = \int_a^\infty X_s^a(x) X_z^a(x)\,dx
\end{equation}
As $a$ is fixed here, we shall drop the superscript $a$ to
lighten the notation.  For the time being we shall restrict
to $\Re(s)>\frac12$ and we define $X_s$ to be the orthogonal
projection to $K_a$ of $\Un_{x>a}(x) x^{-s}$. As a
preliminary to this study we need to say a few words
regarding:
\begin{equation}\label{eq:55}
  g_s(x) := \cH(\Un_{x>a}(x) x^{-s}) = \int_a^\infty J_0(2\sqrt{xy})y^{-s}\,dy
\end{equation}
The integral is absolutely convergent for  $\Re(s)>\frac34$,
semi-convergent for $\Re(s)>\frac14$,
and $g_s$ is defined by the equation as an $L^2$ function
  for $\Re(s)>\frac12$ (it will prove to be entire in $s$
  for each $x>0$). We need the following identity, which
shows also that $g_s(x)$ is
  analytic in $x>0$:
\begin{equation}  \label{eq:50}
  g_s(x) = \chi(s) x^{s-1} - \int_0^a J_0(2\sqrt{xy})
  y^{-s}\,dy
=  \chi(s) x^{s-1} - \sum_{n=0}^\infty (-1)^n
  \frac{x^n\, a^{n+1-s}}{n!^2(n+1-s)}
\end{equation}
This is obtained first in the range $\frac34<\Re(s)<1$:
$\int_a^\infty J_0(2\sqrt{xy})y^{-s}\,dy =
x^{s-1}\int_{ax}^\infty J_0(2\sqrt y) y^{-s}\,dy =
x^{s-1}\left(\chi(s) - \int_0^{ax} J_0(2\sqrt y)
y^{-s}\,dy\right) = x^{s-1}\chi(s) - \int_0^a
J_0(2\sqrt{xy})y^{-s}\,dy$. The poles at $s=1$, $s=2$, \dots
are only apparent. The identity is valid by analytic
continuation in the entire plane $\Re(s)>\frac12$. For each
given $x>0$ we have in fact an entire function of
$s\in\CC$. But we are here more interested in $g_s$ as a
function of $x$ and we indeed see that it is analytic in
$\CC\setminus]-\infty,0]$ (it is an entire function of $x$
if $s\in-\NN$). \footnote{For some other transforms $k(xy)$,
such as the cosine transform, the argument must be slightly
modified in order to accomodate the fact $\int_0^\infty k(y)
y^{-s}\,dy$ has no range of absolute convergence.}


There are unique vectors $u_s, v_s$ in $L^2(0,a;dx)$ such
that
\begin{equation}\label{eq:51}
  \Un_{x>a}(x)x^{-s} = X_s(x) + u_s(x) + \cH(v_s)(x)\,
\end{equation}
and they are the solutions to the system of equations:
\begin{subequations}
  \begin{align}\label{eq:52a}
    u_s + H_a(v_s) &=0\\ \label{eq:52b}
   H_a(u_s) + v_s &= P_a(g_s)
  \end{align}
\end{subequations}
From \eqref{eq:52a} we see that $u_s$ is in fact the
restriction to $(0,a)$ of an entire funtion, and from
\eqref{eq:52b} that $v_s$ is the restriction to $(0,a)$ of a
function which is analytic in $\CC\setminus]-\infty,0]$. Redefining
$u_s$ and $v_s$ to now refer to these analytic functions
their defining equations become (on $(0,+\infty)$):
\begin{subequations}
  \begin{align}\label{eq:53a}
    u_s + \cH P_a(v_s) &=0\\ \label{eq:53b}
   \cH P_a(u_s) + v_s &= g_s
  \end{align}
\end{subequations}
and \eqref{eq:51} becomes (we set $X_s(a) = X_s(a+)$):
\begin{subequations}
\begin{align}
\Un_{x\geq a}(x)x^{-s} &= X_s(x) + \Un_{0<x< a}(x)u_s(x) + \cH P_a(v_s)(x)\\
\Un_{x\geq a}(x)x^{-s} &= X_s(x) - \Un_{x\geq a}(x)u_s(x)\\
X_s(x) &= \Un_{x\geq a}(x) (x^{-s} + u_s(x)) \label{eq:56}
\end{align}
\end{subequations}
The key to the next steps will be the idea to investigate
the distribution $(x\frac{d}{dx} + s)X_s$ on the (positive) real
line. Let $D_s$ be $x\frac{d}{dx} + s$. There holds:
\begin{equation}
  D_s\cH = - \cH D_{1-s}
\end{equation}
To compute $\frac{d}{dx} P_a(v_s)$ we first suppose
$\Re(s)>1$, so (we know the behavior as $x\to0$ from
\eqref{eq:50}) $\frac{d}{dx} P_a(v_s) = P_a(v_s') -
v_s(a)\delta_a(x)$ and $x\frac{d}{dx} P_a(v_s) = P_a(xv_s')
- a v_s(a)\delta_a(x)$. This remains true for
$\Re(s)>\frac12$. Applying $D_s$
  to \eqref{eq:53a} thus
gives $D_s(u_s) - \cH \left(P_a D_{1-s}(v_s) -
av_s(a)\delta_a(x)\right)$. We similarly apply $D_{1-s}$ to
\eqref{eq:53b} and obtain the following system:
\begin{subequations}
  \begin{align}
  D_s(u_s)(x) - (\cH P_a D_{1-s} v_s)(x) &= - a v_s(a) J_0(2\sqrt{ax})\\
   -(\cH P_a D_s u_s )(x) + D_{1-s}(v_s)(x) &= (D_{1-s}
  g_s)(x) - au_s(a) J_0(2\sqrt{ax})
  \end{align}
\end{subequations}
From \eqref{eq:55}, we have
$D_{1-s} g_s = - \cH D_s(\Un_{x>a}x^{-s}) = -
\cH(a^{1-s}\delta_a(x)) = - a^{1-s} J_0(2\sqrt{ax})$. Let us
define
\begin{equation}
  J_0^a(x) = J_0(2\sqrt{ax})
\end{equation}
We have proven:
\begin{subequations}
  \begin{align}
   + D_s u_s - \cH P_a D_{1-s} v_s &= - a v_s(a) J_0^a\\
   -\cH P_a D_s u_s   + D_{1-s} v_s  &= - a (a^{-s} + u_s(a)) J_0^a
  \end{align}
\end{subequations}
Restricting to the interval $(0,a)$ and solving, we find:
\begin{subequations}
  \begin{align}\label{eq:57a}
    P_aD_s u_s &= - a (1 - D_a)^{-1}( v_s(a) J_0^a + (a^{-s}
    + u_s(a)) H_aJ_0^a)\\ \label{eq:57b}
    P_aD_{1-s} v_s &= - a (1 - D_a)^{-1}((a^{-s} +
    u_s(a))J_0^a + v_s(a) H_aJ_0^a)
  \end{align}
\end{subequations}
It is advantageous at this stage to define $\phi_a^+$ and
$\phi_a^-$ to be the solutions of the equations (in $L^2(0,a;dx)$):
\begin{subequations}
  \begin{align}\label{eq:61a}
    \phi_a^+ + H_a \phi_a^+ &= J_0^a\\
\label{eq:61b}
    \phi_a^- - H_a \phi_a^- &= J_0^a
  \end{align}
\end{subequations}
We already know from \eqref{eq:153a} and \eqref{eq:153b}
exactly what $\phi_a^+$ and $\phi_a^-$ are (in this special
case of the $\cH$ transform), but we shall proceed as if we
didn't. We see from \eqref{eq:61a}, \eqref{eq:61b} that
$\phi_a^+$ and $\phi_a^-$ are entire functions, and we can
rewrite the system as:\footnote{in conformity with our
conventions, these are identities on $(0,\infty)$; to see
them as identities on $\CC$ one must read $\int_0^a
J_0(2\sqrt{xy})\phi_a^+(y)\,dy$ rather than $(\cH
P_a\phi_a^+)(x)$.}
\begin{subequations}
  \begin{align}\label{eq:60a}
    \phi_a^+ + \cH P_a \phi_a^+ &= J_0^a\\
\label{eq:60b}
    \phi_a^- - \cH P_a \phi_a^- &= J_0^a
  \end{align}
\end{subequations}
We observe the identities:
\begin{subequations}
  \begin{align}\label{eq:70a}
    (1 - D_a)^{-1} J_0^a &= P_a \frac{\phi_a^+ +
    \phi_a^-}2\\\label{eq:70b}
    (1 - D_a)^{-1} H_a J_0^a &= P_a \frac{- \phi_a^+ +
    \phi_a^-}2
  \end{align}
\end{subequations}
So \eqref{eq:57a} and \eqref{eq:57b} become
\begin{subequations}
  \begin{align}\label{eq:58a}
    D_s u_s &= +a\frac{a^{-s}+u_s(a)-v_s(a)}2\phi_a^+
    -a\frac{a^{-s}+u_s(a)+v_s(a)}2\phi_a^-\\
    D_{1-s} v_s &= -a\frac{a^{-s} + u_s(a) -
    v_s(a)}2\phi_a^+ - a\frac{a^{-s}+u_s(a) + v_s(a)}2\phi_a^-
  \end{align}
\end{subequations}
From \eqref{eq:58a} we compute successively (again,
these are identities on $(0,+\infty)$):
\begin{equation}
  \cH P_a D_s u_s = a\frac{a^{-s}+u_s(a)-v_s(a)}2 (J_0^a -
  \phi_a^+) -a\frac{a^{-s}+u_s(a)+v_s(a)}2 (-J_0^a + \phi_a^-)
\end{equation}
\begin{equation}\label{eq:59}
  P_a D_s u_s = a\frac{a^{-s}+u_s(a)-v_s(a)}2 (\delta_a -
  \cH\phi_a^+) -a\frac{a^{-s}+u_s(a)+v_s(a)}2 (-\delta^a + \cH\phi_a^-)
\end{equation}
In \eqref{eq:59}, $\cH\phi_a^+$ should perhaps be more
precisely written as $\cH(\phi_a^+\Un_{x>0})$. From
\eqref{eq:60a} we know that $\phi_a^+\Un_{x>0}$ is tempered
as a distribution.  From \eqref{eq:56} we compute $D_s X_s =
\Un_{x>a} D_s(u_s) + a(a^{-s}
+ u_s(a))\delta_a(x) = D_s u_s - P_a D_s u_s + a(a^{-s}
+ u_s(a))\delta_a(x)$. From \eqref{eq:58a} and\eqref{eq:59}
then follows:
\begin{equation}
  \begin{split}
    D_s X_s = +a\frac{a^{-s}+u_s(a)-v_s(a)}2(\phi_a^+ + 
    \cH\phi_a^+ - \delta_a)
    -a\frac{a^{-s}+u_s(a)+v_s(a)}2(\phi_a^- - \cH\phi_a^- +
    \delta^a) \\+  a(a^{-s}
+ u_s(a))\delta_a(x)
  \end{split}
\end{equation}
And the  result of the computation is:
\begin{equation}
  D_s X_s = +a\frac{a^{-s}+u_s(a)-v_s(a)}2 (\phi_a^+ + 
    \cH\phi_a^+) 
+a\frac{a^{-s}+u_s(a)+v_s(a)}2(-\phi_a^- + \cH\phi_a^-)
\end{equation}
We then define the remarkable distributions:
\begin{subequations}
  \begin{align}
    A_a &= \frac{\sqrt a}2(\phi_a^+ + \cH\phi_a^+)\\
  -iB_a &= \frac{\sqrt a}2(-\phi_a^-+ \cH\phi_a^-)\\
   E_a &= A_a - i B_a
  \end{align}
\end{subequations}
From \eqref{eq:61a} we observe that $A_a$ 
has its support in $[a,\infty)$. Furthermore it is $\cH$
invariant. Similarly, $-iB_a$, which is $\cH$ anti
invariant, also has its support in $[a,+\infty)$. We recover
$A_a$ and $-iB_a$ from $E_a$ through taking the invariant
and anti-invariant parts. We may also rewrite $D_s X_s$ as:
\begin{equation}\label{eq:72}
  D_s X_s = \sqrt{a}(a^{-s} + u_s(a))E_a - \sqrt{a}\,v_s(a)\cH E_a
\end{equation}
Some other manners of writing
$A_a$ and $-iB_a$ are useful: from \eqref{eq:60a}
$\cH\phi_a^+ = \delta_a - P_a\phi_a^+$ and from \eqref{eq:60b} $\cH\phi_a^- =
\delta_a + P_a\phi_a^-$, so:
\begin{subequations}
  \begin{align}
    A_a &= \frac{\sqrt a}2(\delta_a + \phi_a^+\Un_{x>a})\\
  -iB_a &= \frac{\sqrt a}2(\delta_a - \phi_a^-\Un_{x>a})
  \end{align}
\end{subequations}
And, we take also notice of the following definitions and identities:
\begin{subequations}
  \begin{align}\label{eq:350a}
    j_a &= \sqrt{a}(\delta_a - \phi_a^+\Un_{0<x<a})&\qquad
  j_a &= \sqrt{a}\cH\phi_a^+ &\qquad A_a &= \frac12(j_a +
  \cH j_a)\\
\label{eq:350b}
  -ik_a &= \sqrt{a}(\delta_a + \phi_a^-\Un_{0<x<a})&\qquad
  -ik_a &= \sqrt{a}\cH\phi_a^- &\qquad  B_a &= \frac12(k_a - \cH k_a)
  \end{align}
\end{subequations}

From \eqref{eq:60a} and \eqref{eq:60b} we know that
$\phi_a^+$ and $\phi_a^-$ are bounded, so the right Mellin
transforms are defined directly for $\Re(s)>1$
by:\footnote{the integral for $\wh{E_a}(s)$ is certainly
  absolutely convergent for $\Re(s)>\frac12$ as $\phi_a^+ -
  \phi_a^-$ is square integrable on $(0,\infty)$, and in
    fact it is absolutely convergent for $\Re(s)>\frac14$. As
we know already completely explicitely  $\phi_a^+$
  and $\phi_a^-$, we do not pause on this here. A general argument
suitable to establish in more general cases absolute
convergence for $\Re(s)>\sigma$ for some
$\sigma<\frac12$  will be given later.}
\begin{subequations}
  \begin{align}\label{eq:89a}
    \wh{A_a}(s) &= \frac{\sqrt a}2 \; \Big( a^{-s} + \int_a^\infty
    \phi_a^+(x)x^{-s}\,dx \Big)\\ \label{eq:89b}
    -i \wh{B_a}(s) &= \frac{\sqrt a}2\;  \Big( a^{-s} -
    \int_a^\infty \phi_a^-(x) x^{-s}\,dx \Big)
\end{align}
\begin{align} \label{eq:89c}
    \wh{E_a}(s) &= \sqrt a \; \Big(a^{-s} + \frac12\int_a^\infty
    (\phi_a^+(x) - \phi_a^-(x))x^{-s}\,dx \Big)\\ \label{eq:89d}
   \wh{\cH(E_a)}(s) &= \sqrt a \; \frac12\;\int_a^\infty
    (\phi_a^+(x) + \phi_a^-(x))x^{-s}\,dx
  \end{align}
\end{subequations}
From \eqref{eq:60a} and \eqref{eq:60b} we know that $\phi_a^+
- \phi_a^-$ is square-integrable at $+\infty$, so, using
$\cH\phi_a^+ = \delta_a - P_a\phi_a^+$ and  $\cH\phi_a^- =
\delta_a + P_a\phi_a^-$  we compute:
\begin{equation}
   \int_a^\infty (\phi_a^+(x)-\phi_a^-(x))x^{-s}\,dx =
   \int_0^\infty (\cH\phi_a^+ - \cH\phi_a^-)g_s(x)\,dx =
   -\int_0^a (\phi_a^+(x) + \phi_a^-(x))g_s(x)\,dx
\end{equation}
Then using \eqref{eq:70a}:
\begin{equation}\label{eq:9}
   \int_0^a \frac{\phi_a^+(x) + \phi_a^-(x)}2 g_s(x)\,dx =
   \int_0^a (1 - D_a)^{-1}(J_0^a)(x)g_s(x)\,dx = \int_0^a
   J_0^a(x)((1-D_a)^{-1}(g_s))(x)\,dx 
\end{equation}
Comparing with \eqref{eq:52a} and \eqref{eq:52b} the
  right-most term of \eqref{eq:9} may be
written as $\int_0^a J_0(2\sqrt{ax}) v_s(x)\,dx$ which in
turn we recognize from \eqref{eq:53a} to be $-u_s(a)$. We
have thus proven
the identity:
\begin{equation}\label{eq:100}
  \wh{E_a}(s) =  \sqrt a (a^{-s} + u_s(a))
\end{equation}
In a similar manner we have:
\begin{equation}
   \int_a^\infty \frac{\phi_a^+(x)+\phi_a^-(x)}2 x^{-s}\,dx
   = \int_a^\infty J_0(2\sqrt{ax})x^{-s}\,dx 
   +\int_0^a \frac{-\phi_a^+(x) + \phi_a^-(x)}2 g_s(x)\,dx
\end{equation}
\begin{equation}
\begin{split}
  \int_0^a \frac{-\phi_a^+(x) + \phi_a^-(x)}2 g_s(x)\,dx
= \int_0^a ((1 -D_a)^{-1} H_a J_0^a)(x) g_s(x)\,dx =
-\int_0^a J_0(2\sqrt{ax})u_s(x)\,dx \\
= v_s(a) - g_s(a) = v_s(a) - \int_a^\infty J_0(2\sqrt{ax})x^{-s}\,dx 
\end{split}
\end{equation}
\begin{equation}\label{eq:101}
  \wh{A_a}(s) + i \wh{B_a}(s) = \wh{\cH(E_a)}(s) = \sqrt{a}
  \int_a^\infty \frac{\phi_a^+(x)+\phi_a^-(x)}2 x^{-s}\,dx =
  \sqrt{a} v_s(a)
\end{equation}
Then, we obtain the reformulation of \eqref{eq:72} as:
\begin{equation}
  D_s X_s = \wh{E_a}(s) E_a -  \wh{\cH(E_a)}(s) \cH E_a
\end{equation}
And, noting $\wh{D_s X_s}(z) = (s+z-1)\wh{X_s}(z) =
  (s+z-1)\int_a^\infty X_s(x)X_z(x)\,dx$ we are finally led
  to the remarkable result:
\begin{equation}\label{eq:75}
  X_a(s,z) = \int_a^\infty X_s^a(x)X_z^a(x)\,dx =
    \frac{\wh{E_a}(s) \wh{E_a}(z) - 
  \wh{\cH(E_a)}(s) \wh{\cH E_a}(z)}{s + z -1}
\end{equation}
This equation has been proven under the assumption
  $\Re(s)>1$, and $\Re(z)>\frac12$. To complete the
  discussion we need to know that the evaluators $f\mapsto \wh
  f(s)$, $s\in \CC$ are indeed continuous linear forms on
  $K_a$. For $\Re(s)>\frac12$, we have $\wh f(s) =
  \int_a^\infty f(x) x^{-s}\,dx$. For $\Re(s)<\frac12$ we
  have $\wh f(s) =
  \frac{\Gamma(1-s)}{\Gamma(s)}\wh{\cH(f)}(1-s)$. For
  $\Re(s)=\frac12$ continuity follows by the
  Banach-Steinhaus theorem, and of course more elementary
  proofs exist (as in \cite{cras2001} for the cosine or sine
  transform). So we do have unique Hilbert space vectors
  $X_s^a\in K_a$ such
  that $\forall f\in K_a\forall s\in \CC\ \wh f(s) =
  \int_a^\infty X_s^a(x)f(x)\,dx$. Then \eqref{eq:75} holds
  throughout $\CC\times\CC$ by analytic continuation. 

The vectors $X_s^a$ are zero for $s\in-\NN$, and it is more
  precise to use vectors $\cX_s^a = \Gamma(s)X_s^a$ which are
  non-zero for all $s\in\CC$. These vectors are the
  evaluators\footnote{evaluators for the ``euclidean''
  product $\int fg \,dx$, not the ``hilbertian'' $\int
  f\overline g\,dx$.}  for $f\mapsto \cF(s)$, $\cF(s) =
  \Gamma(s)\wh f(s)$. We recapitulate some of the results in
  the following theorem, whose analog for the cosine (or
  sine) transform was given in \cite{cras2002} (up to
  changes of variables and notations, the first paragraph as
  well as equation \eqref{eq:76} are theorems from
  \cite{bra64}; the equations \eqref{eq:53aa}, \eqref{eq:53bb},
  \eqref{eq:52bis} are our contributions. In this specific
  case of $\cH$ we shall later identify exactly $\phi_a^+$
  and $\phi_a^-$ and $E_a$ and $\cE_a$. As we shall explain
   the analog of the $\cE_a$-function in \cite{bra64}
  has value $1$ at $s=\frac12$, and is not identical with the $\cE_a$ here):
\begin{theorem}
  For a given $a>0$ let $K_a$ be the Hilbert space of square
  integrable functions $f(x)$ on $[a,+\infty)$ whose
  $\cH$-transforms $\int_0^\infty J_0(2\sqrt{xy})f(y)\,dy$
  (in the $L^2$-sense) again vanish for $0<x<a$. The
  completed right Mellin transforms $\Gamma(s)\wh{f}(s) =
  \Gamma(s)\int_a^\infty f(x)x^{-s}\,dx$ are entire
  functions and evaluations at $s\in\CC$ are continuous linear
  forms.\\
Let $\cX_s^a$ for each $s\in\CC$ be the unique vector in $K_a$ such
  that $\forall f\in K_a\; \Gamma(s)\wh{f}(s) =
  \int_a^\infty f(x)\cX_s^a(x)\,dx$. Let $\phi_a^+$ and
  $\phi_a^-$ be the entire functions 
  which are the solutions to:
  \begin{align}\label{eq:53aa}
    \phi_a^+(x)  + \int_0^a J_0(2\sqrt{xy})\phi_a^+(y)\,dy
  &= J_0(2\sqrt{ax})\\  \label{eq:53bb}
    \phi_a^-(x)  - \int_0^a J_0(2\sqrt{xy})\phi_a^-(y)\,dy
  &= J_0(2\sqrt{ax})
  \end{align}
Then 
\begin{equation}
  \label{eq:52bis}
      \wh{E_a}(s) = \sqrt a  \Big(a^{-s} + \frac12\int_a^\infty
    (\phi_a^+(x) - \phi_a^-(x))x^{-s}\,dx \Big)
\end{equation}
is an entire function with trivial zeros at $-\NN$ and, defining
  $\cE_a(s) = \Gamma(s)\wh{E_a}(s)$, we have:
  \begin{equation}
    \label{eq:76}
    \forall s,z\in\CC\quad  \int_a^\infty
  \cX_s^a(x)\cX_z^a(x)\,dx = \frac{\cE_a(s) \cE_a(z) -
  \cE_a(1-s) \cE_a(1-z)}{s + z -1}
  \end{equation}
\end{theorem}

We knew in advance that we had to end up with a
  formula such as \eqref{eq:76} (with a $\cE$ function to be
  discovered\footnote{the method was initially developed by
  the author for the cosine and sine transforms
  \cite{cras2002,cras2003} and leads for them to the only
  known ``explicit'' formulas for $\cE$; for the zero order
  Hankel transform the problem of computing the reproducing
  kernel had been already solved by de~Branges
  \cite{bra64}.}), and this is why we started investigating
  $(x\frac{d}{dx} + s)X_s^a(x)$ in the first place! The
  reason is this: the Hilbert space of the entire functions
  $\cF(s)$, $f\in K_a$ (the Hilbert structure is the one
  from $K_a$, or $(\cF_1,\cF_2) = \frac1{2\pi}\int_{\Re(s)=
  \frac12}
  \cF_1(s)\overline{\cF_2(s)}\frac{|ds|}{|\Gamma(s)|^2}$)
  verifies the de~Branges axioms \cite{Bra}, up to the
  change of variable $s = \frac12 - iz$. Let us recall the
  axioms of \cite{Bra} for a (non-zero) Hilbert space of
  entire functions $F(z)$:
  \begin{enumerate}
  \item[(H1)] for each $z$, evalution at $z$ is a continuous
  linear form,
\item[(H2)] for each $F$,
  $z\mapsto\overline{F(\overline z)}$ belongs 
  to the Hilbert space and has the same norm as $F$,
\item[(H3)] if $F(w) = 0$ then $G(z) =
\frac{z-\overline w}{z - w}
  F(z)$ belongs to the space and has the same norm as $F$
  \end{enumerate}
Let $K(z,w)$ be defined as the evaluator at $z$: $\forall
  F\; F(z) = (F, K(z,\cdot))$. It is anti-analytic in $z$
  and analytic in $w$ (the scalar product is complex linear
  in its first entry, and conjugate linear in its second
  entry). It is a reproducing kernel: $K(z,w) =
  (K(z,\cdot), K(w,\cdot))$. It is proven in \cite{Bra}
that (H1), (H2), (H3)
  entail the existence 
  of an entire function $E(z)$ with $|E(z)|>|E(\overline z)|$ for
  $\Im(z)>0$, such that the space is exactly the set of
    entire functions $F(z)$ such that both
    $\frac{F(z)}{E(z)}$  and $\frac{\overline{F(\overline
    z)}}{E(z)}$ belong to $\HH^2(\Im(z)>0)$, and the Hilbert
    space norm of $F$ is $\frac1{2\pi}\int_\RR
    |F(t)|^2\frac{dt}{|E(t)|^2}$.\footnote{the conditions on
  $F(z)$ are not formulated in \cite{Bra} as Hardy space
  conditions, but they are exactly equivalent.} We
    have incorporated a $2\pi$ for easier comparison with
    our conventions. Then the reproducing kernel is
    expressed as:
  \begin{equation}\label{eq:80}
    K(z,w) = \frac{\overline{E(z)}E(w) - E(\overline
    z)\overline{E(\overline w)}}{i(\overline z -w)}
  \end{equation}
The function $E$ is not unique; if the space has the
isometric
    symmetry $F(z)\mapsto F(-z)$, a function $E$ exists
    which is real on the imaginary axis and writing $E = A -
    i B$ where $A$ and $B$ are real on the real axis, the
    pair $(A,B)$ is unique up to $A \mapsto kA$, $B\mapsto
    k^{-1}B$, $A$ is even and $B$ is odd. If $A(0)\neq0$
    (this happens exactly when the space contains at least
    one element not vanishing at $0$)
    then it may be uniquely normalized so that $A(0)=1$. Then $E$ is
    uniquely determined. 

Model examples are the Paley-Wiener spaces of entire
functions $F(z)$ of exponential type at most $\tau$ with
$\|F\|^2 = \frac1{2\pi} \int_\RR |F(t)|^2\,dt<\infty$. Then
$E(z) = e^{-i\tau z}$ is a possible $E$ function. The
Paley-Wiener spaces are related to the study of the
differential operator $-\frac{d^2}{dx^2}$ on the positive
half-line, and an important class of spaces verifying the
axioms of \cite{Bra} is associated with the theory of the
eigenfunction expansions for Schrödinger operators
$-\frac{d^2}{dx^2} + V(x)$ (\cite{remling}). In these
examples the spaces are indexed by a parameter $\tau$ (the
Schrödinger operator is first studied on a finite interval
$(0,\tau)$) and they are ordered by isometric inclusions
(the $E$-function of a bigger space may be used in the
computation of the norm of an element of a smaller
space). Typically indeed, de~Branges spaces are studied
included in one fixed space $L^2(\RR,\frac{1}{2\pi}d\nu)$,
are ordered by isometric inclusion  and indexed by a
parameter\footnote{the axioms allow for ``jumps'' in the
isometric chain of inclusions, as occur in the theory of the
Krein strings \cite{dymkean}, discrete Schrödinger
equations being special cases.}. Obviously this theory is
intimately related with the Weyl-Stone-Titchmarsh-Kodaira
theory of the spectral measure. The articles of Dym
\cite{dym} and Remling \cite{remling}, the book of Dym and
McKean \cite{dymkean}, will be useful to the interested
reader. In the case of the study of $\cH$ we will have
$d\nu(\gamma) = |\Gamma(\frac12+i\gamma)|^{-2}d\gamma$. It
is an important flexibility of the axioms not to be limited
to functions of finite exponential type, and also the
spectral measures are not necessarily such that
$(1+\gamma^2)^{-1}$ is integrable. It has turned out in our
study of the spaces associated with the $\cH$-transform that
the naturally occurring $E$ function is not the one
normalized to take value $1$ at $z=0$.  Rather the
normalization will prove to be $\lim_{\sigma\to+\infty}
\frac{-iB(i\sigma)}{A(i\sigma)} = 1$. This has an important
impact on the aspect of the differential equations which
will govern the deformation of the $K_a$'s with respect to
$a$: they will take the form of a first order linear
differential system in canonical form (as generally studied
in \cite[\textsection3]{levsarg}.)

The space of the functions $\cF(s) = \Gamma(s)\wh f(s)$,
$f\in K_a$ 
  verify (this is easy) the de~Branges axioms, with $s = \frac12 -
    iz$ and they were defined\footnote{in the variable $z$, and
    associated with the Hankel transform of order zero,
    rather than with the $\cH$ transform.} in
    \cite{bra64}. The spaces $K_a$ have the real structure,
    which is manifest in the $s$ variable through the isometry
  $\cF(s)\mapsto \overline{\cF(\overline s)}$. Rather than with
  the reproducing kernel $K(z_1,z_2)$ we work mainly with
    $\cX(s_1,s_2) = K(-\overline{z_1},z_2)$ which is
    analytic in both variables. Of course then it is
    $\cX(\overline s,s)$ which gives the squared norm of the
    evaluator at $s$. Writing $\cE(s) = E(z)$ we obtain
    from \eqref{eq:80}:
    \begin{equation}
      \label{eq:81}
    \cX(s_1,s_2) = \frac{\cE(s_1)\cE(s_2) - \cE(1-s_1)\cE(1-s_2)}{s_1+s_2 -1}
    \end{equation}
which is indeed what has appeared on the right hand side of
    \eqref{eq:76}. With $\cE(s) = \cA(s) - i \cB(s)$, $\cA$
  (resp. $\cB$) even (resp. odd) under $s\mapsto 1-s$, this
    is also:
    \begin{equation}
      \label{eq:82}
    \cX(s_1,s_2) = 2\frac{-i\cB(s_1)\cA(s_2) +
    \cA(s_1)(-i\cB(s_2))}{s_1+s_2-1}
 \end{equation}
and for $\Re(s)\neq\frac12$, $0 < \cX(\overline s,s) =
    2\frac{\Im({\cB(s)}\overline{\cA(s)})}{\Re(s)
    - \frac12}$ so both $\cA$ and $\cB$ have  all their zeros on the
    critical line.\footnote{this is also seen from $2\cA(s) = \cE(s)+\cE(1-s)$
  as $|\cE(s)|>|\cE(1-s)|$ for $\Re(s)>\frac12$. As
  $\cX(\overline s, s) = \frac{|\cE(s)|^2 -
  |\cE(1-s)|^2}{2\Re(s) - 1}$ this is in fact the same
  argument.}  

The method in this chapter has been developed in
\cite{cras2002, copoisson, cras2003} for the case of the
cosine and sine transforms, and it leads to the currently
only known ``explicit'' formulae\footnote{as ``explicit'' as
the Fredholm determinants of the finite Dirichlet kernels
are ``explicit''.} for the structural elements $\cE$, $\cA$,
$\cB$ and reproducing kernels for the spaces for the cosine
and sine transforms.  So far, almost nothing very specific
to $\cH$ has been used apart from it being self-adjoint
self-reciprocal with an entire multiplicative kernel
$k(xy)$. The next section is still of a very general
validity.

As was mentioned in the Introduction the realization of the
  structural elements of the spaces as right Mellin
  transforms of distributions is a characteristic aspect of
  the method; the Dirac delta's in the expressions for
  $A_a(x)$ and $-iB_a(x)$ could have been overlooked if we
  had only been prepared to use functions, and the whole
  development was based on the computation of
  $(x\frac{d}{dx} + s)X_s(x)$ as a distribution. This aspect
  will be further reinforced in the concluding chapter of
  the paper (section \ref{sec:9}) where it will be seen that
  the distributions $A_a(x)$ and $-i B_a(x)$ are very
  naturally differences of boundary values of analytic
  functions, so they are hyperfunctions \cite{morimoto} in a
  natural manner.

Let us consider the behavior of $\wh{A_a}(s)$,
    $\wh{B_a}(s)$, $\wh{E_a}(s)$ and $\wh{\cH(E_a)}(s)$ for
    $\Re(s)\geq\frac12$. Let us first look at $\wh{E_a}(s) =
    \sqrt a \Big(a^{-s} + \frac12\int_a^\infty (\phi_a^+(x)
    - \phi_a^-(x))x^{-s}\,dx \Big)$. We remark that
    $\phi_a^+(x) - \phi_a^-(x)$ is the $\cH$-transform of
    $-(\phi_a^+(x)+\phi_a^-(x))\Un_{0<x<a}(x)$. 


    \begin{lemma}\label{lem:Ok}
      Let $k(x)$ a continuous function on $[0,+\infty)$ and
      $A\in[0,1]$ be such that $k_1(x) = \int_0^x k(t)\,dt =
      O(x^A)$ as $x\to\infty$. Let $a>0$ and let $f(x)$ be
      an absolutely continuous function on $[0,a]$.  Then
      $\int_0^a k(xy)f(y)\,dy = O(x^{A-1})$ as
      $x\to+\infty$.
    \end{lemma}

There exists $C<\infty$ such that $\forall x>0\;|k_1(x)|\leq
C\,x^A$. Then $\int_0^a k(xy)f(y)\,dy = \frac1x k_1(xa)f(a)
- \frac1x \int_0^a k_1(xy)f'(y)\,dy$, and $|\int_0^a
k_1(xy)f'(y)\,dy|\leq C x^A \int_0^a y^A |f'(y)|\,dy$. This
was easy\dots

With $k(x) = J_0(2\sqrt x)$, one has $k_1(x) =
\sqrt{x}J_1(2\sqrt x) = O(x^{\frac14})$. We have
$\phi_a^+(x) - \phi_a^-(x) = -\int_0^a
J_0(2\sqrt{xy})(\phi_a^+(x) + \phi_a^-(x))\,dy$ and from
Lemma \ref{lem:Ok} this is $O(x^{-\frac34})$. So the
integral in the expression for $\wh{E_a}(s)$ is absolutely
convergent for $\Re(s)>\frac14$. In particular $\wh{E_a}$ is
bounded on
    the critical line. But then $\wh{\cH(E_a)}(s) =
    \chi(s)\wh{E_a}(1-s)$ is also bounded. Hence:
    \begin{proposition}\label{lem:o1}
      The functions $\wh{A_a}$ and $\wh{B_a}$ are bounded on
    the critical line.
    \end{proposition}
 
Let us turn to the situation regarding $\Re(s) = \sigma\to+\infty$.

Let $f(x)$ be a function of class $C^2$ on $[0,a]$ and $e(x) = \int_0^a
J_0(2\sqrt{xy})f(y)\,dy$. It is $O(x^{-\frac34})$. There
holds $\frac{d}{dx} x e(x) = \int_0^a (\frac{d}{dy} y
J_0(2\sqrt{xy}))f(y)\,dy = af(a)J_0(2\sqrt{ax}) - \int_0^a
J_0(2\sqrt{xy})\,yf'(y)\,dy$. Let $k(x) = \int_0^a
J_0(2\sqrt{xy})\,yf'(y)\,dy$. By the Lemma \ref{lem:Ok} it
is $O(x^{-\frac34})$. For $\Re(s)>\frac34$, with absolutely
convergent integrals:
\begin{equation}
  \label{eq:137}
  af(a)\int_a^\infty J_0(2\sqrt{ax})x^{-s}\,dx -
  \int_a^\infty k(x)x^{-s}\,dx = - ae(a)a^{-s} + s
  \int_a^\infty e(x)x^{-s}\,dx
\end{equation}
We show that the left hand side of \eqref{eq:137} is
$O(a^{-s}\frac1s)$ for $\Re(s)>\frac54$. We apply to $k$
what we did for $e$, $\frac{d}{dx} x k(x) =
a^2f'(a)J_0(2\sqrt{ax}) - \int_0^a
J_0(2\sqrt{xy})\,y(yf')'(y)\,dy$. This is $O(1)$ (using
$|J_0|\leq1$). So for $\Re(s)>1$, we can compute
$\int_0^\infty (\frac{d}{dx} x k(x))x^{-s}\,dx$ by
integration by parts, this gives $-ak(a)a^{-s} +
s\int_a^\infty k(x)x^{-s}\,dx$. So for $\Re(s)\geq
1+\epsilon$ we have $\int_a^\infty k(x)x^{-s}\,dx =
O(a^{-s}\frac1s)$. Then regarding $\int_a^\infty
J_0(2\sqrt{ax})x^{-s}\,dx$ we note that $ \frac{d}{dx} x
J_0(2\sqrt{ax}) = J_0(2\sqrt{ax}) -
\sqrt{ax}J_1(2\sqrt{ax})$, so for
$\Re(s)\geq\frac54+\epsilon$ we can apply the same method of
integration by parts, and prove that $\int_a^\infty
J_0(2\sqrt{ax})x^{-s}\,dx = O(a^{-s}\frac1s)$. So the left
hand side of \eqref{eq:137} is indeed $O(a^{-s}\frac1s)$ for
$\Re(s)\geq\frac32$ and we have:
\begin{lemma}
  Let $f(x)$ be a function of class $C^2$ on $[0,a]$ and let
$e(x) = \int_0^a J_0(2\sqrt{xy})f(y)\,dy$. One has
\begin{equation}
  \label{eq:138}
  \int_a^\infty e(x)x^{-s}\,dx =  \frac{a^{-s}}{s}\big(
  ae(a) + O(\frac1s) \big)\qquad (\Re(s)\geq\frac32)
\end{equation}
\end{lemma}
Let us return to
$\int_a^\infty J_0(2\sqrt{ax})x^{-s}\,dx =
\frac1s\big(J_0(2a)a^{1-s}+\int_a^\infty (J_0(2\sqrt{ax}) -
\sqrt{ax}J_1(2\sqrt{ax}))x^{-s}\,dx\big)$. We want to iterate so
we also need $x\frac{d}{dx} \sqrt{ax}J_1(2\sqrt{ax})= \frac14
2\sqrt{ax}\frac{d}{d2\sqrt{ax}} 2\sqrt{ax}J_1(2\sqrt{ax}) =
ax J_0(2\sqrt{ax})$. So we can integrate by parts and obtain
that the last Mellin integral is $O(a^{-s}\frac1s)$ for
$\Re(s)\geq\frac74+\epsilon$. So, certainly:
\begin{equation}
  \label{eq:141}
  \int_a^\infty J_0(2\sqrt{ax}) x^{-s}\,dx =  \frac{a^{-s}}{s}\big(
  a J_0(2a) + O(\frac1s) \big)\qquad (\Re(s)\geq\frac52)
\end{equation}

Using   $\phi_a^+ = J_0^a - \cH P_a \phi_a^+$
and $\phi_a^- = J_0^a + \cH P_a \phi_a^-$ and combining
\eqref{eq:138} and \eqref{eq:141} we obtain:

    \begin{subequations}
    \begin{proposition}\label{prop:10}
      One has for $\Re(s)\geq\sigma_0$ (here $\sigma_0 =
      \frac52$ for example):
      \begin{align} \label{eq:139}
        \wh{E_a}(s) &= a^{\frac12 -s} (1 + \frac{a\phi^+(a)
    - a\phi^-(a)}{2s} + O(\frac1{s^2}))        & 
\wh{A_a}(s) &= \frac{\sqrt a}2 a^{-s} (1 +
    \frac{a\phi_a^+(a)}s + O(\frac1{s^2}))         \\
        \wh{\cH(E_a)}(s) &= a^{\frac12 -s} (\frac{a\phi^+(a)
    + a\phi^-(a)}{2s} + O(\frac1{s^2}))   &
 -i\wh{B_a}(s) &= \frac{\sqrt a}2 a^{-s} (1 -
    \frac{a\phi_a^-(a)}s + O(\frac1{s^2})) 
      \end{align}
    \end{proposition}      
    \end{subequations}

    \begin{theorem}\label{thm:11}
      One has 
      \begin{equation}
        \label{eq:70}
        \lim_{\sigma\to+\infty}
    \frac{-i\cB_a(\sigma)}{\cA_a(\sigma)} = 1 \qquad\text{and}\qquad
    \frac{\cE_a(1-\sigma)}{\cE_a(\sigma)}\sim_{\sigma\to+\infty} 
    \frac{a\phi_a^+(a)+a\phi_a^-(a)}{2\sigma}
      \end{equation}
    \end{theorem}

So the functions $\cA_a$ and $\cB_a$ are not normalized as is
usually done in \cite{Bra} which is to impose
(when possible) to the E function to have value $1$ at the
origin (which for us is $s=\frac12$; the exact value of
$\cA_a(\frac12)$ will be obtained later.) This difference in
normalization is  related to the realization of the
differential equations governing the deformation of the
spaces $K_a$ as a first order differential system in
``canonical'' form, as in the classical spectral theory of linear
  differential equations (\cite{levsarg, coddlevinson}.)
This
  allows to
realize the self-reciprocal scale reversing operator as a
scattering \cite{cras2003}. 

\section{Fredholm determinants, the first order
  differential system, and scattering}
\label{sec:6}
 
\def\ada{a\frac\partial{\partial a}}
\def\xdx{x\frac\partial{\partial x}}
\def\dxx{\frac\partial{\partial x}x}
\def\xdtx{\delta_x}

Let us return to the defining equations for the entire
  functions $\phi_a^+$ and $\phi_a^-$: 
\begin{subequations}
  \begin{align}\label{eq:80a}
    \phi_a^+ + \cH P_a \phi_a^+ &= J_0^a\\
\label{eq:80b}
    \phi_a^- - \cH P_a \phi_a^- &= J_0^a
  \end{align}
\end{subequations}
Either we read these equations as identities on
$(0,\infty)$,
  or we decide that $\cH P_a\phi_a^\pm$ in fact stands for
  $\int_0^a J_0(2\sqrt{xy})\phi_a^\pm(y)\,dy$, and the
  equation holds for $x\in\CC$; the latter option slightly
  conflicts with our earlier definition of $\cH$ as an
  operator on functions or distributions. But whatever
  choice is made this has no impact on what comes next. We shall
  apply to the equations the operators
  $a\frac\partial{\partial a}$ and $x\frac\partial{\partial
  x}$. As $J_0^a(x) = J_0(2\sqrt{ax})$ we have $\ada J_0^a =
  \xdx J_0^a$. We write $\xdtx = \xdx + \frac12 = \dxx -
  \frac12$. First we have:
  \begin{subequations}
    \begin{align}
      \ada\phi_a^+ + \cH P_a
  \ada\phi_a^+ &= - a\phi_a^+(a) J_0^a + \ada J_0^a \\
      \ada\phi_a^- - \cH P_a
  \ada\phi_a^- &= + a\phi_a^-(a) J_0^a + \ada J_0^a 
    \end{align}
  \end{subequations}
Then, as $\xdx \cH = - \cH\dxx$, $\xdtx \cH = - \cH\xdtx$,
$\xdtx P_a f = (P_a \xdtx
  f) - af(a)\delta_a(x)$, $\xdtx J_0^a = \ada J_0^a +
  \frac12 J_0^a$:
\begin{subequations}
  \begin{align}
    \xdtx \phi_a^+ - \cH P_a \xdtx\phi_a^+ &= (\frac12 - a
  \phi_a^+(a))J_0^a + \ada J_0^a \\
    \xdtx \phi_a^- + \cH P_a \xdtx \phi_a^- &= (\frac12 + a
    \phi_a^-(a))J_0^a+\ada J_0^a
  \end{align}
\end{subequations}
Combining we obtain:
\begin{subequations}
  \begin{align}
   \ada\phi_a^+ - \xdtx \phi_a^-  + \cH P_a  (\ada\phi_a^+ -
  \xdtx \phi_a^- ) = -(a\phi_a^+(a) + a \phi_a^-(a) +
  \frac12) J_0^a \\
   \ada\phi_a^- - \xdtx \phi_a^+  - \cH P_a  (\ada\phi_a^- -
  \xdtx \phi_a^+ ) = +(a\phi_a^+(a) + a \phi_a^-(a) -
  \frac12) J_0^a 
  \end{align}
\end{subequations}
Comparing with \eqref{eq:80a} and \eqref{eq:80b}, and as
  there is uniqueness:
\begin{subequations}
  \begin{align}\label{eq:110a}
    \ada\phi_a^+ - \xdtx \phi_a^-  = -(a\phi_a^+(a) + a \phi_a^-(a) +
  \frac12) \phi_a^+\\ \label{eq:110b}
    \ada\phi_a^- - \xdtx \phi_a^+  = + (a\phi_a^+(a) + a \phi_a^-(a) -
  \frac12)\phi_a^-
  \end{align}
\end{subequations}
The quantity $a\phi_a^+(a) + a \phi_a^-(a)$ will play a
  fundamental rôle and we shall denote it by
  $\mu(a)$.\footnote{maybe it would be unfair to hide the
  fact that $\mu(a)=2a$, in
  this study of $\cH$! In a later section a further mu
  function, associated with a variant of $\cH$, will also be
  found explicitely and it will be quite more complicated.} So:
\begin{subequations}
  \begin{align}\label{eq:90a}
    (\ada+\frac12 + \mu(a))\phi_a^+ = \xdtx \phi_a^-\\
\label{eq:90b}
    (\ada+\frac12 - \mu(a))\phi_a^- = \xdtx \phi_a^+ 
  \end{align}
\end{subequations}
It follows easily from this that
  $a\frac\partial{\partial a}(\phi_a^+\phi_a^-) = -
  \phi_a^+\phi_a^- + 
  \frac12\dxx((\phi_a^+)^2 + (\phi_a^-)^2)$. So
  \begin{equation}
    a\frac{d}{da}\int_0^a \phi_a^+(x)
  \phi_a^-(x)\,dx = a\phi_a^+(a)\phi_a^-(a) - \int_0^a \phi_a^+(x)
  \phi_a^-(x)\,dx + \frac12 a(\phi_a^+(a)^2 + \phi_a^-(a)^2)
  \end{equation}
  \begin{equation}
    \label{eq:43}
    a\frac{d}{da}\,a\,\int_0^a \phi_a^+(x)\phi_a^-(x)\,dx = \frac12\mu(a)^2
  \end{equation}
We then compute:
\begin{equation}
  \label{eq:44}
  \int_0^a \phi_a^+(x)\phi_a^-(x)\,dx = \int_0^a ((1 -
  D_a)^{-1}J_0^a)(x) J_0^a(x)\,dx\;,
\end{equation}
where we recall $\phi_a^+ = (1 + H_a)^{-1} J_0^a$, $\phi_a^-
= (1 - H_a)^{-1} J_0^a$, $D_a = H_a^2$.
The operator $D_a$ acts on $L^2(0,a;dx)$ with kernel
$D_a(x,z) = \int_0^a
J_0(2\sqrt{xy})J_0(2\sqrt{yz})\,dy$. After the change of
variables $x=at$, $y=au$, $z=av$ this becomes the operator
$d_a$ on $L^1(0,1;dt)$ with kernel $d_a(t,v) = \int_0^1
aJ_0(2a\sqrt{tu})\,aJ_0(2a\sqrt{uv})\,du$. We compute the
derivative with respect to $a$:
\def\udu{u\frac\partial{\partial u}}
\def\duu{\frac\partial{\partial u}u}
\begin{subequations}
  \begin{align}
  &\frac{\partial}{\partial a} \int_0^1
aJ_0(2a\sqrt{tu})\,aJ_0(2a\sqrt{uv})\,du \\
&= \int_0^1
( (2\udu +1)J_0(2a\sqrt{tu}))\,aJ_0(2a\sqrt{uv})\,du +
\int_0^1
aJ_0(2a\sqrt{tu}))((2\duu -1) J_0(2a\sqrt{uv}))\,du \\
&= 2aJ_0(2a\sqrt t)J_0(2a\sqrt v)
  \end{align}
\end{subequations}
So $\frac{d}{da} d_a$ is a rank one operator, with range
$\CC J_0(2a\sqrt t)\Un_{0<t<1}(t)$. We now use the
well-known formula
\begin{equation}
  \label{eq:38}
  \frac{d}{da} \log\det(1-d_a) = -
\textrm{Tr}((1-d_a)^{-1}\frac{d}{da} d_a)
\end{equation}
The rank one operator $(1-d_a)^{-1}\frac{d}{da} d_a$ has the
function $(1 - d_a)^{-1} 2aJ_0(2a\sqrt t)$ as eigenvector
and the eigenvalue is $\int_0^1 J_0(2a\sqrt t)((1 -
d_a)^{-1} 2aJ_0(2a\sqrt v))(t)\,dt$. Going back to $(0,a)$
we obtain $2\int_0^a J_0(2\sqrt{ax})((1 - D_a)^{-1}
J_0(2\sqrt{az}))(x)\,dx$ and in view of \eqref{eq:44}
we have proven:
\begin{equation}
  \label{eq:39}
\frac{d}{da} \log\det(1-D_a) = - 2\int_0^a
\phi_a^+(x)\phi_a^-(x)\,dx
\end{equation}
Then, using \eqref{eq:43}, we have the important formula:
\begin{equation}
  \label{eq:47}
  \mu(a)^2 = -a\frac{d}{da}\,a\frac{d}{da}\, \log\det(1 -D_a)
\end{equation}
We shall now also relate $\phi_a^+(a)$ and $\phi_a^-(a)$ to
Fredholm determinants. In fact the following holds:
\begin{subequations}
  \begin{align}\label{eq:141a} 
    a\phi_a^+(a) &= + a\frac{d}{da}\log\det(1 + H_a)\\
\label{eq:141b} a\phi_a^-(a) &= - a\frac{d}{da}\log\det(1 - H_a)
  \end{align}
\end{subequations}
This is the application of a well-known general theorem, for
any continuous kernel $k(x,y)$: if
$w(x)+ \int_0^a k(x,y)w(y)\,dy = k(x,a)$ for $0\leq x \leq
a$ then $w(a) = +\frac{d}{da}\log\det_{(0,a)}(\delta(x-y) +
k(x,y))$. A proof may be given which is of a somewhat
similar kind as the one given above for \eqref{eq:39}, or
one may more directly use the Fredholm's formulas for the
determinant and the resolvent.\footnote{let us recall that
for a continuous kernel on a finite interval, the formula of
Fredholm for a determinant as a convergent series always
applies, even if the operator given by the kernel is not
trace class, which may happen.} The theorem is proven in
the book of P.~Lax \cite{lax}, Theorem 12 of Chapter 24 (Lax
treats the case of a kernel on $(a,+\infty)$, here we have
the simpler case of a finite interval $(0,a)$.) This
means that $\mu(a)$ has another expression in terms of
Fredholm determinants:
\begin{equation}
  \label{eq:49}
  \mu(a) = a\frac{d}{da}\log\frac{\det(1 + H_a)}{\det(1 - H_a)}
\end{equation}
Combining \eqref{eq:47} and \eqref{eq:49} we obtain:
\begin{subequations}
  \begin{align}\label{eq:300a}
    -2a\frac{d}{da}\,a\frac{d}{da}\log\det(1+H_a) &=
     \left(a\frac{d}{da}\log\frac{\det(1 + H_a)}{\det(1 -
     H_a)}\right)^2 -
     a\frac{d}{da}a\frac{d}{da}\log\frac{\det(1 +
     H_a)}{\det(1 - H_a)}\\  \label{eq:300b}  
-2a\frac{d}{da}\,a\frac{d}{da}\log\det(1-H_a) &=
     \left(a\frac{d}{da}\log\frac{\det(1 + H_a)}{\det(1 -
     H_a)}\right)^2 +
     a\frac{d}{da}a\frac{d}{da}\log\frac{\det(1 +
     H_a)}{\det(1 - H_a)}\\
     2 a\frac{d}{da} a\phi_a^+(a) &= -\mu(a)^2 + a\mu'(a)\\
     2 a\frac{d}{da} a\phi_a^-(a) &= +\mu(a)^2 + a\mu'(a)\\
\frac{d}{da} a (\phi_a^-(a) - \phi_a^+(a)) &= a (\phi_a^+(a)
     + \phi_a^-(a))^2
  \end{align}
\end{subequations}
These Fredholm determinants identities are reminiscent of
certain well-known Gaudin identities \cite[App. A16]{mehta},
which apply to the even and odd parts of an additive
(Toeplitz) convolution kernel on an interval $(-a,a)$; here
the situation is with kernels $k(xy)$ which have a
multiplicative look, and reduction to the additive case
would give $g(t+u)$ type kernels on semi-infinite intervals.

We have defined  $A_a = \frac{\sqrt a}2(\phi_a^+ +
  \cH\phi_a^+)$ and $-iB_a = \frac{\sqrt a}2(-\phi_a^-+
  \cH\phi_a^-)$. Let us recall that here $\phi_a^\pm$ is
  restricted to $[0,+\infty)$ and is then tempered as a
  distribution. Using the differential equations
\eqref{eq:90a} and \eqref{eq:90b} and the commutation
  property $\xdtx\cH = - \cH\xdtx$, $\xdtx = \xdx+\frac12$,
  we have $\xdtx A_a = \frac{\sqrt a}2(\xdtx\phi_a^+ -
  \cH\xdtx\phi_a^+) = \frac{\sqrt a}2(\ada + \frac12 -
  \mu(a))(\phi_a^- - \cH\phi_a^-) = -(\ada - \mu(a))(-iB_a)$
  and $\xdtx (-iB_a) = \frac{\sqrt a}2(-(\ada + \frac12 +
  \mu(a))\phi_a^+ - (\ada + \frac12 + \mu(a))\cH\phi_a^+) =
  -(\ada+\mu(a))A_a$. The following first order  system of
  differential equations therefore holds:
  \begin{subequations}
    \begin{align}\label{eq:91a}
      \ada A_a &= -\mu(a)A_a - \xdtx(-iB_a)\\
\label{eq:91b}
      \ada(-iB_a) &= +\mu(a)(-iB_a) - \xdtx A_a
    \end{align}
  \end{subequations}
Then we also have the second order differential equations
  $(\ada - \mu(a))(\ada + \mu(a))A_a = +\xdtx^2 A_a$ and
  $(\ada + \mu(a))(\ada - \mu(a))(-iB_a) = + \xdtx^2
  (-iB_a)$, or, taking the right Mellin transforms, and writing
  $s= \frac12 + i \gamma$, $\xdtx = i\gamma$:
\begin{subequations}
  \begin{align}
    -\ada\ada \wh{A_a} + (\mu(a)^2 - a\mu'(a)) \wh{A_a} &= \gamma^2
  \wh{A_a}\\
    -\ada\ada (-i\wh{B_a}) + (\mu(a)^2 + a\mu'(a)) (-i\wh{B_a}) &= \gamma^2
  (-i\wh{B_a})
\end{align}
\end{subequations}
With the new variable $u=\log(a)$ we
  obtain Dirac and Schrödinger equations which are
  associated with this study of $\cH$, modeled on the study
  of the cosine and sine transforms summarized in
  \cite{cras2002,cras2003}. All quantities in the statement of
  the theorem will be completely explicited later in terms
  of Bessel functions, but we keep the notation sufficiently
  general to allow, if an interesting other case arises, to
  write down the identical results:

\begin{subequations}
  \begin{theorem}
    For each $a>0$ let $\phi_a^+$ and $\phi_a^-$ be the entire functions
  which are the solutions to:
  \begin{align}
    \phi_a^+(x)  + \int_0^a J_0(2\sqrt{xy})\phi_a^+(y)\,dy
  &= J_0(2\sqrt{ax})\\
    \phi_a^-(x)  - \int_0^a J_0(2\sqrt{xy})\phi_a^-(y)\,dy
  &= J_0(2\sqrt{ax})
  \end{align}
Let $H_a$ be
  the integral operator on $L^2(0,a;dx)$ with kernel
  $J_0(2\sqrt{xy})$. There holds:
  \begin{align}
    \phi_a^+(a) &= +\frac{d}{da}\log\det(1+H_a)\\
    \phi_a^-(a) &= -\frac{d}{da}\log\det(1-H_a)\;.
  \end{align}
The tempered distributions $A_a = \frac{\sqrt
  a}2(1+\cH)(\phi_a^+\Un_{0<x<\infty})$ and $B_a =
  i\frac{\sqrt a}2 (-1+\cH)(\phi_a^-\Un_{0<x<\infty})$
  vanish on $(-\infty,a)$ and are respectively
  self-reciprocal and skew-reciprocal under $\cH$. Their
  completed right Mellin transforms $\cA_a(s) =
  \Gamma(s)\wh{A_a}(s)$ and $\cB_a(s) =
  \Gamma(s)\wh{B_a}(s)$ are entire functions with all their
  zeros on the critical line, they are respectively even and
  odd for $s\leftrightarrow 1-s$, and they verify the following
  Dirac and Schrödinger types of differential equations in
  the variable $u = \log(a)$, $-\infty<u<+\infty$,
  \begin{align}
\frac{d}{du} \cA_a &= - \mu(a) \cA_a -\gamma \cB_a\\
\frac{d}{du} \cB_a &= + \mu(a) \cB_a + \gamma \cA_a\\
\label{eq:95a}
  \gamma^2\cA_a &= \left(-\frac{d^2}{du^2}   + V_+(u) \right) \cA_a  \\ 
\label{eq:95b}
 \gamma^2\cB_a &= \left(-\frac{d^2}{du^2} + V_-(u)
 \right) \cB_a \\
V_+(\log a) &= \mu(a)^2 - \frac{d\,\mu(a)}{du} =
  -2\frac{d^2\log\det(1+H_a)}{du^2}\\ 
V_-(\log a) &= \mu(a)^2 + \frac{d\,\mu(a)}{du} = - 2\frac{d^2
  \log\det(1-H_a)}{du^2}\\
\mu(a) &= \frac{d}{du} \log\frac{\det(1+H_a)}{\det(1-H_a)} =
  a\phi_a^+(a) + a\phi_a^-(a)
  \end{align}
where  $s = \frac12 + i\gamma$.
  \end{theorem}
\end{subequations}

Let us consider the Hilbert space of pairs $\begin{bmatrix} \alpha(u)\\
    \beta(u) \end{bmatrix}$ on $\RR$ with squared norms
    $\int_{-\infty}^\infty |\alpha(u)|^2
    + |\beta(u)|^2 \frac{du}2$, and the two equivalent
      differential systems in canonical forms:
  \begin{equation}
    \label{eq:36}
    \left(\begin{bmatrix}
      0&1\\-1&0\end{bmatrix}\frac{d}{du} - 
    \begin{bmatrix} 0&\mu(e^u)\\\mu(e^u)&0\end{bmatrix}\right)
      \begin{bmatrix} \alpha(u)\\  \beta(u) \end{bmatrix} 
= \gamma \begin{bmatrix} \alpha(u)\\ \beta(u) \end{bmatrix}
  \end{equation}
  \begin{equation}
    \label{eq:37}
        \left(\begin{bmatrix}
      0&1\\-1&0\end{bmatrix}\frac{d}{du} + 
    \begin{bmatrix} -\mu(e^u)&0\\0&\mu(e^u)\end{bmatrix}\right)
      \begin{bmatrix} \alpha(u)+\beta(u)\\  -\alpha(u)+\beta(u) \end{bmatrix} 
= \gamma \begin{bmatrix} \alpha(u)+\beta(u)\\
-\alpha(u)+\beta(u) \end{bmatrix}
  \end{equation}
The components obey the corresponding Schrödinger equations:
\begin{subequations}
  \begin{align}
\label{eq:400a}
  -\alpha''(u) +  V_+(u) \alpha(u) &= \gamma^2 \alpha
  \qquad\qquad V_+(u) = \mu(e^u)^2 - \frac{d\,\mu(e^u)}{du}\\ 
\label{eq:400b}
  -\beta''(u) +  V_-(u) \beta(u) &= \gamma^2 \beta
  \qquad\qquad V_-(u) = \mu(e^u)^2 + \frac{d\,\mu(e^u)}{du}
  \end{align}
\end{subequations}
Regarding the behavior at $-\infty$, we are in the
  limit-point case for each of the Schrödinger equations
  \eqref{eq:400a} and \eqref{eq:400b} because clearly (say,
  from the defining integral equations for $\phi_a^+$ and
  $\phi_a^-$) one has $\phi_a^+(a)\to_{a\to0}J_0(0) = 1$,
  $\phi_a^-(a)\to_{a\to0}1$, $\mu(a)\sim_{a\to0}2a$, so the
  potentials are exponentially vanishing as
  $u\to-\infty$. Perhaps we should reveal that one has
  exactly $\mu(a) = 2a = 2e^u$ so we are dealing here with
  quite concrete Schrödinger equations and Dirac systems
  whose exact solutions will later be written explicitely in
  terms of modified Bessel functions, but we delay using any
  information which would be too specific of the
  $\cH$-transform.

For each $\gamma\in\CC$
  \begin{equation}
    \label{eq:163}
    u\mapsto \begin{bmatrix} \cA_{\exp(u)}(\frac12+i\gamma)\\
    \cB_{\exp(u)}(\frac12+i\gamma) \end{bmatrix}
  \end{equation} is a (non-zero) solution of the system
\eqref{eq:36}, and we now show that it is square-integrable
(with respect to $du = d\log(a)$) at $+\infty$. Let us recall the equation \eqref{eq:82} $
(s+z-1)\cX_a(s,z) = -2i\cB_a(s)\cA_a(z)
  -2i\cA_a(s)\cB_a(z)$, from which we deduce
  \begin{equation}
    \label{eq:31}
    \ada \cX_a(s,z) = -2\cA_a(s)\cA_a(z) -2 (i\cB_a(s))(i\cB_a(z))
  \end{equation}
We have\footnote{let us recall the notation $\cX_s^a =
  \Gamma(s) X_s^a \in L^2(a,+\infty;dx)$.} $\|\cX_s^a\|^2
  = \cX_a(s,\overline s)$, $\cA_a(\overline s)
  =\overline{\cA_a(s)}$, $i\cB_a(\overline s) =
  \overline{i\cB_a(s)}$, so 
  \begin{equation}
    \ada \|\cX_s^a\|^2 = -2
  |\cA_a(s)|^2 -2 |\cB_a(s)|^2
  \end{equation}
and as of course
  $\lim_{a\to+\infty} \|\cX_s^a\|^2 = 0$ (we have
  $\|\cX_s^a\|^2\leq \int_a^\infty |\cX_s^1|^2(x)\,dx$ for
  $a\geq1$) we obtain:
  \begin{equation}
    \label{eq:35}
     \forall s\in\CC \quad \|\cX_s^a\|^2 = 2 \int_a^\infty
  (|\cA_a(s)|^2 + |\cB_a(s)|^2)\frac{da}a 
  \end{equation}
This establishes the square-integrability at $+\infty$ of
$\left[\begin{smallmatrix} \cA_{\exp(u)}(s)\\
    \cB_{\exp(u)}(s) \end{smallmatrix}\right]$, for any $s\in\CC$.

The solutions
of \eqref{eq:36} with eigenvalue $\gamma=0$ are
$\left[\begin{smallmatrix} \cA_a(\frac12)\\ 0
\end{smallmatrix}\right]$ and $\left[\begin{smallmatrix} 0\\
\cA_a(\frac12)^{-1}
  \end{smallmatrix}\right]$. The former is
square-integrable, so from $2\leq t+t^{-1}$ the latter then
necessarily is not. This confirms that the Dirac system
\eqref{eq:36} is in the limit point case at $+\infty$
(according to a general theorem of Levitan
\cite[\textsection13, Thm 7.1]{levsarg} any first order
differential operator $\left[\begin{smallmatrix}
0&1\\-1&0\end{smallmatrix}\right]\frac{d}{du} +
\left[\begin{smallmatrix}
a(u)&b(u)\\c(u)&d(u)\end{smallmatrix}\right]$ with
continuous coefficients is in the limit point case at
infinity).  So the pair \eqref{eq:163} is in fact, for any
$\gamma\in\CC$, the unique solution of  \eqref{eq:36} which
is square-integrable at $+\infty$. Also the Schrödinger
equation \eqref{eq:400b} is in the limit point case as not
all of its solutions are square integrable at
$+\infty$. Whether the limit-point case at $+\infty$ holds
for equation \eqref{eq:400a} is less evident. Let us recall
from \cite[\textsection9, Thm 2.4]{coddlevinson} and
\cite[\textsection X, Thm X.8]{reedsimon} that a
sufficient condition for this is the existence of a
lowerbound $\liminf_{u\to+\infty} V_+(u)/u^2 > -\infty$. We
will prove in the next chapter that $\mu(a) = 2 a = 2 e^u$
so this is certainly the case here. In the present chapter
only the fact that the Dirac system is known to be in the
limit-point case will be used.

We now take $u_0 = \log a_0$ and apply on $(u_0,\infty)$ the
Weyl-Stone-Titchmarsh-Kodaira theory
(\cite[\textsection9]{coddlevinson},
\cite[\textsection3]{levsarg}). Let $\psi(u,s)$ be the
unique solution of the system \eqref{eq:36} for the
eigenvalue $\gamma$, $s=\frac12+i\gamma$, and with the
initial condition $\psi(u_0,s) = \left[\begin{smallmatrix}
1\\ 0 \end{smallmatrix}\right]$ and let $\phi(u,s)$ be the
unique solution with the initial condition $\phi(u_0,s) =
\left[\begin{smallmatrix} 0\\ 1
\end{smallmatrix}\right]$. Let
$m(\gamma)\psi(u,s)+\phi(u,s)$ for $\Im\gamma>0$ be the
unique solution which is square-integrable on
$(u_0,+\infty)$. So $m(\gamma) =
\frac{\cA_{a_0}(s)}{\cB_{a_0}(s)}$, $s=\frac12+i\gamma$,
$\Re(s)<\frac12$. It is a fundamental general property
of the $m$ function from Hermann Weyl's  theory
that $\Im(m(\gamma))>0$ (for $\Im(\gamma)>0$.) Here, we have
a case where the
  $m$-function is found to be meromorphic on all of $\CC$;
 so we see that its poles and
zeros on $\RR$ are simple. Furthermore, the spectral measure
$\nu$ is obtained via the formula $\nu(a,b) =
\lim_{\epsilon\to 0^+} \frac1\pi
  \int_a^b \Im m(\gamma+i\epsilon)\,d\gamma$ (under the
  condition
$\nu\{a,b\} = 0$). We  obtain:
\begin{equation}
  \label{eq:191}
  d\nu(\gamma) = \sum_{\cB_{a_0}(\rho) = 0}
\frac{\cA_{a_0}(\rho)}{-i\cB_{a_0}'(\rho)}\delta(\gamma-\Im\rho)
\end{equation}
The spectrum is thus purely discrete and the general theory
  tells us further that the finite linear
combinations $\sum_{\rho} c_\rho
\frac{\cA_{a_0}(\rho)}{-i\cB_{a_0}'(\rho)}\psi(u,\rho)$ have
squared norms $\sum_{\rho}
\frac{\cA_{a_0}(\rho)}{-i\cB_{a_0}'(\rho)}|c_\rho|^2$ and
    also that they are
dense in $L^2((u_0,\infty)\to\CC^2;\,du)$.  For
$\cB_{a_0}(\rho) = 0$, $\psi(u,\rho) =
\cA_{a_0}(\rho)^{-1}\left[\begin{smallmatrix}
\cA_{\exp(u)}(\rho)\\
    \cB_{\exp(u)}(\rho) \end{smallmatrix}\right]\Un_{u\geq
u_0}(u)$, so the vectors $Z_\rho^{a_0} =
\left[\begin{smallmatrix} 2\cA_{\exp(u)}(\rho)\\
    2\,\cB_{\exp(u)}(\rho)
    \end{smallmatrix}\right]\Un_{u\geq
u_0}(u)$ are an orthogonal basis of
$L^2((u_0,\infty)\to\CC^2;\, \frac12 du)$ and they satisfy
$\|Z_\rho^{a_0}\|^2 = -
2\,\cA_{a_0}(\rho)\,i\,\cB_{a_0}'(\rho)$.  Similarly
    a
spectral interpretation is given to the zeros of $\cA_{a_0}$
if one looks at the initial condition
$\left[\begin{smallmatrix} 0\\
    1 \end{smallmatrix}\right]$. The factors of $2$ and
$\frac12$, have
    been incorporated so that the statement may be
  translated (taking into account results established
  later) into the fact that the
    evaluators $K_{a_0}(\rho,z)$, for $\cB_{a_0}(\rho) = 0$,
    are an orthogonal basis of the Hilbert space of the
    functions $\Gamma(z)\wh f(z)$, $f\in K_{a_0}$. This last
statement is a
general theorem (under a certain condition) for spaces with
the de~Branges axioms \cite[\textsection22]{Bra}.

To discuss in a  self-contained manner
    the generalized Parseval identity which is associated
    with the differential system on the full line, it is
    convenient to make a preliminary majoration of
    $\|X_s^a\|^2$, $\Re(s) =
  \frac12$. From \eqref{eq:76} we have, for
    $\Re(s)=\frac12$: $\|\cX_s^a\|^2 =
  2\Re(\cE_a(s)\overline{\cE_a'(s)})$. And $\cE_a(s) =
  \Gamma(s)\wh{E_a}(s)$. And $\wh{E_a}(s) = \sqrt a
  \Big(a^{-s} + \frac12\int_a^\infty
    (\phi_a^+(x) - \phi_a^-(x))x^{-s}\,dx \Big)$. We know
    from the discussion of Lemma \ref{lem:o1} that the
    integral in the expression for $\wh{E_a}(s)$ is
    absolutely convergent for $\Re(s)>\frac14$. Hence by the
    Riemann-Lebesgue lemma $\wh{E_a}(\frac12+i\gamma)\sim
    a^{-i\gamma}$ as $|\gamma|\to\infty$,
    $\gamma\in\RR$. Similarly,
    $\wh{E_a}'(\frac12+i\gamma)\sim -\log(a)
    a^{-i\gamma}$. So, with $\|\cX_s^a\|^2 =
    |\Gamma(s)|^2\|X_s^a\|^2$ and using Stirling's formula
    we obtain:
   \begin{lemma}
      For each given $a>0$ one has
    $\|X_s^a\|^2 \sim 2\log|s|$ as $|s|\to\infty$, $\Re(s) =
    \frac12$.
    \end{lemma}

From \eqref{eq:75} expressed using $A_a$ and $B_a$
    we see that $\frac{\wh{B_a}(s)}{s-\frac12}$ is square
    integrable, so $s^{-1} \wh{B_a}(s)$ is square integrable
    on the critical line (with respect to $|ds|$). Then
    using again \eqref{eq:75} we see that $s^{-1}
    \wh{A_a}(s)$ is also square integrable on the critical
    line.\footnote{we know in fact according to proposition
    \ref{prop:10} that $\wh{A_a}$ and $\wh{B_a}$ are bounded on
    the critical line.}
Let us pick a function $F(s)$ on the critical line which is
    such that $sF(s)$ is square integrable. Then
    $F(s)\overline{\wh{A_a}({s})}$ and
    $F(s)\overline{\wh{B_a}({s})}$ are absolutely integrable
    on the critical line and $(\int_{\Re(s)=\frac12}
    |F(s)\overline{\wh{A_a}({s})}|\frac{|ds|}{2\pi})^2 \leq
    C\int_{\Re(s)=\frac12}
    \frac{|\wh{A_a}({s})|^2}{|s|^2}\frac{|ds|}{2\pi}$
    and similarly with $B_a$. If we define $\alpha_F(u) =
    2\int_{\Re(s)=\frac12}
    F(s)\overline{\wh{A_a}({s})}\frac{|ds|}{2\pi}$ and
    $\beta_F(u) = 2\int_{\Re(s)=\frac12}
    F(s)\overline{\wh{B_a}({s})}\frac{|ds|}{2\pi}$ we then
    compute:
 \begin{multline}\label{eq:48}
      \int_{u_0}^\infty |\alpha_F(u)|^2 + |\beta_F(u)|^2\,du
    \leq  C \int_{\Re(s)=\frac12} \frac{\int_{u_0}^\infty
    (|\wh{A_a}(s)|^2 + |\wh{B_a}(s)|^2) \,du}{|s|^2}\;
    \frac{|ds|}{2\pi} \\
=  \frac C2\int_{\Re(s)=\frac12}
    \frac{\|X_s^a\|^2}{|s|^2} \frac{|ds|}{2\pi} < \infty
 \end{multline}
So $\alpha_F(u)$ and $\beta_F(u)$ are square integrable at
$+\infty$. More precisely the above  upper bound holds
as well for $\int_{\Re(s)=\frac12}
    |F(s)\wh{A_a}({s})|\frac{|ds|}{2\pi}$ and $\int_{\Re(s)=\frac12}
    |F(s)\wh{B_a}({s})|\frac{|ds|}{2\pi}$. So the double
    integrals 
    \begin{subequations}
      \begin{align}
        \iint_{u_0 <u<\infty, \Re(s)=\frac12}
    \cA_{\exp(u)}(z)\cA_{\exp(u)}(s)\cF(s)\;\frac{|ds|}{2\pi
    |\Gamma(s)|^2}\;du\\
\iint_{u_0 <u<\infty, \Re(s)=\frac12}
    \cB_{\exp(u)}(z)\cB_{\exp(u)}(s)\cF(s)\;\frac{|ds|}{2\pi|\Gamma(s)|^2}\;du
      \end{align}
    \end{subequations}
where $z\in\CC$ is arbitrary, and $\cF(s) = \Gamma(s)F(s)$,
are absolutely convergent and Fubini may be employed. Using
\eqref{eq:31}:
\begin{equation}
  \label{eq:52}
  \cX_{\exp(u_0)}(z,\overline s) = 2 \int_{u_0}^\infty
  (\cA_{\exp(u)}(z)\overline{\cA_{\exp(u)}(s)} +
  \cB_{\exp(u)}(z)\overline{\cB_{\exp(u)}(s)})\,du
\end{equation}
And we obtain the following identity of absolutely
  convergent integrals, for any $\cF(s) = \Gamma(s)F(s)$
  with $s\, F(s)\in
  L^2(\Re(s)=\frac12;\frac{|ds|}{2\pi})$:
  \begin{equation}
    \label{eq:53}
      \int_{\Re(s)=\frac12} \cX_{\exp(u_0)}(z,\overline
  s)\cF(s)\, \frac{|ds|}{2\pi|\Gamma(s)|^2} = 
  \int_{u_0}^\infty (\cA_{\exp(u)}(z)\alpha_F(u) +
  \cB_{\exp(u)}(z)\beta_F(u))\,du 
  \end{equation}
We shall prove that this identity holds under the weaker
  hypothesis $F(s)\in
  L^2(\Re(s)=\frac12;\frac{|ds|}{2\pi})$. First, still with
  $s\,F(s)$ square integrable we suppose additionally
  that $F = \wh{f}$ with $f\in
    K_{\exp(u_0)}$\footnote{this is certainly possible as
    we know that the $f(x)$ which are  smooth, vanishing
    on $(0,a)$ and of Schwartz decrease
    as $x\to+\infty$ are dense in $K_a$.}.  The
  hilbertian kernel
    $K_{\exp(u_0)}(z,s)$ is $\cX_{\exp(u_0)}(\overline z,s)$
    so $\overline{K_{\exp(u_0)}(z,s)} =
    \cX_{\exp(u_0)}(z,\overline s)$. The equations give
    then:
  \begin{subequations}
   \begin{align}
    \label{eq:57}
    \cF(z) &=  
  \int_{u_0}^\infty (\cA_{\exp(u)}(z)\alpha_F(u) +
  \cB_{\exp(u)}(z)\beta_F(u))\,du \\
\alpha_F(u) &=
    2\int_{\Re(s)=\frac12}
    \cF(s)\cA_a(s)\frac{|ds|}{2\pi\,|\Gamma(s)|^2} \\
    \label{eq:57beta}
\beta_F(u) &=
    2\int_{\Re(s)=\frac12}
    \cF(s)\cB_a(s)\frac{|ds|}{2\pi\,|\Gamma(s)|^2} 
   \end{align}
  \end{subequations}
We have worked under the hypothesis that $sF(s)$ is square
    integrable. To show that the formulae extend in the
    $L^2$ sense, we first examine:
    \begin{equation}
      \label{eq:58}
      |\alpha_F(u)|^2 = 4\int_{\Re(s_1)=\frac12}
    \cF(s_1)\cA_a(s_1)\frac{|ds_1|}{2\pi\,|\Gamma(s_1)|^2}
    \int_{\Re(s_2)=\frac12}
    \overline{\cF(s_2)}\cA_a(s_2)\frac{|ds_2|}{2\pi\,|\Gamma(s_2)|^2}
    \end{equation}
    \begin{equation}
      \label{eq:60}
      \int_{u_0}^\infty |\alpha_F(u)|^2 +
    |\beta_F(u)|^2 \,\frac{du}2=
    \iint_{\Re(s_i)=\frac12}
    \cF(s_1)\overline{\cF(s_2)}\cX_{\exp(u_0)}(s_1,\overline
  s_2)\frac{|ds_1|}{2\pi\,|\Gamma(s_1)|^2}\frac{|ds_2|}{2\pi\,|\Gamma(s_2)|^2}
    \end{equation}
There was absolute convergence in the triple integral used
    as an intermediate. Also $\cX_{\exp(u_0)}(s_1,\overline
    s_2) = \cX_{\exp(u_0)}(s_2,\overline s_1)$ and
    $\int_{\Re(s_1)=\frac12}
    \cF(s_1)\cX_{\exp(u_0)}(s_2,\overline
    s_1)\frac{|ds_1|}{2\pi\,|\Gamma(s_1)|^2} = \cF(s_2)$. Hence:
    \begin{equation}
      \label{eq:61}
      \int_{u_0}^\infty (|\alpha_F(u)|^2 +
    |\beta_F(u)|^2) \,\frac{du}2 = \int_{\Re(s)=\frac12}
    |F(s)|^2\frac{|ds|}{2\pi} =\int_{\exp(u_0)}^\infty |f(x)|^2\,dx
    \end{equation}
So with an arbitrary $f\in K_a$, $F = \wh f$,
    $\cF(s) = \Gamma(s)\wh f(s)$, the assignment $f\mapsto
    (\alpha_F,\beta_F)$ exists in the sense of $L^2$
  convergence when one approximates $f$ by a sequence $f_n$
    in $K_a$ 
  such that $s\wh{f_n}(s)$ is in
  $L^2(\Re(s)=\frac12;\frac{|ds|}{2\pi})$, and $f\mapsto
    (\alpha_F,\beta_F)$ is linear and isometric. We check
    that its range is all of
    $L^2(u_0,\infty;\frac{du}2)\oplus
    L^2(u_0,\infty;\frac{du}2)$. For this let us identify
    the functions $\alpha_w(u)$ and $\beta_w(u)$ which will
    correspond to $\cF(s) = \cX_{a_0}(w,s)$ ($a_0 =
    \exp(u_0)$.)  On one hand from \eqref{eq:53} it must be
    the case that:
    \begin{equation}
      \label{eq:62}
       \forall z\in\CC\quad \int_{\Re(s)=\frac12} \cX_{a_0}(z,\overline
  s)\cX_{a_0}(w,s)\, \frac{|ds|}{2\pi|\Gamma(s)|^2} = 
  \int_{u_0}^\infty (\cA_{\exp(u)}(z)\alpha_w(u) +
  \cB_{\exp(u)}(z)\beta_w(u))\,du      
    \end{equation}
The left hand side is $\cX_{a_0}(w,z)$ which on the other hand
is given by the formula $2\int_{u_0}^\infty
    (\cA_a(z)\cA_a(w)- \cB_a(z)\cB_a(w))\,du$. The functions
    $u\mapsto \left[\begin{smallmatrix} \cA_a(z)\\
    \cB_a(z)\end{smallmatrix}\right]$, $z\in\CC$ are
    certainly dense in $L^2((u_0,\infty)\to
    \CC^2;\frac{du}2)$ as we know in particular that the
    pairs for the $\rho$'s such that $\cB_{a_0}(\rho) = 0$
    give an orthogonal basis. So we have the identification
    on $(u_0,+\infty)$:
    \begin{equation}
      \label{eq:63}
      \alpha_w(u) = 2\,\cA_{\exp(u)}(w)\qquad \beta_w(u) =
      -2\cB_{\exp(u)}(w)
    \end{equation}
This proves that the range is all of $L^2((u_0,\infty)\to
\CC^2;\frac{du}2)$. Let us note that in this identification
the hilbertian evaluator $K_{a_0}(w,\cdot)$ is sent to the
pair $u\mapsto
2\,\Un_{u>u_0}(u)(\overline{\cA_a(w)}\,,\,\overline{\cB_a(w)})$. To
    check if all is coherent we compute the hilbertian
  scalar product $(K_{a_0}(w,\cdot),
  K_{a_0}(z,\cdot))$. We obtain
  $4\int_{u_0}^\infty
  \overline{\cA_a(w)}\cA_a(z)+\overline{\cB_a(w)}{\cB_a(z)}\frac{du}2
= 2\int_{u_0}^\infty
  {\cA_a(\overline w)}\cA_a(z)-{\cB_a(\overline
w)}{\cB_a(z)}\,du = \cX_{exp(u_0)}(\overline w,z)$,
    which is indeed $K_{a_0}(w,z)$.

Let us return to the consideration of a general $F(s)\in
    L^2(\Re(s)=\frac12;\frac{|ds|}{2\pi})$. Under the
    hypothesis that $s\,F(s)$ is square integrable we have
    assigned  to $F$ the functions
\begin{subequations}
   \begin{align}
\alpha_F(u) &=
    2\int_{\Re(s)=\frac12}
    \cF(s)\cA_a(s)\frac{|ds|}{2\pi\,|\Gamma(s)|^2} = \int_{\Re(s)=\frac12}
    F(s)\overline{2\wh{A_a}(s)}\frac{|ds|}{2\pi} \\
\beta_F(u) &=
    2\int_{\Re(s)=\frac12}
    \cF(s)\cB_a(s)\frac{|ds|}{2\pi\,|\Gamma(s)|^2}  
= \int_{\Re(s)=\frac12}
    F(s)\overline{2\wh{\cB_a}(s)}\frac{|ds|}{2\pi} 
   \end{align}
  \end{subequations}
which are square-integrable at $+\infty$. From \eqref{eq:53}
    there holds,
    for any $a_0 = \exp(u_0)$:
    \begin{equation}
            \int_{\Re(s)=\frac12} X_{\exp(u_0)}(z,\overline
  s)F(s)\, \frac{|ds|}{2\pi} = 
  \int_{u_0}^\infty (2\wh{A_{\exp(u)}}(z)\alpha_F(u) + 2
  \wh{B_{\exp(u)}}(z)\beta_F(u))\,\frac{du}2 
    \end{equation}
The function of $z$ on the left side is the orthogonal
    projection $F_{a_0}$ of $F$ to the space
    $\wh{K_{a_0}}$. So, we deduce by unicity
    $\alpha_F(u)\Un_{u\geq u_0}(u) = \alpha_{F_{u_0}}(u)$
    and $\beta_F(u)\Un_{u\geq u_0}(u) =
    \beta_{F_{u_0}}(u)$. We then obtain 
    $\|F_{a_0}\|^2 = \int_{u_0}^\infty (|\alpha_F(u)|^2 +
    |\beta_F(u)|^2)\frac{du}2$  so $\alpha_F$ and $\beta_F$ are
    square-integrable on $(-\infty,+\infty)$, and as $\cup_a
    K_a$ is dense in $L^2(0,\infty;dx)$  the
    assignment $F\mapsto (\alpha_F,\beta_F)$ is isometric,
    and also it has a dense range in
    $L^2(\RR\to\CC^2;\frac{du}2)$. We can then remove the
    hypothesis that $s\,F(s)$ is square integrable and
    define the functions $\alpha_F$ and $\beta_F$ to be the
    limit in the $L^2$ sense of functions $\alpha_n$ and
    $\beta_n$ associated with $F_n$'s such that
    $\|F-F_n\|\to 0$ and the $s\,F_n$ are
    square-integrable. Summing up:
\begin{subequations}
\begin{theorem}\label{thm:expansion}
  There are  unitary identifications 
$L^2(0,\infty;dx)\wt\rightarrow L^2(\Re(s)=\frac12;\frac{|ds|}{2\pi})
\wt{\rightarrow} L^2(\RR\to \CC^2;\frac{du}2)$ 
given in the $L^2$ sense by the formulas, where $\Re(s)=\frac12$:
\begin{align}
F(s) &= \wh f(s) = \int_0^\infty f(x) x^{-s}\,dx\qquad\qquad\qquad
f(x) = \int_{\Re(s)=\frac12} F(s) x^{s-1}\,\frac{|ds|}{2\pi}\\
\alpha(u) &= \lim_{n\to\infty}
    \int_{\Re(s)=\frac12}
    F_n(s)\,
    \overline{2\wh{A_{\exp(u)}}(s)}\;\frac{|ds|}{2\pi}\qquad\quad
    (F_n\to_{L^2} F; \qquad s\,F_n(s)\in L^2) \\
\beta(u) &= \lim_{n\to\infty}
    \int_{\Re(s)=\frac12}
    F_n(s)\, \overline{2\wh{B_{\exp(u)}}(s)}\;\frac{|ds|}{2\pi} \\
\label{eq:500}    F(s) &=  \lim_{a_0\to0}
  \int_{\log(a_0)}^\infty (\alpha(u)\,2\wh{A_{\exp(u)}}(s) +
  \beta(u)\,2\wh{B_{\exp(u)}}(s))\,\frac{du}2 
\end{align}
The orthogonal projection of $f$ to $K_{a_0}$
corresponds to the replacement of $\alpha(u)$ by
$\alpha(u)\Un_{u>u_0}(u)$ and of $\beta(u)$ by
$\beta(u)\Un_{u>u_0}(u)$ ($u_0 = \log(a_0)$.).
The unitary
operators $f\mapsto \cH(f)$,  $F(s)\mapsto \chi(s)F(1-s)$, correspond to
$(\alpha,\beta)\mapsto (\alpha,-\beta)$. For $f = X_z^{a_0}$
   one has $\alpha(u) =
    2\wh{A_{\exp(u)}}(z)\Un_{u>\log(a_0)}(u)$ and $\beta(u) =
    -2\wh{B_{\exp(u)}}(z)\Un_{u>\log(a_0)}(u)$. The
  self-adjoint operator $F(s)\mapsto \gamma F(s)$
  ($s=\frac12+i\gamma$) corresponds to the canonical
  operator:
  \begin{equation}\label{eq:Hoperator}
    H = \begin{bmatrix}
      0&\frac{d}{du}\\-\frac{d}{du}&0\end{bmatrix} - 
    \begin{bmatrix} 0&\mu(e^u)\\\mu(e^u)&0\end{bmatrix}
  \end{equation} which,   in $L^2(\RR\to \CC^2;\frac{du}2)$,
is essentially self-adjoint  when defined on the domain of
the functions of class $C^1$ (or even $C^\infty$) with
compact support. The unitary operator $e^{i\,\tau H}$ acts
on $L^2(0,\infty;dx)$ as $f(x)\mapsto e^{\frac12\tau}
f(e^\tau x)$. 
\end{theorem}
\end{subequations} 

For the statement of self-adjointness we start with
$\alpha$ and $\beta$ of class $C^1$ with compact support,
define $F$ by \eqref{eq:500} and integrate by parts to
confirm that $\gamma F(s)$ corresponds to
$H\left(\left[\begin{smallmatrix}
\alpha\\\beta\end{smallmatrix}\right]\right)$. We know by
Hermann Weyl's theorem that in the limit point case the
pairs $\left[\begin{smallmatrix}
\alpha\\\beta\end{smallmatrix}\right]$ of class $C^1$ with
compact support are a core of self-adjointness (\textit{cf.}
\cite[\textsection13]{levsarg}.) On the other hand we know
that multiplication by $\gamma$ on
$L^2(\Re(s)=\frac12;\frac{|ds|}{2\pi})$ with maximal domain
is a self-adjoint operator. So the two self-adjoint
operators are the same.

Having discussed the matter from the point of view of the
isometric expansion we now turn to another topic, the topic
of the scattering, or rather total reflection against the
potential barrier at $+\infty$.  Another pair of solutions
of the first order system \eqref{eq:36} (hence also of the
second order differential equations) is known. Let us recall
from equations \eqref{eq:350a}, \eqref{eq:350b} that we
defined  $j_a = \sqrt{a}(\delta_a - \phi_a^+\Un_{0<x<a}) =
\sqrt{a} \cH\phi_a^+$ and $ -ik_a = \sqrt{a}(\delta_a +
\phi_a^-\Un_{0<x<a}) = \sqrt{a}\cH\phi_a^-$. Using again
\eqref{eq:90a} and \eqref{eq:90b} it is checked that $j_a$
and $k_a$ verify the exact same differential system as $A_a$
and $B_a$:
  \begin{subequations}
    \begin{align}\label{eq:92a}
      \ada j_a &= -\mu(a)j_a +i \xdtx k_a\\
\label{eq:92b}
      \ada k_a &= +\mu(a)k_a -i \xdtx j_a
    \end{align}
  \end{subequations}
The right Mellin transforms $\wh{j_a}(s)$ and
  $-i\wh{k_a}(s)$ are defined as
  \begin{subequations}
    \begin{align}\label{eq:64a}
      \wh{j_a}(s) &= a^{\frac12 - s} - \sqrt
  a\int_0^a \phi_a^+(x)x^{-s}\,dx\\ \label{eq:65a}
-i\wh{k_a}(s) &= a^{\frac12 -s} + \sqrt
  a\int_0^a \phi_a^-(x)x^{-s}\,dx
    \end{align}
  \end{subequations} As $\phi_a^+$ and $\phi_a^-$ are
analytic these are meromorphic functions in $\CC$ with
possible\footnote{the poles do exist.}  pole locations at $s
= 1$, $2$, \dots From the point of view of the Schrödinger
equation \eqref{eq:400a} and as $u = \log(a)\to-\infty$,
$a\to 0$, we thus see that, for $s = \frac12 + i\gamma$,
$\gamma\in\RR$,  $\wh{j_a}(\frac12+i\gamma)$ and
$\wh{j_a}(\frac12-i\gamma)$ are two (linearly independent
for $\gamma\neq0$)  solutions, differing from $e^{-i\gamma
u}$ and $e^{i\gamma u}$ by an exponentially small (in $u =
\log(a)$) quantity (and similarly with $-i\wh{k_a}$ with
respect to the Schrödinger equation \eqref{eq:95b}). So we
have identified the unique solutions which verify the Jost
conditions at $-\infty$.

As $P_a\phi_a^+$ is square integrable, and also using \eqref{eq:60a}, we
  have on the critical line:
  \begin{subequations}
    \begin{align}
    \label{eq:54}
 &\Gamma(s)\wh{j_a}(s) + \Gamma(1-s)\wh{j_a}(1-s)\\
 &=  \Gamma(s)\wh{j_a}(s)+\Gamma(1-s)a^{s-\frac12} -
  \Gamma(1-s)\sqrt{a}\int_0^a \phi_a^+(x)x^{s-1}\,dx \\
&= \Gamma(s)\wh{j_a}(s) + \Gamma(1-s)a^{s-\frac12} -
  \Gamma(s)\sqrt{a}\int_0^\infty (\cH P_a
    \phi_a^+)(x)x^{-s}\,dx\\
&=\Gamma(s)\wh{j_a}(s) + \Gamma(1-s)a^{s-\frac12} +
\Gamma(s)\sqrt{a}\int_0^\infty
  (\phi_a^+(x)-J_0(2\sqrt{ax}))x^{-s}\,dx\\
    \label{eq:6}
&= \Gamma(s)a^{\frac12-s}+\Gamma(1-s)a^{s-\frac12} -
  \Gamma(s)\sqrt{a}\int_0^a J_0(2\sqrt{ax})x^{-s}\,dx\notag \\
&\qquad\qquad
  +\Gamma(s)\sqrt{a}\int_a^\infty (\phi_a^+(x)-J_0(2\sqrt{ax}))x^{-s}\,dx      
    \end{align}
  \end{subequations}
As $J_0(2\sqrt{ax})-\phi_a^+(x)$ is square integrable both
  integrals are simultaneously absolutely convergent at
  least for $\frac12<\Re(s)<1$ (the $\frac12$ can be
  improved, but this does not matter). As the boundary
  values on the critical line coincide we have an identity
  of analytic functions. We recognize in $\int_a^\infty
  J_0(2\sqrt{ax})x^{-s}\,dx$, which is absolutely convergent
  for $\Re(s)>\frac34$, the quantity $g_s(a)$ (equation
  \eqref{eq:55}). And from equation \eqref{eq:50} we know
  $g_s(a) = \chi(s)a^{s-1} - \int_0^a
  J_0(2\sqrt{ax})x^{-s}\,dx$. So
  \begin{equation}
    \label{eq:27}
    \Gamma(s)\wh{j_a}(s) + \Gamma(1-s)\wh{j_a}(1-s)
 = \Gamma(s)a^{\frac12-s}+\Gamma(s)\sqrt{a}\int_a^\infty
  \phi_a^+(x)x^{-s}\,dx\;, 
  \end{equation} which is indeed $2\cA_a(s)$. From the
equation \eqref{eq:64a} $\wh{j_a}(s) = a^{\frac12}(a^{-s} -
\int_0^a \phi_a^+(x)x^{-s}\,dx)$ (valid as is for
$\Re(s)<1$) the function $u\mapsto \wh{j_a}(s)$ differs from
$u\mapsto e^{-i\gamma u}$ by an error which is relatively
exponentially smaller (we write $s=\frac12 + i\gamma$,
$\Im(\gamma)>-\frac12$). So $\wh{j_a}$ is the Jost solution
at $-\infty$ of the Schrödinger equation \eqref{eq:95a}. The
identity relating $\cB_a(s)$ and $\wh{k_a}(s) =
i\,a^{\frac12}(a^{-s} +  \int_0^a \phi_a^-(x)x^{-s}\,dx)$ is
proven similarly.

  \begin{theorem} The unique\footnote{here we make
use of the fact
      that \eqref{eq:95a} is in the limit point case at
      $+\infty$, because it is proven in the next chapter,
      or known from \eqref{eq:153a}, \eqref{eq:153b}, that
      $\mu(a) = 2a= 2\,e^u$, in this study of the $\cH$
      transform.} solution, square integrable at
$u=+\infty$, of the Schrödinger equation \eqref{eq:95a}
(resp.  \eqref{eq:95b}; $\gamma\neq0$) is expressed in terms
of the functions $\wh{j_a}(\frac12+i\gamma)$
(resp. $-i\,\wh{k_a}(\frac12+i\gamma)$) satisfying at
$-\infty$ the Jost condition $\wh{j_a}\sim_{u\to-\infty}
e^{-i\gamma u}$ (resp. $-i\wh{k_a}\sim_{u\to-\infty}
e^{-i\gamma u}$) as:
  \begin{subequations}
    \begin{align}\label{eq:66}
      \cA_a(s) &= \frac12\left(\Gamma(s)\wh{j_a}(s) +
  \Gamma(1-s)\wh{j_a}(1-s)\right) \\
     \cB_a(s) &= \frac12\left(\Gamma(s) \wh{k_a}(s) -
  \Gamma(1-s) \wh{k_a}(1-s)\right) 
    \end{align}
  \end{subequations}
  \end{theorem}

Let us add a time parameter $t$ and consider
  the wave equation:
  \begin{equation}
    \label{eq:65}
    \left( \frac{\partial^2}{\partial t^2} -
  \frac{\partial^2}{\partial u^2} + \mu^2 -
  \frac{d\mu}{du} \right) \Phi(t,u) = 0
  \end{equation}
Then $\Phi(t,u) = e^{i\gamma t}\wh{j_{\exp(u)}}(\frac12 +
  i\gamma)$ is a solution which behaves as $e^{i\gamma
  (t-u)}$ as $u\to-\infty$. This wave is thus right-moving,
  it is an incoming wave from $u=-\infty$ at
  $t=-\infty$. For a given frequency $\gamma$ there is a
  unique, up to multiplicative factor, wave which respects
  the condition of being at each time square integrable at
  $u = +\infty$. This wave is $e^{i\gamma t}
  \cA_{\exp(u)}(\frac12 + i\gamma)$. So the equation
  \eqref{eq:66} represents the decomposition in incoming and
  reflected components. There is in the reflected component
  a phase shift $\theta_\gamma = \arg \chi(s)$, the solution
  behaving approximatively at $u\to-\infty$ as
  $C(\gamma)\cos(\gamma u + \frac12\theta_\gamma)$. This is
  an absolute scattering, as there is nothing a priori to
  compare it too. We will thus declare that equation
  \eqref{eq:65} has realized $\chi(s)$ as an (absolute)
  scattering. Similarly the Schrödinger equation
  \eqref{eq:95b} realizes $-\chi(s)$ as an absolute
  scattering.

We have $2\cA_a(\frac12) = 2\Gamma(\frac12)\wh{j_a}(\frac12)$
and $\wh{j_a}(\frac12) = 1 - a^{\frac12}\int_0^a
\phi_a^+(x)x^{-\frac12}\,dx$. So $\lim_{a\to0}
\cA_a(\frac12) = \Gamma(\frac12) =\sqrt\pi$. On the other
hand $a\frac{d}{da} \cA_a(\frac12) = -\mu(a) \cA_a(\frac12)$
and $\mu(a) = a\frac{d}{da}
\log\frac{\det(1+H_a)}{\det(1-H_a)}$. so:
\begin{equation}
  \label{eq:71}
  \cA_a(\frac12) = \sqrt\pi\, \wh{A_a}(\frac12) = \sqrt\pi
  \;\frac{\det(1-H_a)}{\det(1+H_a)} 
\end{equation} 
We have $a\frac{d}{da} \|\cX_{\frac12}^a\|^2 = -
2\cA_a(\frac12)^2$. And $\cX_{\frac12}^a =
\Gamma(\frac12)X_{\frac12}^a$. So:
\begin{theorem}\label{thm:x12}
  The squared-norm of the evaluator $f\mapsto X_{\frac12}^a(f)
  = \int_a^\infty
  \frac{f(x)}{\sqrt x}\,dx$ on the Hilbert space $K_a$ of
  square integrable functions vanishing on $(0,a)$ and with
  $\cH$ transforms again vanishing on $(0,a)$ is:
  \begin{equation}
    \label{eq:73}
    \| X_{\frac12}^a \|^2 = 2 \int_a^\infty
    \left(\det\frac{1-H_b}{1+H_b}\right)^2\,\frac{db}b
  \end{equation}
  where $H_a$ is the restriction of $\cH$ to $L^2(0,a;dx)$.
\end{theorem}

It will be seen that $\det(1+H_a) = e^{a - \frac12 a^2}$ and
$\det(1-H_a) = e^{-a - \frac12 a^2}$. Having spent a long
time in the general set-up we now turn to determine
explicitely what the functions  $\phi_a^+$, $\phi_a^-$,
etc\dots are.

\section{The $K$-Bessel function in the theory of the $\cH$ transform} 

Let us recall that we may define the $\cH$ transform on all
of $L^2(\RR;dx)$ through the formula $\wt{\cH(f)}(\lambda) =
\frac{i}\lambda \wt{f}(\frac{-1}\lambda)$. This anticommutes
with $f(x)\mapsto f(-x)$, and $\cH$ leaves separately invariant
$L^2(0,+\infty;dx)$ and $L^2(-\infty,0;dx)$. We defined the
groups $\tau_a:f(x)\mapsto f(x-a)$ and $\tau_a^\# =
\cH\tau_a\cH$. We observed that the two groups are mutually
commuting, and that if the leftmost point of the support of
$f$ is at $\alpha(f)\geq0$ then the leftmost point of the
support of $\tau_b^\#(f)$, for any $b\geq0$, more precisely
for any $b\geq -\alpha(\cH(f))$, is still exactly
at $\alpha(f)$. From this we obtain the exact description of
$K_a$:
\begin{lemma}
  One has $K_a = \tau_a \tau_a^\# L^2(0,+\infty;dx)$.
\end{lemma}
Let now $Q$ be the orthogonal projection $L^2(\RR;dx)\mapsto
L^2(0,+\infty;dx)$. The orthogonal projection $Q_a$ from
$L^2(0,\infty;dx)$ to $K_a$ is thus exactly $\tau_a \tau_a^\#
Q \tau_{-a} \tau_{-a}^\#$. It will be easier to work with
$R_a = Q\tau_{-a}\tau_{-a}^\#$, especially as we are
interested in scalar products so we can skip the
$\tau_a\tau_a^\#$ isometry. First, we obtain $g_t^a(x) =
R_a(f_t)(x)$, for $f_t(x) = e^{-tx}$. The part of
$\tau_{-a}(f_t)$ supported in $x<0$ will be sent by
$\tau_{-a}^\#$ to a function supported again in $x<0$. We can
forget about it and we have thus first
$e^{-ta}e^{-tx}\Un_{x>0}(x)$, whose $\cH$ transform is
$e^{-ta}\frac1t \exp(-\frac xt)$, which we translate to the
left, again we cut the part in $x<0$, and we reapply $\cH$,
this gives $g_t^a(x) =
e^{-a(t+\frac1t)}e^{-tx}\Un_{x>0}(x)$. In other words we
have used in this computation:
\begin{equation}
  \label{eq:74}
  Q\tau_{-a}\tau_{-a}^\# = \cH Q\tau_{-a}\cH Q\tau_{-a}\qquad(a\geq0)
\end{equation}
The orthogonal projection $f_t^a := Q_a(f_t)$ of $f_t(x) =
e^{-tx}\Un_{x>0}(x)$ to $K_a$ is thus
$\tau_a\tau_a^\#(g_t^a)$. We can then compute exactly the
Fourier transform of $f_t^a$ as $f_t^a(i\tau) = (e^{-\tau
x}, \tau_a\tau_a^\#(g_t^a))_{L^2(\RR)}$ which is also
$(\tau_{-a}^\#\tau_{-a} e^{-\tau x}, g_t^a)_{L^2(\RR)} =
(g_\tau,g_t^a) =
e^{-a(t+\frac1t)}e^{-a(\tau+\frac1\tau)}\frac1{t+\tau}$. Hence:
\begin{lemma}
  The orthogonal projection $f_t^a$ to $K_a$ of
  $e^{-tx}\Un_{x>0}(x)$ has its Fourier transform
  $\wt{f_t^a}(\lambda)$ which is given as:
  \begin{equation}
    \label{eq:64}
    \wt{f_t^a}(i\tau) =
    e^{-a(t+\frac1t+\tau+\frac1\tau)}\frac1{t+\tau} 
  \end{equation}
\end{lemma}
The Gamma completed right Mellin transform $\cF_t^a(s)$ of
$f_t^a$ is the left Mellin transform of
$\wt{f_t^a}(i\tau)$. 
\begin{equation}
  \label{eq:77}
  \int_a^\infty f_t^a(x) \cX_s^a(x)\,dx = \cF_t^a(s) = 
  e^{-a(t+\frac1t)}\int_0^\infty
  e^{-a(\tau+\frac1\tau)}\frac{\tau^{s-1}}{t+\tau}\,d\tau 
\end{equation}
Let us write $W_s^a$ for the
element of $L^2(0,+\infty;dx)$ such that $\tau_a\tau_a^\#
W_s^a = \cX_s^a$. We have $\cF_t^a(s) = (\cX_s^a, f_t^a) =
(W_s^a, g_t^a)= e^{-a(t+\frac1t)}\int_0^\infty
W_s^a(x)e^{-tx}\,dx$. So the Laplace transform of $W_s^a(x)$
is exactly:
\begin{equation}
  \label{eq:78}
  \int_0^\infty W_s^a(x)e^{-tx}\,dx = \int_0^\infty
  e^{-a(\tau+\frac1\tau)}\frac{\tau^{s-1}}{t+\tau}\,d\tau 
\end{equation}
Writing $\frac1{t+\tau} = \int_0^\infty
e^{-(t+\tau)x}\,dx$, we recover $W_s^a(x)$ as:
\begin{equation}
  \label{eq:79}
  W_s^a(x) = \int_0^\infty
  e^{-a(\tau+\frac1\tau)}\tau^{s-1}e^{-\tau x}\,d\tau
\end{equation}
Then we obtain $\int_0^\infty W_s^a(x)W_z^a(x)\,dx$ which is
nothing else than $\cX_a(s,z)$:
\begin{theorem}
  The (analytic) reproducing kernel associated with the
  space of the
  completed right Mellin transforms of the elements of $K_a$
  is 
  \begin{equation}
  \label{eq:83}
  \cX_a(s,z) = \iint_{[0,+\infty)^2}
  e^{-a(t+\frac1t+u+\frac1u)}\frac{t^{s-1}u^{z-1}}{t+u}\,dtdu 
  \end{equation}
\end{theorem}
Here is a shortened argument: the analytic reproducing
kernel $\cX_a(s,z)$ is the completed right Mellin transform
of $\cX_s^a(x)$, so this is $\int_0^\infty (\cX_s^a, e^{-t
x})t^{s-1}\,dt$. But for $\Re(s)>\frac12$, $(\cX_s^a, e^{-t
x}) = \Gamma(s) (Q_a(x^{-s}\Un_{x>a}), e^{-t x}) = \Gamma(s)
(x^{-s}\Un_{x>a}, f_t^a) = \cF_t^a(s)$ ($Q_a$ is the
orthogonal projection to $K_a$). This gives again
\eqref{eq:83}.

To proceed further, we compute $(s+z-1)\cX_a(s,z)$. Using
integration by parts, multiplication by $s$ (resp. $z$) is converted
into $-t\frac{d}{dt}$ (resp.  $-u\frac{d}{du}$; there are no
  boundary terms.)
\begin{subequations}
\begin{align}
  s\cX_a(s,z) &= \iint_{[0,+\infty)^2}
  (a(t-\frac1t)+\frac{t}{t+u})e^{-a(t+\frac1t+u+\frac1u)}
  \frac{t^{s-1}u^{z-1}}{t+u}\,dtdu \\
  z\cX_a(s,z) &= \iint_{[0,+\infty)^2}
  (a(u-\frac1u)+\frac{u}{t+u})e^{-a(t+\frac1t+u+\frac1u)}
  \frac{t^{s-1}u^{z-1}}{t+u}\,dtdu 
\end{align}
\begin{align}
(s+z-1)\cX_a(s,z) &= a \iint_{[0,+\infty)^2}
  (t-\frac1t + u - \frac1u)e^{-a(t+\frac1t+u+\frac1u)}
  \frac{t^{s-1}u^{z-1}}{t+u}\,dtdu\notag\\
&= a\int_0^\infty e^{-a(t+\frac1t)}t^{s-1}\,dt\int_0^\infty
  e^{-a(u+\frac1u)}u^{z-1}\,du\notag \\
&\quad-a\int_0^\infty e^{-a(t+\frac1t)}t^{s-2}\,dt\int_0^\infty
  e^{-a(u+\frac1u)}u^{z-2}\,du
\end{align}
\end{subequations}
The $K$-Bessel function is $K_s(x) = \frac12\int_0^\infty
  e^{-x\frac12(t+\frac1t)}t^{s-1}\,dt = \int_0^\infty
  e^{-x\cosh u}\cosh(su)\,du$. It is an even function of
  $s$. It has, for each $x>0$, all its zeros on the imaginary
  axis, and was used by P\'olya in a famous work on 
  functions inspired by the Riemann $\xi$-function and for
  which he proved the validity of the Riemann hypothesis
  \cite{polya1,polya2}. We have obtained the formula
  \begin{equation}
    \label{eq:85}
    \cX_a(s,z) = \frac{E(s)E(z)- E(1-s)E(1-z)}{s+z-1}\qquad
  E(s) = 2\sqrt a K_s(2a)
  \end{equation}
To confirm $\cE_a(s) = 2\sqrt a K_s(2a)$, let us define
  temporarily $A(s) = \frac12 \sqrt a \int_0^\infty
  e^{-a(t+\frac1t)} (1+\frac1t)t^{s-1}\,dt$ and $-iB(s)
  =\frac12 \sqrt a \int_0^\infty e^{-a(t+\frac1t)}
  (1-\frac1t)t^{s-1}\,dt$ which are respectively even and
  odd under $s\mapsto 1-s$ and are such that $E(z) = A(z) -
  i B(z)$. We have $
    \forall s,z\in\CC\quad {-iB(s)A(z) + A(s)(-iB(z))}
    ={-i\cB_a(s)\cA_a(z) + \cA_a(s)(-i\cB_a(z))}$ and
    considering separately the
  even and odd parts in $z$, we find that there exists a
  constant $k(a)$ such that $A(s) = k(a)\cA_a(s)$ and $B(s)
  = k(a)^{-1} \cB_a(s)$. Let us check that
  $\lim_{\sigma\to\infty} \frac{-iB(\sigma)}{A(\sigma)} =
  1$. It is a corollary to $\lim_{\sigma\to\infty}
  K_{\sigma}(x)/K_{\sigma+1}(x) = 0$ which is elementary:
  $\int_{-\infty}^0 \exp(-x\cosh u)e^{\sigma u}\,du = O(1)$
  ($\sigma\to+\infty$), and for each $T>0$, $\int_T^\infty
  \exp(-x\cosh u)e^{\sigma u}\,du \geq T\exp(-x\cosh
  3T)e^{2\sigma T}$, $\int_0^T \exp(-x\cosh u)e^{\sigma
  u}\,du\leq Te^{\sigma T}$, and combining we get
  $K_\sigma(x) = (1+ o(1)) \frac12\int_T^\infty \exp(-
  x\cosh u)e^{\sigma u}\,du$. So $\limsup_{\sigma\to+\infty}
  \frac{K_\sigma(x)}{K_{\sigma+1}(x)} \leq e^{-T}$ for each
  $T>0$.  Using \eqref{eq:70}, we then conclude $k(a) = 1$.

Let us examine the equality $\cE_a(s) = 2\sqrt{a}K_s(2a)= \sqrt a\int_0^\infty
  \exp({-a(t+\frac1t)})t^{s-1}\,dt$. It exhibits $\cE_a$ as the
  left Mellin transform of $\sqrt
  a\exp(-a(t+\frac1t))$, so the distribution $E_a$ is
  determined as the distribution whose Fourier transform is
  $\sqrt a\exp(i\,a(\lambda -\lambda^{-1}))$. Using  $\tau_a$ and
  $\tau_a^\#$, this means exactly:
  \begin{equation}
    \label{eq:87}
    E_a = \sqrt{a}\,\tau_a^\#\,\tau_a\delta = \sqrt{a}\,\cH\tau_a\cH\delta(x-a)
  \end{equation}
We exploit the symmetry $K_s = K_{-s}$, which corresponds to
$\wt{E_a}(\lambda) = \wt{E_a}(-\lambda^{-1}) =
  -i\lambda\wt{\cH E_a}(\lambda)$, so the unexpected identity appears:
  \begin{equation}
    \label{eq:88}
    E_a = \frac{d}{dx} \cH E_a
  \end{equation}
From \eqref{eq:87} we read $\cH E_a = \sqrt{a}\tau_a
  J_0(2\sqrt{ax}) =
  \sqrt{a}J_0(2\sqrt{a(x-a)})\Un_{x>a}(x)$. Using
  \eqref{eq:88}, and recalling
  equations \eqref{eq:89c} and \eqref{eq:89d} we deduce:
\begin{subequations}
  \begin{align}
    \label{eq:90}
    \frac{\phi_a^+(x)+\phi_a^-(x)}2 &= J_0(2\sqrt{a(x-a)}) =
    I_0(2\sqrt{a(a-x)}\\
    \label{eq:91}
    \frac{\phi_a^+(x)-\phi_a^-(x)}2 &=
   \frac{\partial}{\partial x} J_0(2\sqrt{a(x-a)}) = 
  \frac{\partial}{\partial x} I_0(2\sqrt{a(a-x)})
  \end{align}
\end{subequations}
We knew already from equations \eqref{eq:153a},
\eqref{eq:153b}! Summing up we have proven:
 
\begin{subequations}
\begin{theorem}
Let $\cH$ be the self-reciprocal operator with kernel
$J_0(2\sqrt{xy})$ on $L^2(0,\infty;dx)$. Let $H_a$ be the
restriction of $\cH$ to $L^2(0,a;dx)$. The solutions to
the integral equations $\phi_a^+ + H_a\phi_a^+ =
J_0(2\sqrt{ax})$ and $\phi_a^- -H_a\phi_a^- =
J_0(2\sqrt{ax})$ are the entire functions:
\begin{align}\label{eq:100a}
\phi_a^+(x) &= I_0(2\sqrt{a(a-x)}) +
\frac{\partial}{\partial x} I_0(2\sqrt{a(a-x)})\\ 
\label{eq:100b}
\phi_a^-(x) &= I_0(2\sqrt{a(a-x)}) -
\frac{\partial}{\partial x} I_0(2\sqrt{a(a-x)})
\end{align}
One has $1-a = \phi_a^+(a) = \frac{d}{da}\log\det(1+H_a)$ and
$1+ a = \phi_a^-(a) =  -\frac{d}{da}\log\det(1-H_a)$, and 
\begin{align}
  \det(1+H_a) &= e^{+a - \frac12 a^2}\\
  \det(1-H_a) &= e^{-a -\frac12 a^2} 
\end{align}
The  tempered distributions $A_a =
\frac{\sqrt{a}}2(1+\cH)\phi_a^+$ and  $-i\,B_a =
\frac{\sqrt{a}}2(-1+\cH)\phi_a^-$, 
respectively invariant and anti-invariant under 
$\cH$, are also given as:
\begin{align}\label{eq:200a}
  A_a(x) &= \frac{\sqrt a}2\;\Big( \delta_a(x) +
\Un_{x>a}(x)\big(J_0(2\sqrt{a(x-a)}) + 
\frac{\partial}{\partial x} J_0(2\sqrt{a(x-a)})\big)\Big)\\
\label{eq:200b}
  -i B_a(x) &= \frac{\sqrt a}2\;\Big( \delta_a(x) -
\Un_{x>a}(x)\big(J_0(2\sqrt{a(x-a)}) - 
\frac{\partial}{\partial x} J_0(2\sqrt{a(x-a)})\big)\Big)
\end{align}
Their Fourier transforms  are $\int_\RR e^{i\lambda
x}A_a(x)\,dx = \frac{\sqrt{a}}2 (1 + \frac
i\lambda)\exp(ia(\lambda-\frac1\lambda))$ and $-i\,\int_\RR
e^{i\lambda x}B_a(x)\,dx = \frac{\sqrt{a}}2 (1 - \frac
i\lambda)\exp(ia(\lambda-\frac1\lambda))$. The 
Gamma completed right Mellin transforms are:
\begin{align}
  \Gamma(s)\wh{A_a}(s) &= \cA_a(s) = \sqrt{a}(K_s(2a) + K_{s-1}(2a))\\
-i\Gamma(s)\wh{B_a}(s) &= -i\cB_a(s) = \sqrt{a}(K_s(2a) - K_{s-1}(2a))\\
  \cA_a(s) - i \cB_a(s) &= \cE_a(s)  = 2\sqrt{a}\;K_s(2a) =
  \sqrt{a}\int_0^\infty e^{-a(t+\frac1t)}t^{s-1}\,dt
\end{align}
They verify the
first order system, where $\mu(a) = a\phi^+(a)+a\phi^-(a) = 2a$:
\begin{equation}
      \left(\begin{bmatrix}
      0&1\\-1&0\end{bmatrix}\,a\frac{d}{da} - 
    \begin{bmatrix} 0&\mu(a)\\\mu(a)&0\end{bmatrix}\right)
      \begin{bmatrix} \cA_a(s)\\  \cB_a(s) \end{bmatrix} 
= -i(s-\frac12) \begin{bmatrix} \cA_a(s)\\ \cB_a(s) \end{bmatrix}
\end{equation} The pair $\left[\begin{smallmatrix}
\cA_a(s)\\  \cB_a(s) \end{smallmatrix}\right]$ is the unique
solution  of the first order system which is
square-integrable with respect to $d\log(a)$ at
$+\infty$. The total reflection against the exponential
barriers at $\log(a)\to+\infty$ of the associated Schr\"odinger
equations realizes $+\frac{\Gamma(1-s)}{\Gamma(s)}$ and
$-\frac{\Gamma(1-s)}{\Gamma(s)}$ ($\Re(s) = \frac12$) as
scattering matrices.
\end{theorem}

\end{subequations}

From \eqref{eq:71} we have $\cA_a(\frac12) = \sqrt \pi\,
e^{-2a}$. To normalize $\cA_a$ according to $\cA_a(\frac12)
=1$, we would have to make the replacement $\cA_a\to
\pi^{-\frac12}\, e^{2a}\cA_a$ and $\cB_a\to \sqrt \pi\,
e^{-2a}\cB_a$ and the expression of $\cE_a$ in terms of the
$K$-Bessel function would be less simple.  Let us
also note that according to \eqref{eq:70} we must have
$\frac{K_{\sigma-1}(2a)}{K_\sigma(2a)}
\sim_{\sigma\to+\infty} \frac a\sigma$.

Regarding the isometric expansion, as given in theorem
\ref{thm:expansion}, we apply it to a function $F(s) =
\int_0^\infty k(x)x^{-s}\,dx$ such that $\frac{d}{dx} x\,
k(x)$ as a distribution on $\RR$ is in $L^2$.  Using the
$L^2$-function $\frac1s \wh{A_a}(s)$, which is the Mellin
transform of the function $C_a(x) = \frac1x \int_0^x
A_a(x)\,dx$, and the Parseval identity we obtain $\alpha(u)
= 2\int_0^\infty (-x\frac{d}{dx} k(x))C_a(x)\,dx$ as an
absolutely convergent integral. The square
integrable function $C_a(x)$ is explicitely:
\begin{equation}
  \label{eq:93}
  C_a(x) = \frac{\sqrt
  a}{2x}\Un_{x>a}(x)\Big(J_0(2\sqrt{a(x-a)}) 
  +\sqrt{\frac{x-a}a}J_1(2\sqrt{a(x-a)}) \Big)
\end{equation}
And under the hypothesis made on $x\frac{d}{dx} k(x)$ we
obtain the existence of:
\begin{equation}
  \label{eq:94}
  \alpha(u) = \lim_{X\to\infty} \sqrt a k(a) - 2
  Xk(X)C_a(X) + \sqrt a\int_a^X k(x) (1+\frac{\partial}{\partial
  x}) J_0(2\sqrt{a(x-a)})\,dx
\end{equation}
Let us observe that $k(x) = \int_x^\infty \frac{l(y)}y\,dy$
with $l(y)\in L^2$, so  $|k(x)|^2 \leq
\frac Cx$ and $Xk(X)C_a(X) = O(X^{-\frac14})$. Hence:
\begin{equation}
  \label{eq:95}
    \alpha(u) = \sqrt a k(a)  + \sqrt a\int_a^{\to\infty} k(x)
  (1+\frac{\partial}{\partial 
  x}) J_0(2\sqrt{a(x-a)})\,dx
\end{equation}
Comparing with equation \eqref{eq:9a}
we see that the $f(y)$ defined there is related to
$\alpha(u)$ , $u=\log(a)$ by the formula $f(y) =
\frac1{2\sqrt a}
\alpha(\log(a))$, $a = \frac y2$, so $|f(y)|^2\,dy =
\frac1{4a}|\alpha(\log(a))|^2 2\,da = |\alpha(\log(a))|^2
\frac12 d\log(a)$. Similarly we obtain $\beta(u)$:
\begin{equation}
  \label{eq:96}
    \beta(u) = \sqrt a k(a)  + \sqrt a\int_a^{\to\infty} k(x)
  (-1+\frac{\partial}{\partial 
  x}) J_0(2\sqrt{a(x-a)})\,dx
\end{equation}
and comparing with \eqref{eq:9b}, $g(y) = \frac1{2\sqrt a}
\beta(\log(a))$, $|g(y)|^2\,dy  = |\beta(\log(a))|^2
\frac12 d\log(a)$. So according to theorem
\ref{thm:expansion} we do have equation \eqref{eq:9d}:
\begin{equation}
  \label{eq:97}
  \int_0^\infty (|f(y)|^2 + |g(y)|^2)\,dy = \int_0^\infty |k(x)|^2\,dx
\end{equation}
From \ref{thm:expansion} the assignment $k\to(\alpha,\beta)$
extends to a unitary identification
$L^2(\Re(s)=\frac12;\frac{|ds|}{2\pi}) \wt{\rightarrow}
L^2(\RR\to \CC^2;\frac{du}2)$, which has the property
$\cH(k)\to(\alpha,-\beta)$. In order to complete the proof
of the isometric expansion, it 
remains to check the equation $\eqref{eq:9c}$ which expresses $k$
in terms of $f$ and $g$. According to \ref{thm:expansion} we
recover $k(x)$ has the inverse Mellin transform of $
  \int_\RR (\alpha(u)\,2\wh{A_{a}}(s) +
  \beta(u)\,2(-i\wh{B_{a}}(s)))\,\frac{du}2 $. Expressing
this in terms of $f(y)$ and $g(y)$, $y = 2a$, $u=\log(a)$,
this means the identity of distributions, where we suppose for
simplicity that $f(y)$ and $g(y)$ have compact support in
$(0,+\infty)$ (as usual, this means having support away from
$0$ as well as $\infty$.):
\begin{equation}
  \label{eq:84}
\begin{split}
  k(x)  &= \int_0^\infty \left( 2\sqrt{\frac
  y2}f(y)\,2 A_{\frac y2}(x) 
  + 2\sqrt{\frac y2}g(y)\,2(-i B_{\frac y2}(x)) \right)
  \frac{dy}{2y}\\
&= 2\int_0^\infty \left( \sqrt{y}f(2y)\,A_y(x) 
  + \sqrt{y}g(2y)\,2(-i B_y(x)) \right)
\frac{dy}{y}
\end{split}
\end{equation}

Then imagining that we
are integrating against a test function $\psi(x)$ and using
Fubini we obtain:
\begin{equation}
  \label{eq:86}\begin{split}
  &2\int_0^\infty  \sqrt{y}f(2y)\,A_y(x) \frac{dy}{y} \\
  &= \int_0^\infty f(2y)\left(\delta(x-y) +
  \Un_{x>y} \Big(J_0(2\sqrt{y(x-y)}) -
  \sqrt{\frac{y}{x-y}}J_1(2\sqrt{y(x-y)})\Big)\right) \,dy\end{split}
\end{equation}
\begin{equation}\begin{split}
 &=  f(2x) +  \int_0^x \Big(J_0(2\sqrt{y(x-y)}) -
  \sqrt{\frac{y}{x-y}}J_1(2\sqrt{y(x-y)})\Big) f(y)\,dy\\
 &=  f(2x) + \frac12\int_0^{2x} \Big(J_0(\sqrt{y(2x-y)}) -
  \sqrt{\frac{y}{2x-y}}J_1(\sqrt{y(2x-y)})\Big) f(y)\,dy
\end{split}
\end{equation}
Proceeding similarly with $g(y)$ one obtains for
$2\int_0^\infty  \sqrt{y}g(2y)(-i \,B_y(x)) \frac{dy}{y}$:
\begin{equation}
  \label{eq:89}
  g(2x)
  - \frac12\int_0^{2x} \Big(J_0(\sqrt{y(2x-y)}) + 
  \sqrt{\frac{y}{2x-y}}J_1(\sqrt{y(2x-y)})\Big) f(y)\,dy
\end{equation}
Combining \eqref{eq:86} and \eqref{eq:89} in the formula
\eqref{eq:84} for
$k(x)$ we obtain equation \eqref{eq:9c}.

\section{The reproducing kernel and differential equations
    for the extended spaces}
\label{sec:8}

\kern-12pt

Let $L_a\subset L^2(0,\infty;dx)$ be the Hilbert space of
square integrable functions $f$ which are constant in
$(0,a)$ and with their $\cH$-transforms again constant in
$(0,a)$. The distribution $x\frac{d}{dx}\,\frac{d}{dx}x\,f
=\frac{d}{dx}x\,x\frac{d}{dx} f$ vanishes in $(0,a)$ and its
$\cH$ transform does too. So $s(s-1)\wh f(s)$ is an entire
function with trivial zeros at $-\NN$. The Hilbert space of
the functions $s(s-1)\Gamma(s)\wh f(s)$ satisfies the axioms
of \cite{Bra}; we prove everything according to the methods
developed in the earlier chapters.  Our goal is to determine
the evaluators and reproducing kernel for the spaces $L_a$.

\kern-6pt

For $f\in L_a$,  $\wh f(s)$ is a meromorphic function with at
most a pole at $s=1$ and also $\wh f(0)$ does not
necessarily vanish. The Mellin-Plancherel transform
$\int_0^\infty f(x)x^{-s}\,dx = \int_0^a c(f)x^{-s}\,dx +
\int_a^\infty f(x) x^{-s}\,dx$ has polar part
$-\frac{c(f)}{s-1}$. Let us write $(f,Y_1) = -c(f) =
\Res(\wh f(s),s=1) = s(s-1)\Gamma(s)\wh f(s)|_{s=1}$. This
defines an element $Y_1\in L_a$. We define also $\cY_1^a =
\Gamma(1)Y_1 = Y_1$. Then $(f,\cY_1^a) = s(s-1)\Gamma(s)\wh
f(s)|_{s=1}$. We also define $\cY_0^a$ as the vector such that
$(f,\cY_0^a) = s(s-1)\Gamma(s)\wh f(s)|_{s=0} = -\wh
f(0)$. One observes $(f,\cH(Y_1)) = (\cH(f), Y_1) =
s(s-1)\Gamma(s)\wh{\cH(f)}(s)|_{s=1} = s(s-1)\Gamma(1-s)\wh
f(1-s)|_{s=1} = -\wh f(0) = (f, \cY_0^a)$ so $\cY_0^a =
\cH(\cY_1^a)$. To lighten the notation we sometimes write $\cY_1$
  and $\cY_0$ instead $\cY_1^a$ and $\cY_0^a$ when no
confusion can arise.

\kern-6pt

We will also consider the vectors $X_s^\times\in L_a$ such that
$\forall f\in L_a\; \wh f(s) =
  (X_s^\times,f)$.\footnote{sometimes written
  $X_s^{a\times}$.} The orthogonal
projection of $X_s^\times$ to $K_a\subset L_a$ is $X_s$. Let us
look more closely at this orthogonal projection. First let
$N_a$ be the (closed) vector space sum  $L^2(0,a;dx) +
\cH L^2(0,a;dx)$. Inside $N_a$ we have the codimension two
space $M_a$ defined as the sum of $(\Un_{0<x<a})^\perp\cap
L^2(0,a;dx)$ and of its image under $\cH$. Finally, let
$R_a$ be the orthogonal complement in $N_a$ of $M_a$, which
has dimension two. For a
function $f$ to belong to $L_a$ it is necessary and
sufficient that its orthogonal projection to $N_a$ be
perpendicular to the functions in $L^2(0,a;dx)$ which are
perpendicular to $\Un_{0<x<a}$, and the same for the
$\cH$-transform, so this means exactly that its orthogonal
projection to $N_a$ belongs to $R_a$. So we have the
orthogonal decomposition of 
$L^2(0,\infty;dx)$ into the sum of the three spaces $K_a$,
$R_a$ and  $M_a$ and $L_a = K_a \oplus R_a$. For $f\in L_a$
to be in $K_a$ it is necessary and sufficient that $c(f) =
-(f,\cY_1^a) = 0$ and the same for $c(\cH(f))$, so this means
that $\{\cY_1^a, \cY_0^a\}$ is a basis of $R_a$. The function
$\cY_1^a$ belongs to $N_a = L^2(0,a;dx) +
\cH L^2(0,a;dx)$ and as such is uniquely written as $u_1 +
\cH v_1$. As $\cY_1^a\in L_a$ we have constants $\alpha,
\beta\in\CC$ such that:
\begin{subequations}
  \begin{align}
    u_1 + H_a v_1 &= -\alpha\\
    H_a u_1 + v_1 &= -\beta
  \end{align}
where we recall that $P_a$ is the restriction to
$(0,a)$ and $H_a = P_a \cH P_a$, $D_a = H_a^2$. From what
was said previously $\alpha = 
(\cY_1^a,\cY_1^a)$ and $\beta = (\cH(\cY_1^a),\cY_1^a) =
(\cY_0^a,\cY_1^a)$. We have thus:
  \begin{align}
    u_1  &= (1 - D_a)^{-1}(-\alpha\Un_{0<x<a} + \beta
    H_a(\Un_{0<x<a}))\\
\label{eq:120}
    v_1  &= (1 - D_a)^{-1}(+\alpha H_a(\Un_{0<x<a}) - \beta\Un_{0<x<a})
  \end{align}
\end{subequations}

Also the function $\cY_1^a$ may be obtained as the orthogonal
projection to $L_a$ of $-\frac1a \Un_{0<x<a}$. Indeed it
follows from what has been seen above that for
any element $f\in L_a$, $(f,\cY_1^a) = -\frac1a \int_0^a
f(x)\,dx$. As the function $-\frac1a \Un_{0<x<a}$ already
belongs to $N_a$, we have
\begin{equation}
  -\frac1a \Un_{0<x<a} = u_1 + \cH v_1 + u_2 + \cH v_2
\end{equation}
where $u_2+\cH v_2$ belongs to $M_a$, which means that
$u_2\in L^2(0,a;dx)$ verifies $\int_0^a u_2(x)\,dx = 
0$ and $v_2 \in L^2(0,a;dx)$ verifies $\int_0^a v_2(x)\,dx =
0$. But there is unicity so we have exactly
\begin{equation}
  u_1 + u _2 = -\frac1a \Un_{0<x<a}\qquad\qquad v_1+ v_2 = 0
\end{equation}
And we deduce:
\begin{equation}
  \label{eq:98}
  \int_0^a u_1(x)\,dx = -1\qquad\qquad \int_0^a v_1(x)\,dx = 0
\end{equation}
So $\alpha$ and $\beta$ are determined as the solutions
of the system:
\begin{subequations}
  \begin{align}
    \alpha (\Un_{0<x<a},(1-D_a)^{-1}\Un_{0<x<a}) -
    \beta(\Un_{0<x<a},(1-D_a)^{-1}H_a\Un_{0<x<a}) &= 1\\
   \alpha(\Un_{0<x<a},(1-D_a)^{-1}H_a\Un_{0<x<a}) -
    \beta(\Un_{0<x<a},(1-D_a)^{-1}\Un_{0<x<a}) &= 0
  \end{align}
\end{subequations}

We thus have:
\begin{subequations}
\begin{proposition}\label{prop:pq}
  Let $p(a)$ and $q(a)$ be defined as 
  \begin{align}
    p(a) &= \int_0^a (1 - D_a)^{-1}(\Un_{0<x<a})(x)\,dx\\
    q(a) &= \int_0^a (1 - D_a)^{-1}H_a(\Un_{0<x<a})(x)\,dx
  \end{align}
  then:
  \begin{equation}
    \label{eq:92}
    \begin{bmatrix} p(a)&-q(a)\\-q(a)&p(a) \end{bmatrix}
    \begin{bmatrix} (\cY_1,\cY_1)&(\cY_1,
\cY_0)\\(\cY_0,\cY_1) & (\cY_0,\cY_0) \end{bmatrix}
=
    \begin{bmatrix} 1&0\\0&1 \end{bmatrix}
  \end{equation}
 The evaluators $\cY_1^a$ and $\cY_0^a = \cH(\cY_1^a)$ are given as $u_1 +
 \cH v_1$ and $\cH u_1 + v_1$ with:
  \begin{align}
    u_1 &= -(\cY_1,\cY_1)(1 - D_a)^{-1}(\Un_{0<x<a}) +
    (\cY_0,\cY_1) (1 - D_a)^{-1} H_a(\Un_{0<x<a})\\ v_1 &=
    - (\cY_0,\cY_1)(1 - D_a)^{-1}(\Un_{0<x<a})
    +(\cY_0,\cY_0) (1 - D_a)^{-1} H_a(\Un_{0<x<a})
  \end{align}
\end{proposition}
\end{subequations}

We have introduced, for $s\neq 0,1$, $X_s^\times$ as the
evaluator $\wh f(s)$ for functions in $L_a$. We shall write
$\cX_s^\times = \Gamma(s)X_s^\times$ and then $\cY_s = s(s-1)
\cX_s^\times$. This is compatible with our previous definitions
of $\cY_1^a$ and $\cY_0^a$. We note that the orthogonal
projection of $\cX_s^\times$ to $K_a$ is $\cX_s$. So we may write
\begin{equation}
  \label{eq:99}
  \cX_s^\times = \cX_s + \lambda(s) \cY_1^a + \mu(s) \cY_0^a
\end{equation}
We shall also write $\cY_a(s,z) = \int_0^\infty
\cY_s^a(x)\cY_z^a(x)\,dx = z(z-1)\Gamma(z)\wh{\cY_s^a}(z)$. One has
$\cH(\cY_s^a) = \cY_{1-s}^a$ and $\cY_a(s,z) = \cY_a(1-s,1-z) =
\cY_a(z,s)$. Taking
the scalar products with $\cY_1^a$ and $\cY_0^a$ in
\eqref{eq:99} we obtain
\begin{subequations}
  \begin{align}
    \frac1{s(s-1)} \cY_a(1,s)  &= \lambda(s) \alpha +
    \mu(s) \beta\\
    \frac1{s(s-1)} \cY_a(1,1-s)  &= \lambda(s) \beta +
    \mu(s) \alpha
\end{align}
\begin{align}
   \lambda(s) &=  \frac1{s(s-1)} (p \cY_a(1,s) - q
    \cY_a(1,1-s))\\
   \mu(s) &= \frac1{s(s-1)} (-q\cY_a(1,s) + p
    \cY_a(1,1-s)) = \lambda(1-s)
  \end{align}
\end{subequations}
Combining with \eqref{eq:99}, this gives:
\begin{equation}
  \label{eq:102}
     \cY_s = s(s-1)\cX_s + \cY_a(1,s)(p \cY_1^a - q \cY_0^a) +
    \cY_a(1,1-s) (-q \cY_1^a + p \cY_0^a)
\end{equation}
\begin{equation}
  \label{eq:124}
  \text{Let}\quad T_a(s) = p(a) \cY_a(1,s) - q(a) \cY_a(1,1-s)
\end{equation}

\begin{subequations}
\begin{proposition}\label{prop:y}
  The (analytic) reproducing kernel $\cY_a(s,z)$ of the extended space
  $L_a$ is given by each of the following expressions:
  \begin{gather}
    \label{eq:103}
    s(s-1)z(z-1)\cX_a(s,z) + 
    \begin{bmatrix} \cY_a(1,s) & \cY_a(1,1-s) \end{bmatrix}
    \begin{bmatrix} p(a)&-q(a)\\-q(a)&p(a) \end{bmatrix}
    \begin{bmatrix} \cY_a(1,z) \\ \cY_a(1,1-z) \end{bmatrix}\\
  \label{eq:104}
    = s(s-1)z(z-1)\cX_a(s,z) + 
    \begin{bmatrix} T_a(s) & T_a(1-s) \end{bmatrix}
    \begin{bmatrix} \alpha(a)&\beta(a)\\\beta(a)&\alpha(a) \end{bmatrix}
    \begin{bmatrix} T_a(z) \\ T_a(1-z) \end{bmatrix}\\
    = s(s-1)z(z-1)\cX_a(s,z) + T_a(s)\cY_a(1,z) +
    T_a(1-s)\cY_a(1-z)
\end{gather}
\end{proposition}
\end{subequations}

A very important observation, before turning to the
determination of the quantities $p(a)$ and $q(a)$ shall now
be made. Let $L$ be the unitary operator: 
\begin{equation}
  \label{eq:105}
  L(f)(x) = f(x)  - \frac1x \int_0^x f(y)\,dy
\end{equation}
It is the operator of multiplication by $\frac{s-1}s$ at the
level of right Mellin transforms. Obviously it converts
functions constant on $(0,a)$ into functions vanishing on
$(0,a)$. Let us now consider the
operator
\begin{equation}
  \label{eq:106}
  \cH^\etendu = L\;\cH\; L^{-1} = L^2\,\cH = \cH L^{-2}
\end{equation}
It is a unitary, self-adjoint, self-reciprocal, scale
reversing operator whose kernel is easily computed to be
\begin{equation}
  \label{eq:107}
  k^\etendu(xy) = J_0(2\sqrt{xy}) - 2 \frac{J_1(2\sqrt{xy})}{\sqrt{xy}} +
  \frac{1 - J_0(2\sqrt{xy})}{xy} = \sum_{n=0}^\infty (-1)^n
  \frac{n^2 x^n y^n}{(n+1)!^2}
\end{equation}
It has $L(e^{-x}) = (1+\frac1x) e^{-x} -\frac1x$ as
self-reciprocal function; the Mellin transform is
$\frac{s-1}s\Gamma(1-s)$ which, multiplied by $s(s-1)$ gives
$(1-s)^2\Gamma(1-s)$ which is the Mellin transform of a more
convenient invariant function for $\cH^\etendu$, the function
$x(x-1)e^{-x}$. This function is the analog for $\cH^\etendu$
of $e^{-x}$ for $\cH$. Let us now consider the space
$L(L_a)$. It consists of the square integrable functions
vanishing identically on $(0,a)$ and having $\cH^\etendu$
transforms also identically zero on $(0,a)$. But then the
entire theory applies to $\cH^\etendu$ exactly as it did for
$\cH$, up to some minor details in the proofs where the
function $J_0$ was really used like in Lemma \ref{lem:o1} or
proposition \ref{prop:10}. We have for $\cH^\etendu$ versions
of all quantities previously considered for $\cH$. To check
that the proof of \ref{prop:10} may be adapted, we need to
look at $k_1^\etendu(x) = \int_0^x k^\etendu(t)\,dt =
\sqrt{x}J_1(2\sqrt{x}) + 2(J_0(2\sqrt
x) - 1) + 2 \int_0^{2\sqrt x} \frac{1 -
  J_0(u)}u\,du = O_{x\to+\infty}(x^{\frac14})$. So we may
  employ lemma
\ref{lem:Ok} as was done for $\cH$. Proposition
\ref{prop:10} and Theorem \ref{thm:11} thus hold. We must be
careful that the operator $L^{-1}$ is always involved when
comparing functions or distributions related to $\cH^\etendu$
with those related to $\cH$. For example, one has
$X_s^\times = \frac{s}{s-1} L^{-1} X_s^\etendu$ and $\cY_s =
s(s-1)\cX_s^\times = L^{-1} \cX_s^\etendu$. The two types of
Gamma completed Mellin transforms differ: for $\cH$ we
consider $\Gamma(s)\wh f(s)$ while for $\cH^\etendu$ we
consider $s^2\Gamma(s)\wh g(s)$. Indeed this is quite the
coherent thing to do in order that:
\begin{equation}
  \label{eq:108}
  s^2\Gamma(s)\wh g(s)  = s(s-1)\Gamma(s) \wh f(s)
\end{equation}
for $g = L(f)$. The bare Mellin transforms of elements of
spaces $K_a^\etendu$ are not always entire in the complex plane:
they may have a pole at $s = 0$. After multiplying by
$s^2\Gamma(s)$ which is the \emph{left} Mellin transform of
the self-invariant function $x(x-1)e^{-x}$, as $\Gamma(s)$
is the \emph{left} Mellin transform of $e^{-x}$, we do
obtain entire functions, whose trivial zeros are at $-1$,
$-2$, \dots ($0$ is not a trivial zero anymore.) From
equation \eqref{eq:108} we see that the (analytic)
reproducing kernel $\cX_a^\etendu(s,z)$ exactly coincides with the
function $\cY_a(s,z)$ whose initial computation has been
given in Proposition \ref{prop:y}. Also the 
Schrödinger equations  will realize $\pm\left(\frac{1-s}s\right)^2
\frac{\Gamma(1-s)}{\Gamma(s)}$ as scattering matrices, and
there will be an isometric expansion generalizing the
de~Branges-Rovnyak expansion to the spaces $L_a$. We will
determine exactly the functions
$\cA_a^\etendu(s)$, $\cB_a^\etendu(s)$, $\cE_a^\etendu(s)$ and especially the
function $\mu^\etendu(a)$. It will be seen that this is a  more
complicated function than the simple-minded $\mu(a) =
2a$\dots 

The key now is to obtain the functions $p(a)$ and $q(a)$
defined in Proposition \ref{prop:pq}, and the function
$\cY_a(1,s)$. It turns out that their computation also
involves the quantities (we recall that $J_0^a(x) =
J_0(2\sqrt{ax})$):
\begin{subequations}
  \begin{align}
    r(a) &= 1 + \int_0^a ((1 - D_a)^{-1} H_a\cdot J_0^a)(x)
    \,dx\\
    s(a) &= \int_0^a ((1 - D_a)^{-1} \cdot J_0^a)(x) \,dx
  \end{align}
\end{subequations}
In order to compute $r$, $s$, $p$, $q$ we shall need the
already defined
functions $\phi_a^+$ ($=(1+H_a)^{-1}J_0^a$ on $(0,a)$),
$\phi_a^-$, ($=(1-H_a)^{-1}J_0^a$) as well as the entire functions
$\psi_a^+$ and  $\psi_a^-$ verifying:
\begin{subequations}
  \begin{align}\label{eq:111a}
    \psi_a^+ + \cH P_a \psi_a^+ &= 1\\ \label{eq:111b}
    \psi_a^- - \cH P_a \psi_a^- &= 1
  \end{align}
\end{subequations}
We have $r(a) = 1 + \frac12 \int_0^a (-\phi_a^+(x) +
\phi_a^-(x))\,dx$, and we know explicitely
$\phi_a^\pm$. But, we shall proceed in a more general
manner. First we recall the differential equations
\eqref{eq:110a}, \eqref{eq:110b} which are verified
by $\phi_a^\pm$ (where $\xdtx = x\frac\partial{\partial x}+\frac12$):
\begin{subequations}
\begin{align}
    \ada\phi_a^+ &= +\xdtx \phi_a^-  -(\mu(a) +
  \frac12) \phi_a^+\\ 
    \ada\phi_a^- &= +\xdtx \phi_a^+  + (\mu(a) -
  \frac12)\phi_a^-
  \end{align}
\end{subequations}
We compute $ar'(a)  = a\frac{-\phi_a^+(a) + \phi_a^-(a)}2 + \frac12
\int_0^a (x\frac\partial{\partial x}+1)(\phi_a^+(x) - \phi_a^-(x))
+\mu(a) (\phi_a^+(x)+\phi_a^-(x)) \,dx$ and simplifying this
gives exactly $ar'(a) = \mu(a) \frac12 \int_0^a
(\phi_a^+(x)+\phi_a^-(x)) \,dx = s(a)$. Similarly starting
with 
$s(a) = \frac12 \int_0^a
(\phi_a^+(x)+\phi_a^-(x)) \,dx $ we obtain $as'(a) = \frac12
\mu(a) + \frac12 \int_0^a
(x\frac\partial{\partial x}(\phi_a^+(x)+\phi_a^-(x)) + \mu(a)
(-\phi_a^+(x)+\phi_a^-(x)))\,dx$ which gives $\mu(a)r(a) -
s(a)$. So the quantities $r$ and $s$ verify the system:
\begin{subequations}
  \begin{align}
    ar'(a) &= \mu(a) s(a)\\
    (as)'(a) &= \mu(a) r(a)
  \end{align}
\end{subequations}
Either solving the system taking into account the behavior
as $a\to0$ or using the explicit formulas for $\phi_a^\pm$
we obtain in this specific instance of the study of $\cH$,
for which $\mu(a) = 2a$, that $r(a) = I_0(2a)$ and $s(a) =
I_1(2a)$. 

From \eqref{eq:111a} and \eqref{eq:111b} we obtain two types
of differential equations, either involving
$x\frac{\partial}{\partial x}$ or $a\frac{\partial}{\partial
a}$. 
From
$\psi_a^+(x) + \int_0^a J_0(2\sqrt{xy})\psi_a^+(x)\,dx = 1$, we
obtain $(1 + \cH P_a) a\frac{\partial}{\partial a}\psi_a^+(x)
= - a \psi_a^+(a) J_0^a$. We do similarly with $\psi_a^-$
and deduce:
\begin{subequations}
  \label{eq:109}
\begin{align}
  a\frac{\partial}{\partial a}\psi_a^+(x) &= - a\psi_a^+(a)
  \phi_a^+(x) \\
  a\frac{\partial}{\partial a}\psi_a^-(x) &= + a\psi_a^-(a) \phi_a^-(x)
\end{align}
\end{subequations}
Regarding the differential equations with
$x\frac{\partial}{\partial x}$, which we shall actually not
use, the computation is done using only the
fact that the kernel is a function of
$xy$ so $x\frac{\partial}{\partial x} J_0(2\sqrt{xy}) =
y\frac{\partial}{\partial y} J_0(2\sqrt{xy})$. We  only state
the result: 
\begin{subequations}
\begin{align}
  (x\frac{\partial}{\partial x}+\frac12)\psi_a^+(x) &=
  \frac12 \psi_a^-(x) - a\psi_a^+(a)  \phi_a^-(x) \\
  (x\frac{\partial}{\partial x}+\frac12)\psi_a^-(x) &=
  \frac12 \psi_a^+(x) + a\psi_a^-(a)  \phi_a^+(x) 
\end{align}
\end{subequations}
Let us now turn to the quantities $p(a)$ and $q(a)$. We have
$p(a) = \int_0^a (1 - D_a)^{-1}(\Un_{0<x<a})(x)\,dx =
\frac12 \int_0^a (\psi_a^+(x) + \psi_a^-(x))\,dx$. So $p'(a)
=  \frac12 (\psi_a^+(a) + \psi_a^-(a)) -\frac12 \psi_a^+(a)
\int_0^a \phi_a^+(x)\,dx + \frac12 \psi_a^-(a)\int_0^a
\phi_a^-(x)\,dx$. Reorganizing this gives:
\begin{equation}
  \label{eq:110}
  p'(a) = \frac{\psi_a^+(a) + \psi_a^-(a)}2 (1 +
  \int_0^a \frac{-\phi_a^+(x) + \phi_a^-(x)}2\,dx) +
  \frac{-\psi_a^+(a) + \psi_a^-(a)}2\int_0^a \frac{+\phi_a^+(x)
  + \phi_a^-(x)}2\,dx
\end{equation}
We remark that from the integral equations defining
$\psi_a^\pm$ we have $\psi_a^+(a) = 1 - \int_0^a
J_0(2\sqrt{ax})\psi_a^+(x)\,dx = 1 - \int_0^a
\phi_a^+(x)\,dx$ and $\psi_a^-(a) = 1 + \int_0^a
J_0(2\sqrt{ax})\psi_a^-(x)\,dx = 1 + \int_0^a
\phi_a^-(x)\,dx$. So $\frac{\psi_a^+(a) + \psi_a^-(a)}2 =
r(a)$ and $\frac{-\psi_a^+(a) + \psi_a^-(a)}2 = s(a)$. Hence
the quantity $p(a)$ verifies:
\begin{equation}
  \label{eq:111}
  p'(a) = r(a)^2 + s(a)^2
\end{equation}
With exactly the same method one obtains:
\begin{equation}
  \label{eq:112}
  q'(a) = 2 r(a) s(a)
\end{equation}
Let us observe that $q(a) = \frac12\int_0^a
((1-H_a)^{-1}-(1+H_a)^{-1})(1)\,dx = O(a^2)$ and $p(a) =
\frac12 \int_0^a((1+H_a)^{-1}+(1-H_a)^{-1})(1)\,dx
\sim_{a\to0} a$. So $(p\pm q)\sim_{a\to0} a$. Also
$r(a)\sim_{a\to0} 1$ and $s(a)\sim_{a\to0} a$. The equation
for $p(a)$ can be integrated:
\begin{equation}
  \label{eq:113}
  p(a) = a (r(a)^2 - s(a)^2)
\end{equation}
Indeed this has the correct derivative. Regarding $q(a)$ the
situation is different, one has $q' = 2rs = \frac{2a}\mu
rr'$ so in the special case considered here, and only in
that case we have $q(a) = \frac12 (r(a)^2 -1)$. Summing up:

\begin{proposition}
The quantities $r(a)$, $s(a)$, $p(a)$ and $q(a)$  verify the
differential equations $ar'(a) = \mu(a)s(a)$, $as'(a) + s(a)
= \mu(a)r(a)$,  $p'(a) = r(a)^2 +s(a)^2$, $q'(a) =
2r(a)s(a)$, $p(a) = a (r(a)^2 - s(a)^2)$. In the special
case of the $\cH$ transform one has:
\begin{subequations}
  \begin{align}
    r(a) &= I_0(2a)\\
    s(a) &= I_1(2a)\\
    p(a) &= a( I_0^2(2a) - I_1^2(2a))\\
    q(a) &= \frac12 (I_0^2(2a) - 1)
  \end{align}
\end{subequations}
\end{proposition}

We now need to determine $\cY_a(1,s) =
s(s-1)\Gamma(s)\wh{Y_1}(s)$. There holds $Y_1 = u_1 + \cH
v_1 = -\alpha \Un_{0<x<a} + \Un_{x>a} \cH v_1$. So
$\wh{Y_1}(s) = \alpha \frac{a^{1-s}}{s-1} + \int_a^\infty
(\cH v_1)(x)x^{-s}\,dx$. Then $\int_a^\infty (\cH v_1)(x)
x^{-s}\,dx = \int_0^\infty v_1(x)g_s(x)\,dx = \int_0^a
v_1(x)g_s(x)\,dx$, where the function $g_s$ from
\eqref{eq:55} has been used. Recalling from \eqref{eq:52a},
\eqref{eq:52b} the analytic functions $u_s$, equal to $- (1
- D_a)^{-1} H_a(g_s)$ on $(0,a)$, and $v_s$, equal to $(1 -
D_a)^{-1} P_a(g_s)$ on $(0,a)$, and using \eqref{eq:120} and
self-adjointness we obtain
\begin{equation}
  \label{eq:114}
  \int_0^a v_1(x)g_s(x)\,dx = - \alpha \int_0^a u_s(x)\,dx -
  \beta \int_0^a v_s(x)\,dx
\end{equation}
Let us now recall that we computed (\eqref{eq:57a}) $(\xdx +
s) u_s$ and found it to be on the interval $(0,a)$ given as $ -
av_s(a) (1 - D_a)^{-1}(J_0^a) - a(a^{-s} + u_s(a)) (1 -
D_a)^{-1}H_a(J_0^a)$. Integrating and also using equations
\eqref{eq:100} and \eqref{eq:101} we obtain
\begin{equation}
  \label{eq:115}
  \sqrt a\wh{E_a}(s)- a^{1-s} + (s-1)\int_0^a u_s(x)\,dx = - \sqrt a
  \wh{\cH(E_a)}(s) s(a) - \sqrt{a} \wh{E_a}(s) (r(a) - 1)
\end{equation}
\begin{equation}
  \label{eq:116}
  \int_0^a u_s(x)\,dx = \sqrt{a}\;\frac{a^{\frac12-s} -
  \wh{E_a}(s) r(a) - \wh{\cH(E_a)}(s)  s(a)}{s-1}
\end{equation}
We have similarly (\eqref{eq:57b}) $(\xdx + 1 - s)v_s =
  -\sqrt a\wh{E_a}(s) (1 - D_a)^{-1}(J_0^a) - \sqrt{a} \wh{\cH(E_a)}(s) (1 -
D_a)^{-1}H_a(J_0^a)$ so integration gives $av_s(a) - s \int_0^a
  v_s(x)\,dx = -\sqrt a\wh{E_a}(s) s(a) - \sqrt{a}
  \wh{\cH(E_a)}(s)(r(a) - 1)$ hence
  \begin{equation}
    \label{eq:117}
    \int_0^a v_s(x)\,dx = \sqrt a\;\frac{\wh{E_a}(s) s(a) +
    \wh{\cH(E_a)}(s) r(a)}{s}
  \end{equation}
Combining \eqref{eq:116}, \eqref{eq:117} with
  \eqref{eq:114}, and using $\cY_a(1,s) =
  s(s-1)\Gamma(s)\wh{Y_1^a}(s)$: 
\begin{subequations}
\begin{align}
\label{eq:118}
\wh{Y_1}(s) &= \sqrt{a}\, \wh{E_a}(s)
\Big(\frac{\alpha(a)r(a)}{s-1} - \frac{\beta(a)s(a)}{s}\Big)
+ \sqrt{a}\, \wh{\cH(E_a)}(s)\Big(\frac{\alpha(a)s(a)}{s-1}
- \frac{\beta(a)r(a)}{s}\Big)\\
\cY_a(1,s)  &= \sqrt{a}\, \Big( s\alpha(a) 
(\cE_a(s)r(a) + \cE_a(1-s)s(a)) + (1-s)\beta(a) (\cE_a(s)s(a) +
\cE_a(1-s)r(a))\Big) 
\end{align}
\end{subequations}

\begin{proposition}
  The functions $\cY_a(1,s)$ and $\cY_a(1,1-s)$ verify
  \begin{equation}
    \label{eq:119}
     \begin{bmatrix} \cY_a(1,s) \\ \cY_a(1,1-s) \end{bmatrix}
  = \sqrt{a}\;  \begin{bmatrix}
  \alpha(a)&\beta(a)\\\beta(a)&\alpha(a) \end{bmatrix} 
    \begin{bmatrix} s(\cE_a(s)r(a) + \cE_a(1-s)s(a)) \\ (1-s)(\cE_a(s)s(a) +
\cE_a(1-s)r(a)) \end{bmatrix}   
  \end{equation}
\end{proposition}

Comparing with equation
\eqref{eq:124} we get:  $T_a(s) = \sqrt{a}
s(\cE_a(s)r(a) + \cE_a(1-s)s(a))$. So:

\begin{subequations}
 \begin{theorem}\label{thm:cy}
  The analytic reproducing kernel $\cY_a(s,z)$ associated
  with the extended spaces $L_a$ is:
  \begin{equation}
    \cY_a(s,z) = s(s-1)z(z-1)\cX_a(s,z) + \begin{bmatrix}
    T_a(s) & T_a(1-s) \end{bmatrix}
    \begin{bmatrix} \alpha(a)&\beta(a)\\\beta(a)&\alpha(a) \end{bmatrix}
    \begin{bmatrix} T_a(z) \\ T_a(1-z) \end{bmatrix} 
   \end{equation}
\begin{align}
       \cX_a(s,z) &= \frac{\cE_a(s)\cE_a(z) -
    \cE_a(1-s)\cE_a(1-z)}{s+z-1} & \cE_a(s) &=
    2\sqrt{a}K_s(2a)\\
    T_a(s) &= \sqrt{a}\,
s\;(\cE_a(s)r(a) + \cE_a(1-s)s(a)) & r(a) &= I_0(2a)\quad
    s(a) = I_1(2a)\\
    \alpha(a) &= \frac{p(a)}{p(a)^2 - q(a)^2}& p(a) &= a
    (I_0^2(2a) - I_1^2(2a))\\ 
    \beta(a) &= \frac{q(a)}{p(a)^2 - q(a)^2}& q(a) &= \frac12 (I_0^2(2a) - 1)
  \end{align}
\end{theorem}
\end{subequations}

We proceed now to the determination of $\cA_a^\etendu$,
$\cB_a^\etendu$ and $\cE_a^\etendu = \cA_a^\etendu(s) - i \cB_a^\etendu(s)$. The
function $\cA_a^\etendu(s)$ is even under $s\to 1-s$ and
$\cB_a^\etendu(s)$ is odd. We must have:
\begin{equation}
  \label{eq:123}
  z \cY_a(1,z) = 2 (-i\cB_a^\etendu(1))\cA_a^\etendu(z) + 2\cA_a^\etendu(1)(-i\cB_a^\etendu(z))
\end{equation}
On the other hand from \eqref{eq:124} we have $\cY_a(1,z) =
\alpha T_a(z) + \beta T_a(1-z)$. Let us write
\begin{align}
  z\,T_a(z) = \sqrt{a}& (z(z-1) r(a)\cE_a(z) + z
                    r(a)\cE_a(z) \notag \\ & +z(z-1)
                    s(a)\cE_a(1-z) + z s(a)\cE_a(1-z) ) \\
  z\,T_a(1-z) = \sqrt{a}& (-z(z-1)s(a)\cE_a(z) -
                    z(z-1)r(a)\cE_a(1-z))
\end{align}
\begin{align}
 z\cY_a(1,z) = \sqrt{a}& \Big( z(z-1) ( (\alpha r -
                    \beta s) \cE_a(z) + (\alpha s - \beta r)
                    \cE_a(1-z))\notag\\
                       & + z \alpha (r\, \cE_a(z) +
                          s\cE_a(1-z)) \Big)
\end{align}
Extracting the even part $(z\cY_a(1,z))^+$ and the odd part $(z\cY_a(1,z))^-$:
\begin{equation}
  \label{eq:126}
  (z\cY_a(1,z))^+ = \sqrt{a}\Big( z(z-1)(\alpha
  -\beta)(r+s)\cA_a + (z -\frac12)\alpha (r-s)(-i \cB_a) +
  \frac12 \alpha (r+s) \cA_a\Big)
\end{equation}
\begin{equation}
  \label{eq:127}
  (z\cY_a(1,z))^- = \sqrt{a}\Big( z(z-1)(\alpha
  +\beta)(r-s)(-i\cB_a) + (z -\frac12)\alpha (r+s) \cA_a  +
  \frac12 \alpha (r-s) (-i\cB_a)\Big)
\end{equation}
We have $(z\cY_a(1,z))^+ = 2 (-i\cB_a^\etendu(1))\cA_a^\etendu(z)$ and
$(z\cY_a(1,z))^- = 2\cA_a^\etendu(1)(-i\cB_a^\etendu(z))$. Let us define
$K(a) = (2(-i\cB_a^\etendu(1)))^{-1}$ and $L(a) =
(2\cA_a^\etendu(1))^{-1}$. We know that:
\begin{equation}
  \label{eq:128}
  \lim_{\sigma\to+\infty}
  \frac{-i\cB_a^\etendu(\sigma)}{\cA_a^\etendu(\sigma)} = 1
\end{equation}
So it must be that 
\begin{equation}
  \label{eq:129}
  K(a) (\alpha  -\beta)(r+s) = L(a) (\alpha  +\beta)(r-s)
\end{equation}
Also, taking $z=1$ in \eqref{eq:126} we have $\frac1{KL} =
2\sqrt{a} \frac12\alpha\, (r \cE_a(1) + s \cE_a(0)) = \alpha
T_a(1)$. But referring to \eqref{eq:124} one has $T_a(1) =
p\alpha - q\beta = 1$. So: 
\begin{equation}
  \label{eq:130}
  K(a) L(a) = \frac1{\alpha(a)}
\end{equation}
Then:
\begin{equation}
  \label{eq:131}
  K(a)^2 = \frac1{\alpha(a)} \frac{(\alpha
  +\beta)(r-s)}{(\alpha -\beta)(r+s)} =\frac{\alpha^2 -
  \beta^2}{\alpha(a)} \frac{r^2-s^2}{(\alpha
  -\beta)^2(r+s)^2} = \frac1{p} \frac{p}a \frac1{(\alpha
  -\beta)^2(r+s)^2}
\end{equation}
\begin{equation}
  \label{eq:132}
  (\alpha - \beta)(r+s) K(a) = a^{-\frac12}
\end{equation}
We conclude:
\begin{subequations}
  \begin{align}\label{eq:150a}
    \cA_a^\etendu &= z(z-1) \cA_a +
  (z-\frac12)\frac{\alpha(r-s)}{(\alpha-\beta)(r+s)}(-i\cB_a)
  + \frac{\alpha}{2(\alpha-\beta)}\cA_a\\ \label{eq:150b}
    -i\cB_a^\etendu &= z(z-1) (-i\cB_a) +
  (z-\frac12)\frac{\alpha(r+s)}{(\alpha+\beta)(r-s)}\cA_a
  + \frac{\alpha}{2(\alpha+\beta)}(-i\cB_a)
  \end{align}
\end{subequations}
Let us now observe that $\frac\alpha{\alpha\pm\beta} =
  \frac{p}{p\pm q}$ and further:
\begin{subequations}
\begin{equation}
  \label{eq:133}
  \frac\alpha{\alpha-\beta} \frac{r-s}{r+s} = \frac{p}{p-q}
  \frac{a(r-s)^2}p = a\frac{p'-q'}{p-q} = a\frac{d}{da} \log(p-q)
\end{equation}
\begin{equation}
  \label{eq:134}
  \frac\alpha{\alpha+\beta} \frac{r+s}{r-s} = \frac{p}{p+q}
  \frac{a(r+s)^2}p = a\frac{p'+q'}{p+q} = a\frac{d}{da} \log(p+q)
\end{equation}
\end{subequations}
\begin{subequations}
  \begin{align}
    \cA_a^\etendu(z) &= ( z(z-1) + \frac12 \frac{p}{p-q} )\cA_a(z) +
  a\frac{d}{da} \log(p-q) (z-\frac12)(-i\cB_a(z))\\
    -i\cB_a^\etendu(z) &= ( z(z-1) +   \frac12 \frac{p}{p+q} ) (-i\cB_a(z)) +
   a\frac{d}{da} \log(p+q) (z-\frac12)\cA_a(z)
  \end{align}
\end{subequations}
Combining we get finally:

\begin{theorem}
  The E function associated with the entire
  functions $s(s-1)\Gamma(s)\wh f(s)$, $f\in L_a$ is:
  \begin{equation}
  \label{eq:135}
\begin{split}
  \cE_a^\etendu(z) = \Big( z(z-1) + \frac12 a\frac{d}{da} \log(p(a)^2
  - q(a)^2) (z-\frac12) + \frac12 p(a)\alpha(a)\Big)&
  \cE_a(z)\\
+ \Big( \frac12 a\frac{d}{da} \log\frac{p(a) + q(a)}{p(a)-
  q(a)} (z- \frac12) + \frac12 p(a) \beta(a) \Big)&
  \cE_a(1-z)
\end{split}
  \end{equation}
where $p(a) = a(I_0^2(2a)- I_1^2(2a))$, $q(a) =
  \frac12(I_0^2(2a)-1)$, $\alpha(a) = \frac{p(a)}{p(a)^2 -
  q(a)^2}$, $\beta(a) = \frac{q(a)}{p(a)^2 - q(a)^2}$, and
  $\cE_a(z) = 2\sqrt a K_z(2a)$.
\end{theorem}

We shall now obtain by two methods the function
  $\mu^\etendu(a)$. First, we compute $\cE_a^\etendu(\frac12) =
  (-\frac14 + \frac12 p (\alpha + \beta))\cE_a(\frac12) =
  \frac14\frac{p+q}{p-q}\cE_a(\frac12)$ and invoke
  $a\frac{d}{da} \cE_a^\etendu(\frac12) =  -
  \mu^\etendu(a)\cE_a^\etendu(\frac12)$. We thus have:

  \begin{theorem}
    The mu function for the chain of spaces $L_a$,
  $0<a<\infty$ is 
  \begin{align}
    \label{eq:136}
    \mu^\etendu(a) &= \mu(a) + a\frac{d}{da} \log\frac{p-q}{p+q}\\
             &= 2a + a\frac{d}{da} \log\frac{(2a-1) I_0^2(2a) - 2a
  I_1^2(2a) + 1}{(2a+1) I_0^2(2a) - 2a  I_1^2(2a) - 1}\\
             &= 2a - 2 + o(1)\qquad (a\to+\infty)
  \end{align}
  \end{theorem}

The asymptotic behavior is a corollary to $\lim_{a\to\infty}
  2a\frac{-pq' + qp'}{p^2 - q^2} = -2$ which itself follows
  from $p^2 - q^2 \sim \frac14 I_0(2a)^4\frac1{16 a^2}$ and
  $(\frac qp)'\sim +\frac1{16 a^3}$ which are easily deduced
  from the asymptotic expansion $I_0(x) =
  \frac{e^x}{\sqrt{2\pi x}}(1 + \frac1{8x} + \frac9{128
  x^2}+\dots)$ (\cite{watson}). Of course the $o(1)$ is in
  fact an $O(a^{-1})$.

The second method to obtain $\mu^\etendu(a)$ relies on 
$\frac{\cE_a^\etendu(1-\sigma)}{\cE_a^\etendu(\sigma)}\sim_{\sigma\to+\infty} 
    \frac{\mu^\etendu(a)}{2\sigma}$ (\eqref{eq:70}). We have:
    \begin{subequations}
      \begin{align}
      \frac{\cE_a^\etendu(\sigma)}{\sigma^2 \cE_a(\sigma)}
    &=  1 + \frac{\frac12 a\frac{d}{da} \log(p^2 - q^2) -
    1}\sigma + O(\frac1{\sigma^2})\\
 \frac{\cE_a^\etendu(1-\sigma)}{\sigma^2 \cE_a(1-\sigma)}
    &\to_{\sigma\to\infty}  1 - \frac12 (a\frac{d}{da}
    \log\frac{p+q}{p-q} ) \frac2{\mu(a)}
      \end{align}
    \end{subequations}
so $\frac{\mu^\etendu(a)}{\mu(a)} = 1 -  \frac1{\mu(a)} a\frac{d}{da}
    \log\frac{p+q}{p-q} $ and \eqref{eq:136} is
    confirmed. We can use this method to gather more
    information. From \eqref{eq:139} we have, as $\Re(s)\to+\infty$:
    \begin{subequations}
      \begin{align}
         \wh{E_a}(s) &= a^{\frac12 -s} (1 + \frac{a\phi^+(a)
    - a\phi^-(a)}{2s} + O(\frac1{s^2}))\\
 \wh{E_a^\etendu}(s) &= a^{\frac12 -s} (1 + \frac{a\phi^{\etendu+}(a)
    - a\phi^{\etendu-}(a)}{2s} + O(\frac1{s^2}))
      \end{align}
    \end{subequations}
Let us be careful  that $\cE_a(s) =
    \Gamma(s)\wh{E_a}(s)$ while $\cE_a^\etendu(s) = s^2 \Gamma(s)
    \wh{E_a^\etendu}(s)$.  We obtain:
    \begin{subequations}
      \begin{align}
      a\phi^{\etendu+}(a) - a\phi^{\etendu-}(a) &= a\phi^+(a) -
    a\phi^-(a) + a\frac{d}{da} \log(p^2 - q^2) - 2    \\ \label{eq:140b}
        a\phi^{\etendu+}(a) &= a\phi^+(a) + a\frac{d}{da} \log\frac{p-q}a
        \\ \label{eq:140c}
        a\phi^{\etendu-}(a) &= a\phi^-(a) - a\frac{d}{da}
        \log\frac{p+q}a
      \end{align}
    \end{subequations}
We recall that $(p\pm q)\sim_{a\to0} a$.  We integrate
\eqref{eq:140b} and \eqref{eq:140c} using \eqref{eq:141a},
\eqref{eq:141b} and this gives $\det(1+H_a^\etendu) =
\frac{p-q}a \det(1+H_a)$ and $\det(1 - H_a^\etendu) =
\frac{p+q}a \det(1 - H_a)$.
    \begin{subequations}
      \begin{align}
      \det(1+H_a^\etendu) &=\frac{p-q}a \det(1+H_a) = \det(1+H_a)\;\frac1a
      \int_0^a (r-s)^2\,da\\
\det(1 - H_a^\etendu) &= \frac{p+q} a\det(1 - H_a) = \det(1-H_a)\;\frac1a
      \int_0^a (r+s)^2\,da
      \end{align}
    \end{subequations}

    \begin{subequations}
\begin{theorem}
  Let $\cH^\etendu = L\,\cH\,L^{-1}$ be the self-reciprocal
  operator on $L^2(0,\infty;dx)$  with kernel:
  \begin{equation}
    J_0(2\sqrt{xy}) - 2 \frac{J_1(2\sqrt{xy})}{\sqrt{xy}} +
  \frac{1 - J_0(2\sqrt{xy})}{xy} = \sum_{n=0}^\infty (-1)^n
  \frac{n^2 x^n y^n}{(n+1)!^2}
  \end{equation}
  and let $H_a^\etendu$ be the restriction to $L^2(0,a;dx)$. Then:
  \begin{align}
    \det(1+H_a^\etendu) &= e^{+a-\frac12 a^2}\; \frac1a\int_0^a (I_0(2a) -
    I_1(2a))^2\,da  = e^{+a-\frac12 a^2}\big( I_0^2(2a) -
    I_1^2(2a) - \frac{I_0^2(2a) - 1}{2a}\big)\\
    \det(1-H_a^\etendu) &= e^{-a-\frac12 a^2}\;\frac1a\int_0^a (I_0(2a) +
    I_1(2a))^2\,da = e^{-a-\frac12 a^2}\big( I_0^2(2a) -
    I_1^2(2a) + \frac{I_0^2(2a) - 1}{2a}\big)
  \end{align}
\end{theorem}
    \end{subequations}

From theorem \ref{thm:cy} $\|\cY_{\frac12}\|^2 =
\frac1{16}\|\cX_{\frac12}\|^2  + 2(\alpha +
\beta)T_a(\frac12)^2$ and $T_a(\frac12) =
\frac12\sqrt{a}\,(r+s)\cE_a(\frac12)$. Also,
$(\alpha+\beta)(r+s)^2 = \frac1{p-q} (r+s)^2$. Furthermore
$\|\cY_{\frac12}\|^2 = \frac1{16} \Gamma(\frac12)^2
\|X_{\frac12}^\etendu\|^2 = \frac\pi{16}\|X_{\frac12}^\etendu\|^2$ and 
$\|\cX_{\frac12}\|^2 = \pi\|X_{\frac12}\|^2$. And also from \eqref{eq:71}
$\cE_a(\frac12) = \sqrt\pi
\;\frac{\det(1-H_a)}{\det(1+H_a)}$ and from theorem
\ref{thm:x12} one has $\| X_{\frac12}^a \|^2 = 2 \int_a^\infty
    \left(\det\frac{1-H_b}{1+H_b}\right)^2\,\frac{db}b$ and
    the analog holds for $X_{\frac12}^{\etendu a}$. Let us observe
    that $X_{\frac12}^{\etendu a} = - L X_{\frac12}^\times$
    so $\|X_{\frac12}^{\etendu a}\| = \|X_{\frac12}^{a\times}\|$. So 
    \begin{equation}
      \label{eq:125}
      \|X_{\frac12}^{a\times}\|^2 = 2 \int_a^\infty
    \left(\det\frac{1-H_b^\etendu}{1+H_b^\etendu}\right)^2\,\frac{db}b = 2 \int_a^\infty
    \left(\det\frac{1-H_b}{1+H_b}\right)^2\,\frac{db}b +
    8a\,\frac{p'+q'}{p-q}
    \left(\det\frac{1-H_a}{1+H_a}\right)^2 
    \end{equation}

\begin{subequations}
\begin{theorem}
Let $L_a$ be the Hilbert space of square integrable
functions on $f\in L^2(0,\infty;dx)$ such that both $f$ and
$\cH(f) = \int_0^\infty J_0(2\sqrt{xy})f(y)\,dy$ are
constant on $(0,a)$. Then the squared norm of the linear
form $f\mapsto \int_0^\infty \frac{f(x)}{\sqrt x}\,dx$ is
given by either one of the following two
expressions:
\begin{equation}
  \label{eq:121}
  2\int_a^\infty \left(\frac{(2b+1)I_0^2(2b) - 2b I_1^2(2b)
  -1}{(2b-1)I_0^2(2b) - 2b I_1^2(2b) + 1}\right)^2\,\frac{e^{-4b}}b\,db
\end{equation}
\begin{equation}
  \label{eq:122}
  = 2\int_a^\infty \frac{e^{-4b}}b\,db +
  2\frac{8a(I_0(2a)+I_1(2a))^2}{(2a-1)I_0^2(2a) - 2a I_1^2(2a)
  + 1} e^{-4a}
\end{equation}
The squared norm of the restriction of the linear form to
the subspace $K_a$ of functions vanishing on $(0,a)$ and
with $\cH(f)$ also vanishing on $(0,a)$ is $2\int_a^\infty
\frac{e^{-4b}}b\,db$.
\end{theorem}
\end{subequations}

One may express the wish to verify explicitely from
equations \eqref{eq:150a} and \eqref{eq:150b}, or in
the equivalent form
\begin{subequations}
  \begin{align}
    \cA_a^\etendu(z) &= \Big( (z-\frac12)^2
    +\frac14\frac{p+q}{p-q}\Big) \cA_a(z) +
  (a\frac{d}{da}\log(p-q))(z-\frac12)(-i\cB_a(z))\\
    -i\cB_a^\etendu(z) &= \Big( (z-\frac12)^2
    +\frac14\frac{p-q}{p+q}\Big)(-i\cB_a(z)) +
    (a\frac{d}{da}\log(p+q))(z-\frac12)\cA_a(z)
  \end{align}
\end{subequations}
the
differential system:
\begin{subequations}
  \begin{align}
    a\frac\partial{\partial a} \cA_a^\etendu(z) &= -
    \mu^\etendu(a)\cA_a^\etendu(z) - (z-\frac12)(-i\cB_a^\etendu(z))\\
    a\frac\partial{\partial a} (-i\cB_a^\etendu(z)) &= +
    \mu^\etendu(a)(-i\cB_a^\etendu(z)) - (z-\frac12)\cA_a^\etendu(z)
  \end{align}
\end{subequations}
and also to verify explicitely the reproducing kernel
formula
\begin{align}
  \cY_a(s,z) = \frac{\cE_a^\etendu(s)\cE_a^\etendu(z) -
    \cE_a^\etendu(1-s)\cE_a^\etendu(1-z)}{s+z-1} 
\end{align}
The interested reader will see that the algebra has a tendency to
become slightly involved if one does not benefit from the
following preliminary observations: using $p' = r^2 + s^2$,
$q' = 2rs$, $ar' = \mu r$, $as' = \mu r -s$, $p = a(r^2
-s^2)$ one first establishes $aq''+q' = 2\mu p'$, $ap''+p' -
\frac pa  = 2\mu q'$. Using this one checks easily:
\begin{subequations}
  \begin{align}
    \left(\frac{p'+q'}{p+q}\right)' +
    \left(\frac{p'+q'}{p+q}\right)^2 &= \frac1{a^2}\frac
    p{p+q} + \frac{2\mu - 1}a \frac{p'+q'}{p+q}\\
    \left(\frac{p'-q'}{p-q}\right)' +
    \left(\frac{p'-q'}{p-q}\right)^2 &= \frac1{a^2}\frac
    p{p-q} - \frac{2\mu + 1}a \frac{p'-q'}{p-q}
  \end{align}
Also the identity
\begin{equation}
  \frac{p^2}{p^2 - q^2} = a\frac{p'+q'}{p+q} a\frac{p' -
  q'}{p-q} = p\alpha
\end{equation}
is useful. The verifications may then be done. 
\end{subequations}


\section{Hyperfunctions in the study of the $\cH$ transform}
\label{sec:9}

In this final section we return to the equation
\eqref{eq:21}:
\begin{equation}
  \label{eq:18}
   \wt{\psi(f)}(i\,t) = \frac{t+1}{2t} \wt{f}(i\,\frac{t+\frac1t}2)\;, 
\end{equation}
Let us recall that $f\in L^2(0,\infty;dx)$ and
$\psi:L^2(0,\infty;dx)\to L^2(0,\infty;dx)$ is the isometry
which corresponds to $F(w)\mapsto F(w^2)$, where
$F(w) = \sum_{n=0}^\infty c_n w^n$, $f(x) =
\sum_{n=0}^\infty c_n P_n(x)e^{-x}$, $P_n(x) =
L_n^{(0)}(2x)$. Let $g = \psi(f)$. Using $\lambda =
it$, in the
$L^2$ sense:
\begin{equation}
  \label{eq:67}
  g(x) = \frac1{2\pi} \int_{-\infty}^\infty \frac{\lambda +
  i}{2\lambda}\wt
  f(\frac{\lambda-\frac1\lambda}2) e^{-i\lambda x}\,d\lambda
\end{equation}
It is natural to consider separately $\lambda>0$ and
$\lambda<0$. So let us define:
\begin{subequations}
\begin{align}
  \label{eq:68}
  G_+(x) &= \frac1{2\pi} \int_{-\infty}^0 \frac{\lambda +
  i}{2\lambda}\wt
  f(\frac{\lambda-\frac1\lambda}2) e^{-i\lambda x}\,d\lambda\\
  \label{eq:69}
  G_-(x) &= - \frac1{2\pi} \int_{0}^\infty  \frac{\lambda +
  i}{2\lambda}\wt
  f(\frac{\lambda-\frac1\lambda}2) e^{-i\lambda x}\,d\lambda
\end{align}  
\end{subequations}
We observe that $G_+$ is in the Hardy space of $\Im(x)>0$
and $G_-$ is in the Hardy space of $\Im(x)<0$. Their boundary
values must coincide on $(-\infty,0)$ as $g\in
L^2(0,+\infty;dx)$. So we have a single analytic function
$G(z)$ on $\CC\setminus[0,+\infty)$ with $G = G_+$ for
$\Im(x)>0$ and $G = G_-$ for $\Im(x)<0$. Then $g = \psi(f) =
G_+ - G_-$ is computed as
\begin{equation}
  \label{eq:140}
  g(x) = G(x+i0) - G(x-i0)
\end{equation}
In other words $g$ is most naturally seen as a hyperfunction
\cite{morimoto}, as a difference of boundary values of
analytic functions. We shall now compute it explicitely, and
also we will show later that this observation extends to the
distributions $A_a(x)$, $-iB_a(x)$, $E_a(x)$ which are
associated with the study of the $\cH$ transform. The point
of course is that the corresponding functions $G$ will for
them have a simple natural expression.

\begin{subequations}
We have, for $\Im(z)>0$:
\begin{equation}
  \label{eq:152}
  G(z) = \frac1{2\pi} \int_0^{\infty} \frac{\lambda -
  i}{2\lambda}\wt
  f(\frac{-\lambda+\frac1\lambda}2) e^{+i\lambda z}\,d\lambda
\end{equation}
\begin{equation}
  \label{eq:153}
   G(z) = \frac1{2\pi} \int_0^{\infty} \frac{\lambda -
  i}{2\lambda}\left(\int_0^\infty 
  e^{i\frac12 x(-\lambda+\frac1\lambda)} f(x)\,dx\right)
  e^{+i\lambda z}\,d\lambda 
\end{equation}
Let $\mu = \frac12 (\lambda - \frac1\lambda)$, $\lambda =
  \mu + \sqrt{1+\mu^2}$, $\frac{\lambda -
  i}{2\lambda}\,d\lambda = \frac{\lambda}{\lambda +
  i}\,d\mu$,  with, for $0<\lambda<\infty$, $-\infty<\mu<\infty$.
\begin{equation}
  \label{eq:160}
   G(z) = \frac1{2\pi} \int_{-\infty}^{\infty} \frac{\lambda}{\lambda +
  i}\left(\int_0^\infty 
  e^{- i\mu x} f(x)\,dx\right)
  e^{+i\lambda z}\,d\mu  
\end{equation}
For $\mu\to+\infty$, $\lambda = 2\mu + \frac1{2\mu} +
  \dots$, and for $\mu\to-\infty$, $\lambda = -\frac1{2\mu} +
  \dots$, and $\frac{\lambda}{\lambda +
  i} \sim \frac{i}{2\mu}$ as $\mu\to-\infty$. So far the
  inner integral is in the $L^2$ sense. We shall now 
  suppose that $f$ and $f'$ are in $L^1$ (so  $\lim_{x\to\infty} f(x)
  = 0$)  and write $\int_0^\infty 
  e^{- i\mu x} f(x)\,dx  = \int_0^\infty 
  e^{- i\mu x - x} e^{x}f(x)\,dx =  \frac{f(0)}{i\mu +1} +
  \frac1{i\mu +1}
  \int_0^\infty e^{- i\mu x} (f(x)+f'(x))\,dx$. 
  \begin{equation}
    \label{eq:161}
    G(z) = \frac1{2\pi} \left( f(0) \int_{-\infty}^{\infty}
  \frac{\lambda}{\lambda + 
  i}\frac{e^{+i\lambda z}}{i\mu + 1}
  \,d\mu  + \int_{-\infty}^{\infty} \frac{\lambda}{\lambda +
  i}\frac{e^{+i\lambda z}}{i\mu +1}\left(\int_0^\infty 
  e^{- i\mu x} (f(x)+f'(x))\,dx\right)
  \,d\mu  \right)
  \end{equation} 
In this manner, with $f\in L^1$, $f'\in
L^1$, $\Im(z)>0$, we have an absolutely convergent double integral.
  \begin{equation}
    \label{eq:162}
    G(z) = \frac1{2\pi} \left( f(0) \int_{-\infty}^{\infty}
  \frac{\lambda}{\lambda + 
  i}\frac{e^{+i\lambda z}}{i\mu + 1}
  \,d\mu  + \int_0^\infty 
  \left(\int_{-\infty}^{\infty} \frac{\lambda}{\lambda +
  i}\frac{e^{-i\mu x  +i\lambda z}}{i\mu +1} \,d\mu\right)
  (f(x)+f'(x))\,dx\right) 
 \end{equation}
Observing $\frac1{2\pi} \int_{-\infty}^{\infty} \frac{\lambda}{\lambda +
  i}\frac{e^{+i\lambda z}}{i\mu + 1} \,d\mu = \frac1{2\pi}\int_0^\infty
  \frac{\lambda-i}{2\lambda} \frac1{i\mu +1} e^{+i\lambda
  z}\,d\lambda = \frac1{2\pi i}\int_0^\infty \frac{1}{\lambda-i}
  e^{+i\lambda z}\,d\lambda$, we then
suppose $\Re(z)<0$, $\Im(z)>0$ (or
$\Im(z)\geq0$) so that we may rotate the contour to
$\lambda=-it$, $0\leq t<\infty$. This procedure gives thus:
  \begin{equation}
    \frac1{2\pi} \int_{-\infty}^{\infty} \frac{\lambda}{\lambda +
  i}\frac{e^{+i\lambda z}}{i\mu + 1} \,d\mu  = \frac1{2\pi
  i}\int_0^\infty \frac{e^{tz}}{1+t}\,dt 
  \end{equation}
Also:
\begin{equation}
  \label{eq:166}
    \frac1{2\pi} \int_{-\infty}^{\infty} \frac{\lambda}{\lambda +
  i}\frac{e^{-i\mu x +i\lambda z}}{i\mu + 1} \,d\mu =
  \frac1{2\pi i}\int_0^\infty \frac{1}{\lambda-i} 
  e^{-i\frac{x}2(\lambda - \frac1\lambda)+i\lambda z}\,d\lambda 
\end{equation}
We rotate the contour to $\lambda\in i[0,-\infty)$,
  which is licit as $x\geq0$ and, for $\Re(z)<0$, $x\geq0$, we obtain:
  \begin{equation} 
    \frac1{2\pi i}\int_0^\infty \frac{e^{zt -
  \frac{x}2(t+\frac 1t)}}{1+t}\,dt
  \end{equation}
Going back this allows to write \eqref{eq:162}, for
$\Re(z)<0$, $\Im(z)>0$ as:
\begin{equation}
  \label{eq:165}
  G(z) = \frac1{2\pi i} \left( f(0) \int_0^\infty
  \frac{e^{zt}}{1+t}\,dt + \int_0^\infty 
  \left(\int_0^\infty \frac{e^{zt - \frac{x}2(t+\frac
  1t)}}{1+t}\,dt\right) (f(x)+f'(x))\,dx\right)
\end{equation}
and finally, after integrating by parts:
\begin{equation}
  \label{eq:167}
  G(z) = \frac1{2\pi i} \int_0^\infty \left(\int_0^\infty
  \frac12(1+\frac1t) e^{zt - \frac12 y(t+\frac
  1t)}\,dt\right) f(y)\,dy
\end{equation}
This last expression (still temporarily under the hypothesis
$f,f'\in L^1$) is certainly a priori absolutely convergent
for $\Re(z)<0$ and gives $G(z)$ in this half-plane.
\end{subequations}

\begin{subequations}
We are led to the study of:
\begin{equation}
  \label{eq:159}
  a(z,y) = \frac1{2\pi i} \int_0^\infty \frac12(1+\frac1t)
  e^{zt - \frac12 y(t+\frac 1t)} \,dt
\end{equation}
We still temporarily assume $\Re(z)<0$. We even suppose
$z<0$ and make a change of variable:
\begin{equation}
  \label{eq:164}
  a(z,y) =  \frac1{2\pi i}\left( \sqrt{\frac y{y-2z}}\; \frac12\int_0^\infty
  e^{-\frac12\sqrt{y(y-2z)}(u+\frac1u)}\,du + \frac12\int_0^\infty
  e^{-\frac12\sqrt{y(y-2z)}(u+\frac1u)}\frac1u\,du \right)
\end{equation}
\begin{equation}
  \label{eq:168}
  a(z,y) =  \frac1{2\pi i}\left( \sqrt{\frac 1{y-2z}}\; \frac12\int_0^\infty
  e^{-\frac12\sqrt{y-2z}(v+ y \frac1v)}\,dv + \frac12\int_0^\infty
  e^{-\frac12\sqrt{y-2z}(v+ y \frac1v)}\frac1v\,dv \right)
\end{equation}
\begin{equation}
  \label{eq:154}
   a(z,y) =  \frac1{2\pi i}\left( \sqrt{\frac y{y-2z}}\;
   K_1(\sqrt{ y(y-2z)}) + K_0(\sqrt{y(y-2z)})\right) 
\end{equation}
\end{subequations}
For any $z\in \CC\setminus[0,+\infty)$ and any $y\geq 0$ the
integrals in \eqref{eq:168} converge absolutely and define
an analytic function of $z$. Furthermore the $K$ Bessel
functions decrease exponentially as $y\to+\infty$ in
\eqref{eq:154}. For fixed $z$, $a(z,y)$ is certainly a
square-integrable function of $y$ (also at the origin),
locally uniformly in $z$ so the equation \eqref{eq:167}
defines $G$ as an analytic function on the entire domain
$\CC\setminus[0,+\infty)$. Then by an approximation argument
\eqref{eq:167} applies to any $f\in L^2(0,\infty;dx)$ and
any $z\in \CC\setminus[0,+\infty)$.

We now study the boundary values $a(x+i0,y)$,
$a(x-i0,y)$, $x,y\geq0$. We could use the expression of
the $K$ Bessel functions in terms of the Hankel functions
$H^{(1)}$ and $H^{(2)}$, go to the boundary, and then recover the
Bessel functions $J_0$ and $J_1$. But we shall proceed in a
more direct manner. Let us first examine
\begin{subequations}
\begin{equation}
  \label{eq:169}\begin{split}
  d(z,y) &= \frac1{2\pi i} \frac12 \int_0^\infty e^{zt -
  \frac12 y(t+\frac 1t)} \, \frac1t\,dt \\ &= \frac1{2\pi i}
  \frac12\int_0^\infty
  e^{-\frac12\sqrt{y(y-2z)}(u+\frac1u)}\frac1u\,du =
  \frac1{2\pi i} \int_1^\infty
  e^{-\sqrt{y(y-2z)}\,t}\frac{dt}{\sqrt{t^2-1}}\end{split}
\end{equation}
\begin{equation}
  \label{eq:170}
  = \frac1{2\pi i} \int_1^\infty
  e^{-\sqrt{y(y-2z)}\,t}\;t^{-1}\,dt + \frac1{2\pi i}
  \int_1^\infty
  e^{-\sqrt{y(y-2z)}\,t}(\frac1{\sqrt{t^2-1}} - \frac1t)\,dt
\end{equation}
\begin{equation}
  \label{eq:171}
  = \frac1{2\pi i}
  \frac{e^{-\sqrt{y(y-2z)}} - \int_1^\infty
  e^{-\sqrt{y(y-2z)}\,t}\,t^{-2}\,dt}{\sqrt{y(y-2z)}} + \frac1{2\pi i}
  \int_1^\infty
  e^{-\sqrt{y(y-2z)}\,t}(\frac1{\sqrt{t^2-1}} - \frac1t)\,dt
\end{equation}
\end{subequations}
We now look at the (distributional) boundary values $z\to x$ with $z=
x+i\epsilon$, $\epsilon\to 0^+$ or $z= x-i\epsilon$ and
$\epsilon \to 0^+$. We shall take $x>0$. Here the
singularities at $y=2x$ and at $y=0$ are integrable and we
need only take the limit in the naive sense. We distinguish
$y>2x$ from $0<y<2x$. In the former case, nothing happens:
\begin{subequations}
\begin{equation}
  \label{eq:172}
  d(x+i0,y) = d(x-i0,y) = \frac1{2\pi i} \int_1^\infty
  e^{-\sqrt{y(y-2x)}\,t}\frac{dt}{\sqrt{t^2-1}} 
\end{equation}
In the latter case:
\begin{equation}
  \label{eq:173}
  d(x+i0,y) = \frac1{2\pi}
  \frac{e^{+i\sqrt{y(2x-y)}} - \int_1^\infty
  e^{+i\sqrt{y(2x-y)}\,t}\,t^{-2}\,dt}{\sqrt{y(2x-y)}} + \frac1{2\pi i}
  \int_1^\infty
  e^{+i\sqrt{y(2x-y)}\,t}(\frac1{\sqrt{t^2-1}} - \frac1t)\,dt
\end{equation}
\begin{equation}
  \label{eq:174}
  d(x-i0,y) = -\frac1{2\pi}
  \frac{e^{-i\sqrt{y(2x-y)}} - \int_1^\infty
  e^{-i\sqrt{y(2x-y)}\,t}\,t^{-2}\,dt}{\sqrt{y(2x-y)}} + \frac1{2\pi i}
  \int_1^\infty
  e^{-i\sqrt{y(2x-y)}\,t}(\frac1{\sqrt{t^2-1}} - \frac1t)\,dt
\end{equation}
So $d(x+i0,y)- d(x-i0,y)$ is supported in $(0,2x)$ and has
values there
\begin{equation}
  \label{eq:175}
  \frac1{\pi}
  \frac{\cos\sqrt{y(2x-y)} - \int_1^\infty
  \cos(\sqrt{y(2x-y)}\,t)\,t^{-2}\,dt}{\sqrt{y(2x-y)}} + \frac1{\pi}
  \int_1^\infty
  \sin(\sqrt{y(2x-y)}\,t)(\frac1{\sqrt{t^2-1}} - \frac1t)\,dt
\end{equation}
\end{subequations}
We used this method to have a clear control not only of the
pointwise behavior but also of the limit as a
distribution. There is no necessity now to keep working with
absolutely convergent integrals and we have the simple
result, using the very classical Mehler formula:\footnote{we
are mainly  interested in the boundary value as a
distribution and we skip the discussion of the pointwise
behavior at the borders $y=0$ and $y=2x$.}
\begin{equation}
  \label{eq:176}
  d(x+i0,y)-d(x-i0,y) = \Un_{0<y<2x}(y) \frac1{\pi}
  \int_1^\infty
  \frac{\sin(\sqrt{y(2x-y)}\,t)}{\sqrt{t^2-1}}\,dt =
  \frac12\, \Un_{0<y<2x}(y) J_0(\sqrt{y(2x-y)})
\end{equation}

Let us now consider the behavior of 
\begin{equation}
  \label{eq:177}
  e(z,y) = \frac1{2\pi i} \frac12 \int_0^\infty e^{zt -
  \frac12 y(t+\frac 1t)} \, dt =\frac1{2\pi
  i}\sqrt{\frac y{y-2z}}\; \frac12\int_0^\infty
  e^{-\frac12\sqrt{y(y-2z)}(u+\frac1u)}\,du
\end{equation}
We make the simple observation that $e(z,y) =
\frac{\partial}{\partial z} d(z,y)$. So we shall have (as is
confirmed by a more detailed examination):
\begin{equation}
  \label{eq:178}\begin{split}
  e(x+i0,y)-e(x-i0,y) &= \frac{\partial}{\partial x}
  \frac12\, \Un_{0<y<2x}(y) J_0(\sqrt{y(2x-y)}) \\ &=
  \delta_{2x}(y) - \frac12\, \Un_{0<y<2x}(y)
  \sqrt{\frac{y}{2x-y}} J_1(\sqrt{y(2x-y)})
  \end{split}
\end{equation}

Combining all those elements we obtain that the function
$k(x) = \psi(f)(x)$ is given as:
\begin{equation}
  \label{eq:179}
  k(x) = f(2x) + \frac12 \int_0^{2x}
  J_0(\sqrt{y(2x-y)})f(y)\,dy - \frac12 \int_0^{2x}
  \sqrt{\frac{y}{2x-y}} J_1(\sqrt{y(2x-y)}) f(y)\,dy
\end{equation}
Some pointwise regularity of $f$ at $x$ is necessary to
fully justify the formula; in order to check if continuity
of $f$ at $2x$ is enough we can not avoid examining $e(z,y)$
more closely as $z\to x$.  
\begin{equation}
  \label{eq:181}\begin{split}
    e(z,y) &= \sqrt{\frac{y}{y-2z}}\frac1{2\pi i}
    \int_1^\infty
    e^{-\sqrt{y(y-2z)}\,t}\frac{t\,dt}{\sqrt{t^2-1}} \\ &=
     \frac1{2\pi i}\frac{e^{-\sqrt{y(y-2z)}}}{y -2z} +
    \sqrt{\frac{y}{y-2z}}\frac1{2\pi i} \int_1^\infty
    e^{-\sqrt{y(y-2z)}\,t}\;\frac{t -
    \sqrt{t^2-1}}{\sqrt{t^2-1}}\,dt\end{split}
\end{equation}
The integral term on the right causes no problem at all. And
    writing $\frac{e^{-\sqrt{y(y-2z)}}}{y -2z} =
    \frac1{y-2z} + \frac{e^{-\sqrt{y(y-2z)}}-1}{y -2z}$,
    again the term on the right has no problem, so there
    only remains $\frac1{y-2z}$, and of course, this is very
    well-known, the difference between $+i0$ and $-i0$ gives
    the Poisson kernel, so for non-tangential convergence,
    continuity of $f$ at $2x$ is enough. Of course this
    discussion was quite superfluous if we wanted to
    understand $k$ as an $L^2$ function, here we have the
    information that non tangential boundary value of
    $G(x+i0)-G(x-i0)$ does give pointwise the formula
    \eqref{eq:179} if $f$ is continuous at $y=2x$. We can
    also
rewrite \eqref{eq:179} as:
\begin{equation}
  \label{eq:180}
  k(x) = (1 + \frac{d}{dx}) \frac12 \int_0^{2x}
  J_0(\sqrt{y(2x-y)})f(y)\,dy
\end{equation}
This is exactly one half of equation \eqref{eq:9c}, where
$k$ was obtained from $(f,g)$ as
$\psi(f)+w\cdot\psi(g)$. Let us observe that $w =
\frac{\lambda - i}{\lambda + i}$ verifies, as an operator,
$(\frac{d}{dx} + 1)\cdot w = w\cdot(\frac{d}{dx} + 1)=
\frac{d}{dx} - 1$. So the isometry corresponding to
$g(w)\mapsto w\,G(w^2)$, which is the composite $w\cdot
\psi$, sends $g$ to  $(-1 + \frac{d}{dx}) \frac12 \int_0^{2x}
  J_0(\sqrt{y(2x-y)})f(y)\,dy$. This is indeed the second
half of equation \eqref{eq:9c}.

The formulas \eqref{eq:9a} and \eqref{eq:9b} may be
established in an exactly analogous manner (taking $k$ with
compact support to simplify the discussion). But this would
be a repetition of the arguments we just went
through,
so rather I will conclude the paper with
a method allowing to go directly from $\cA_a(s)$,
$-i\cB_a(s)$, $\cE_a(s)$ to the distributions $A_a(x)$,
$-iB_a(x)$, $E_a(x)$, and this will show that they are in a
natural manner (differences of) boundary values of an
analytic function.

From the expression $\cE_a(s) = \Gamma(s)\wh{E_a}(s) =
2\sqrt{a} K_s(2a) = \sqrt{a} \int_0^\infty e^{-a(t+\frac1t)}
t^{s-1}\,dt$, we shall recover $\wh{E_a}(s)$ as a right
Mellin transform with the help of the
Hankel formula $\Gamma(s)^{-1} = \int_\cC e^{v} v^{-s}\,dv$,
where $\cC$ is a contour coming from $-\infty$ along the
lower edge of the cut along $(-\infty,0]$ turning
counterclockwise around the origin and going back to $-\infty$
along or slightly above the upper edge of the cut. Let us
write the Hankel formula as
\begin{equation}
  \label{eq:182}
  \frac{t^{s-1}}{\Gamma(s)} = \frac1{2\pi i}\,\int_\cC
  e^{tv} v^{-s}\,dv\qquad (t>0)
\end{equation}
So we have:
\begin{equation}
  \label{eq:183}
  \wh{E_a}(s) = \sqrt{a} \frac1{2\pi i} \int_0^\infty
  \left(\int_\cC e^{tv} v^{-s}\,dv\right)e^{-a(t+\frac1t)} \,dt
\end{equation}
Let us suppose $\Re(s)>1$. Then the contour $\cC$ can be
deformed into the contour $\cC_\epsilon$ coming from $-i\infty$ to
$-i\epsilon$, then turning counterclockwise from
$e^{-i\frac\pi2}\epsilon$ to $e^{i\frac\pi2}\epsilon$, then
going to $+i\infty$. Also we impose $0<\epsilon<a$. The
integrals may then be permuted:
\begin{equation}
  \label{eq:184}
  \wh{E_a}(s) = \sqrt{a} \frac1{2\pi i} \int_{\cC_\epsilon}
  \left(\int_0^\infty e^{tv} e^{-a(t+\frac1t)} \,dt\right)
  v^{-s}\,dv
\end{equation}
and using $e(z,y)$ from \eqref{eq:177} this gives:
\begin{equation}
  \label{eq:185}
  \Re(s)>1\implies\quad \wh{E_a}(s) = \sqrt{a} \,
  \int_{\cC_\epsilon} 2e(v,2a) v^{-s}\,dv
\end{equation}
We have previously studied $e(z,y)$, which is also expressed
  as in \eqref{eq:181}. We see on this basis and simple
  estimates that we may deform $\cC_\epsilon$ into a contour
  $\cC_{a,\eta}$ going from $+\infty$ to $a+\eta$ along the lower
  border, turning clockwise around $a$ from $a+\eta -i0$ to
  $a+\eta+i0$, then going from $a+\eta$ to $+\infty$ on the
  upper border ($\eta\ll1$). We will have in particular from
  $\eqref{eq:181}$ a term $\frac2{2\pi
  i}\int_{a+\eta-i0}^{a+\eta+i0} \frac{v^{-s}}{2a-2v}\,dv$
  which is $a^{-s}$. The final result is obtained:
  \begin{equation}
    \label{eq:186}
    \Re(s)>1\implies\quad \wh{E_a}(s) = \sqrt{a}\left(
  a^{-s} - \int_a^\infty  \sqrt{\frac{a}{x-a}}
  J_1(2\sqrt{a(x-a)})x^{-s}\,dx\right) 
  \end{equation}
This identifies $\wh{E_a}(s)$ as the right Mellin transform
  of the distribution
  \begin{equation}
    \label{eq:187}
    E_a(x) = \sqrt{a}\left( \delta_a(x) +
  \Un_{x>a}(x)\frac{\partial}{\partial x}
  J_0(2\sqrt{a(x-a)}) \right) = \sqrt{a}
  \frac{\partial}{\partial x} \left(
  \Un_{x>a}(x)J_0(2\sqrt{a(x-a)})\right)
  \end{equation}
This proof reveals that the distribution $E_a(x)$
  is expressed in a natural manner as the difference of
  boundary values $\sqrt{a} (2 e(x+i0,2a) - 2 e(x-i0,2a))$,
  with
  \begin{equation}
    \label{eq:190}
    \sqrt{a}\,2e(z,2a) = \sqrt{a}\,\frac1{2\pi i} \int_0^\infty e^{z t -
 a(t+\frac 1t)} \, dt =   \sqrt{a}\,\frac1{2\pi i} 2\sqrt{\frac a{a-z}}\; 
   K_1(2\sqrt{ a(a-z)})
  \end{equation}

The formulas \eqref{eq:200a} and \eqref{eq:200b} are
  recovered in the same manner. 
  \begin{theorem}
The distribution $A_a(x) =
  \frac{\sqrt{a}}2(1+\cH)(\phi_a^+\Un_{0<x<\infty})$,
  $\phi_a^+(x)+\int_0^a J_0(2\sqrt{xy})\phi_a^+(y)\,dy =
  J_0(2\sqrt{ax})$, is
  the difference of boundary values $\sqrt{a} (a(x+i0,2a)
  - a(x-i0,2a))$, with:
  \begin{equation}
    \label{eq:188}\begin{split}
    \sqrt{a}\, a(z,2a) &= \sqrt{a}\,\frac1{2\pi i} \int_0^\infty
  \frac12(1+\frac1t) e^{z t - a(t+\frac 1t)} \, dt \\&=
  \sqrt{a}\,\frac1{2\pi i}\left( \sqrt{\frac a{a-z}}\; 
   K_1(2\sqrt{ a(a-z)}) + K_0(2\sqrt{a(a-z)})\right) \end{split}
  \end{equation}
The distribution $-iB_a(x) =
  \frac{\sqrt{a}}2(-1+\cH)(\phi_a^-\Un_{0<x<\infty})$,
  $\phi_a^-(x)-\int_0^a J_0(2\sqrt{xy})\phi_a^-(y)\,dy =
  J_0(2\sqrt{ax})$,  is the difference of
  boundary values $\sqrt{a}(-ib(x+i0,2a) - (-ib(x-i0,2a)))$,
  with:
  \begin{equation}
    \label{eq:189}\begin{split}
\sqrt{a}(-ib(z,2a)) &= \sqrt{a}\,\frac1{2\pi i} \int_0^\infty
  \frac12(1-\frac1t) e^{z t - a(t+\frac 1t)} \,
  dt    \\ &=
  \sqrt{a}\,\frac1{2\pi i}\left( \sqrt{\frac a{a-z}}\; 
   K_1(2\sqrt{ a(a-z)}) - K_0(2\sqrt{a(a-z)})\right) \end{split}
  \end{equation}
  \end{theorem}

\section{Appendix: a remark on the resolvent of the
Dirichlet kernel}

In this paper we have studied a special transform on the
positive half-line with a kernel of a multiplicative type
$k(xy)$, following the method summarized in
\cite{cras2002,cras2003}. We have associated to the kernel the
investigation of its Fredholm determinants on finite
intervals $(0,a)$, and have related them with first and
second order differential equations leading to problems of
spectral and scattering theory. There is a vast literature
on kernels of the additive type $k(x-y)$, and on the related
Fredholm determinants on finite intervals. The Dirichlet
kernel on $L^2(-s,s;dx)$:
\begin{equation}
  K_s(x,y) = \frac{\sin(x-y)}{\pi(x-y)}
\end{equation}
has been the subject of many works (only a few references
will be mentioned here.) The Fredholm determinant $\det(1 -
K_s)$, as a function of $s$ (or more generally as a function
of the endpoints of finitely many intervals), has many
properties, and is related to the study of random matrices
\cite{mehta}. The Fredholm determinants of the
even and odd parts
\begin{equation}
  K_s^\pm(x,y) = \frac{\sin(x-y)}{\pi(x-y)} \pm \frac{\sin(x+y)}{\pi(x+y)}
\end{equation}
on $L^2(0,s;dx)$ have been studied by Dyson \cite{dyson}. He
used the second derivatives of their logarithms to construct
potentials for Schrödinger equations on the half-line, and
studied their asymptotics with the tools of scattering
theory. Jimbo, Miwa, M\^ori, and Sato \cite{jmms} related
$\det(1 - K_s)$ to a Painlevé equation. Widom \cite{widom}
obtained the leading asymptotics using the Krein continuous
analog of orthogonal polynomials. Deift, Its, and Zhou
\cite{deiftitszhou} justified  the Dyson
asymptotic expansions using tools developed for
Riemann-Hilbert problems.  Tracy and Widom \cite{tracywidom}
established partial differential equations for the Fredholm
determinants of integral operators arising in the study of
the scaling limit of the distribution functions of
eigenvalues of random matrices. We refer the reader to the
cited references and to \cite{deiftitskrazhou} for recent
results and we apologize for not providing any more detailed
information here.

We have, in the present paper, been talking a lot of
scattering and determinants and one might wonder if this is
not a re-wording of known things. In fact, our work is with
the multiplicative kernels $k(xy)$, and (direct) reduction
to additive kernels would lead to (somewhat strange)
$g(t+u)$ kernels on semi-infinite intervals
$(-\infty,\log(a)]$. So we are indeed doing something
different; one may also point out that the entire functions
arising in the present study are not of finite exponential
type; and the scattering matrices do not at all tend to $1$
as the frequency goes to infinity. In the case of the cosine
and sine kernels the flow of information will presumably go
from the additive to the multiplicative, as the additive
situation is more flexible, and has stimulated the
development of powerful tools, with relation to Painlevé
equations, Riemann-Hilbert problems, Integrable systems
\cite{deiftitszhou}.

Nevertheless, one may ask if the framework of reproducing
kernels in Hilbert spaces of entire functions also may be
used in the additive situation. This is the case indeed and
it is very much connected to the method of Krein in inverse
scattering theory, and his continuous analog of orthogonal
polynomials (used by Widom in the context of the Dirichlet
kernel in \cite{widom}.) The Gaudin identities for
convolution kernels (\cite[App. A16]{mehta}) play a rôle
very analogous to the identities in the present paper
\eqref{eq:300a}, \eqref{eq:300b} involved in the study of
multiplicative kernels.  Widom in his proof \cite{widom} of
the main term of the asymptotics as $s\to+\infty$ studied
the Krein functions associated with the complement of the
interval $(-1,+1)$ and he mentioned the interest of extremal
properties. In this appendix, I shall point out that the
resolvent of the Dirichlet kernel indeed does have an
extremal property: it coincides exactly (up to complex
conjugation in one variable) with the reproducing kernel of
a certain (interesting) Hilbert space of entire
functions. This could be a new observation, obviously
closely related to the method of Widom \cite{widom}.

The space $mPW_s$ we shall use is, as a set, the
Paley-Wiener space $PW_s$, but the norm is different:
\begin{equation}\begin{split}
  mPW_s &= \{f(z) \text{ entire of exponential type at most
  } s\textrm{ with }\|f\|<\infty\}\\
\|f\|^2 &= \int_{\RR\setminus(-1,1)} |f(t)|^2\,dt
\end{split}
\end{equation}
Let
  $X_s(z,w)$ be the element of $mPW_s$ which is the evaluator
  at $z$: $\forall f\in mPW_s\; (f,X_s(z,\cdot))=f(z)$. We
  shall compare $X_s(z,w)$ with the resolvent of the kernel
\begin{equation}
  D_s(x,y) = \frac{\sin(s(x-y))}{\pi(x-y)}
\end{equation}
on $L^2(-1,1;dx)$.

Let $f\in mPW_s$. It belongs to $PW_s$ so 
\begin{equation}\label{eq:260}
  f(z) = \int_\RR f(t) \frac{e^{is (t-z)} -
  e^{-is(t-z)}}{2\pi i (t-z)}\,dt = \int_\RR
  f(t) \frac{\sin(s(t-z))}{\pi(t-z)}\,dt
\end{equation}
On the other hand:
\begin{equation}\label{eq:261}
  f(z) = \left(\int_{-\infty}^{-1} + \int_1^\infty\right)
  f(t)\overline{X_s(z,t)}\,dt 
\end{equation}
As $\overline{f(\overline z)} = \left(\int_{-\infty}^{-1} +
  \int_1^\infty\right)
  \overline{f(t)}\;\overline{X_s(z,t)}\,dt =
  \left(\int_{-\infty}^{-1} + \int_1^\infty\right)
  \overline{f(t)}{X_s(\overline z,t)}\,dt $ one has
  $\overline{X_s(z,t)} = X_s(\overline z,t)$ for
  $t\in\RR$. We have for $y_1$ and $y_2$ real
  \begin{equation}
    X_s(y_1,y_2) =
  \int_{\RR\setminus(-1,1)}
  X_s(y_1,t)\overline{X_s(y_2,t)}\,dt =
  \int_{\RR\setminus(-1,1)} X_s(y_1,t)X_s(y_2,t)\,dt =
  X_s(y_2,y_1)
  \end{equation}
so more generally $X_s(\overline{z_1},z_2)
  = X_s(\overline{z_2},z_1)$. 

We apply \eqref{eq:261} to $f(z) = \frac{\sin(s(z-y))}{\pi(z-y)}$ for some
  $y\in\CC$:
  \begin{equation}
    \label{eq:262}
    \frac{\sin(s(z-y))}{\pi(z-y)} = \int_{\RR\setminus(-1,1)}
  \frac{\sin(s(t-y))}{\pi(t-y)}X_s(\overline z,t)\,dt
  \end{equation}
We apply \eqref{eq:260} to $f(y) = X_s(\overline z,y)$ for
some $z\in\CC$:
\begin{equation}
  \label{eq:263}
  X_s(\overline z,y) = \int_\RR
  X_s(\overline z,t) \frac{\sin(s(t-y))}{\pi(t-y)}\,dt
\end{equation}
Combining we obtain:
\begin{equation}
  \label{eq:264}
  X_s(\overline z,y) -  \frac{\sin(s(z-y))}{\pi(z-y)} =
  \int_{-1}^1 X_s(\overline z,t)\,\frac{\sin(s(t-y))}{\pi(t-y)} \,dt
\end{equation}
Restricting to $y\in(-1,1)$, $z=x\in(-1,1)$, this says
  exactly:
  \begin{equation}
    \label{eq:265}
    X_s(x,y) = R_s(x,y)
  \end{equation}
where $R_s(x,y)$ is the kernel of the resolvent:
  $1+R_s = (1 - D_s)^{-1}$, $R_s - D_s = R_s D_s$. The
  resolvent $R_s(x,y)$ is entire in $(x,y)$ and the general
  formula is thus:
  \begin{equation}
    \label{eq:fin}
    \forall z,w\in \CC\qquad R_s(z,w) = X_s(\overline z,w)\;.
  \end{equation}

\end{document}